\documentclass[leqno,12pt,a4paper]{article}
\usepackage{setspace} \usepackage{etoolbox} \usepackage{amsfonts,amssymb}
\usepackage[margin=1in]{geometry}
\usepackage{microtype}
\usepackage{times}
\usepackage[amsthm,thmmarks]{ntheorem}
\usepackage{amsmath,xypic,mathtools}
\usepackage{multicol}
\usepackage{hyphenat}
\usepackage[usenames,dvipsnames]{color}
\usepackage{rotating}
\usepackage{longtable}
\usepackage{caption}
\usepackage{tcolorbox}
\usepackage{stackengine}

\usepackage{epigraph}
\usepackage{enumitem}
\setlist[enumerate]{listparindent=0.5in}
\usepackage[cal=boondoxo,calscaled=.96,scr=rsfs]{mathalpha}
\DeclareMathAlphabet{\mathscrbf}{OMS}{mdugm}{b}{n} \newcommand{\be}{\begin{equation}}
\newcommand{\ee}{\end{equation}}
\newcommand{\bes}{\begin{equation*}}
\newcommand{\ees}{\end{equation*}}
\newcommand{\bea}{\begin{eqnarray}}
\newcommand{\eea}{\end{eqnarray}}
\newcommand{\beas}{\begin{eqnarray}}
\newcommand{\eeas}{\end{eqnarray}}
\newcommand{\ben}{\begin{note}}
\newcommand{\een}{\end{note}}
\newcommand{\bexl}{\vskip0.1em\noindent\hrulefill\vskip1em\begin{ExerciseList}}
\newcommand{\eexl}{\end{ExerciseList}\hrulefill}

\newcommand{\bthm}{\begin{theorem}}
\newcommand{\ethm}{\end{theorem}}
\newcommand{\bpro}{\begin{prop}}
\newcommand{\epro}{\end{prop}}
\newcommand{\bcor}{\begin{corollary}}
\newcommand{\ecor}{\end{corollary}}
\newcommand{\bcon}{\begin{conjecture}}
\newcommand{\econ}{\end{conjecture}}
\newcommand{\bp}{\begin{proof}}
\newcommand{\ep}{\end{proof}}
\newcommand{\blem}{\begin{lemma}}
\newcommand{\elem}{\end{lemma}}
\newcommand{\bn}{\begin{note}}
\newcommand{\en}{\end{note}}
\newcommand{\benum}{\begin{enumerate}}
\newcommand{\eenum}{\end{enumerate}}
\newcommand{\bed}{\begin{defn}}
\newcommand{\eed}{\end{defn}}
\newcommand{\brem}{\begin{remark}}
\newcommand{\erem}{\end{remark}}

\newcommand{\btik}{\begin{tikzpicture}\begin{axis}[scale=0.5,axis y line=center, axis x line=middle]}
\newcommand{\etik}{\end{axis}\end{tikzpicture}}

\let\into=\hookrightarrow
\let\mapsto=\longmapsto

\newcommand{\upperRomannumeral}[1]{\uppercase\expandafter{\romannumeral#1}}

 \usepackage{tikz}
\usetikzlibrary{mindmap,backgrounds,cd}
\usepackage{graphicx}
\usepackage[pagewise]{lineno}
\usepackage{stmaryrd}
\usepackage{url}
\usepackage{natbib}
\let\cite=\citep
	\usepackage{fancyhdr}
	\usepackage[colorlinks,citecolor=blue,pagebackref]{hyperref}
	\usepackage{cleveref}
\usepackage{titlesec}
\usepackage{xr}
\newcommand{\constrone}[2]{\cite[{#1\ref{I-#2}}]{joshi-teich}}
\newcommand{\constrtwo}[2]{\cite[{#1\ref{II-#2}}]{joshi-teich-estimates}}

\usepackage{tocloft}
\advance\cftsecnumwidth 0.5em\relax
\advance\cftsubsecindent 0.5em\relax
\advance\cftsubsecnumwidth 0.5em\relax

\newtoggle{draft}
\toggletrue{draft}
\newtheorem{theorem}[equation]{Theorem}      \newtheorem{theoremdef}[equation]{Theorem-Definition}
\newtheorem{lemma}[equation]{Lemma}          \newtheorem{corollary}[equation]{Corollary}  \newtheorem{proposition}[equation]{Proposition}

\theoremstyle{definition}
\newtheorem{conj}[equation]{Conjecture}
\newtheorem{example}[equation]{Example}

\theoremstyle{definition}
\newtheorem{defn}[equation]{Definition}
\theoremstyle{remark}

\theoremsymbol{\ensuremath{\bullet}}

\theoremstyle{definition}
\newtheorem{remark}[equation]{Remark}
\theoremsymbol{\ensuremath{\bullet}}

\newcommand{\bthmdef}{\begin{theoremdef}}
\newcommand{\ethmdef}{\end{theoremdef}}

\numberwithin{equation}{subsection}
\newcommand{\parat}[1]{\subsection{#1}}

\newcommand{\subparat}[1]{\subsection{#1}}
\newcommand{\subpara}{\parat{}}

\let\into=\hookrightarrow
\let\isom=\simeq

\let\tensor=\otimes
\newcommand{\A}{\mathcal{A}}
\newcommand{\abs}[1]{\left\vert#1\right\vert}

\newcommand{\bF}{{\bar{F}}}
\newcommand{\bQ}{{\bar{\Q}}}

\newcommand{\C}{{\mathbb C}}

\newcommand{\End}{\rm{End}}

\newcommand{\Ext}{{\rm Ext}\,}
\newcommand{\F}{{\mathbb F}}

\newcommand{\gal}{{\rm Gal}}

\newcommand{\N}{\mathcal{N}}

\newcommand{\Q}{{\mathbb Q}}
\newcommand{\R}{{\mathbb R}}

\newcommand{\Spec}{{\rm Spec}}

\newcommand{\Z}{{\mathbb Z}}

\renewcommand{\O}{{\mathcal O}}
\renewcommand{\P}{{\mathbb P}}
\renewcommand{\wp}{{\mathfrak p}}

\newcommand{\fm}{{\mathfrak{m}}}

\newcommand{\invlim}{\varprojlim}
\newcommand{\mapright}[1]{{\xymatrix{{}\ar[r]^{#1}&{}}}}

\renewcommand{\bpro}{\begin{proposition}}
	\renewcommand{\epro}{\end{proposition}}
\newcommand{\bdefn}{\begin{defn}}
	\newcommand{\edefn}{\end{defn}}
\renewcommand{\bcon}{\begin{conj}}
	\renewcommand{\econ}{\end{conj}}

\titleformat{\subsection}[runin]{\normalfont\bfseries}{\S\ \thesubsection}{.5em}{}[{\ \ }]
\titlespacing{\subsection}{0pt}{1.5ex plus .1ex minus .2ex}{0pt}
\titleformat{\subsubsection}[runin]{\normalfont\bfseries}{\S\ \thesubsubsection}{.5em}{}[{\ \ }]
\titlespacing{\subsubsection}{0pt}{1.5ex plus .1ex minus .2ex}{0pt}

\title{Construction of Arithmetic Teichmuller Spaces II$\frac{1}{2}$: Deformations of {N}umber {F}ields
}
\author{Kirti Joshi}

\begin{document}
\maketitle

\lhead{}

\epigraphwidth0.55\textwidth
\epigraph{I learned very early the difference between knowing the name of something and knowing something.}{\citeauthor{feynman-quote}}
\newcommand{\act}{\curvearrowright}
\newcommand{\lmp}{{\Pi\act\Ot}}
\newcommand{\lmpi}{{\lmp}_{\int}}
\newcommand{\lmpf}{\lmp_F}
\newcommand{\Om}{\O^{\times\mu}}
\newcommand{\Omf}{\O^{\times\mu}_{\bF}}
\renewcommand{\N}{\mathbb{N}}
\newcommand{\yoga}{Yoga}
\newcommand{\gl}[1]{{\rm GL}(#1)}
\newcommand{\bK}{\overline{K}}
\newcommand{\reptrip}{\rho:G_K\to\gl V}
\newcommand{\reptripp}[1]{\rho\circ\alpha:G_{\ifstrempty{#1}{K}{{#1}}}\to\gl V}
\newcommand{\benumlab}{\begin{enumerate}[label={{\bf(\arabic{*})}}]}
\newcommand{\ord}{\mathop{\rm ord}\nolimits}	
\newcommand{\kcs}{K^\circledast}
\newcommand{\lcs}{L^\circledast}
\renewcommand{\A}{\mathbb{A}}
\newcommand{\bfq}{\bar{\mathbb{F}}_q}
\newcommand{\tripod}{\P^1-\{0,1728,\infty\}}

\newcommand{\vseq}[2]{{#1}_1,\ldots,{#1}_{#2}}
\newcommand{\anab}[4]{\left({#1},\{#3 \}\right)\anabelmap\left({#2},\{#4 \}\right)}

\newcommand{\cpt}{\C_p^\flat}

\newcommand{\gln}{{\rm GL}_n}
\newcommand{\glo}[1]{{\rm GL}_1(#1)}
\newcommand{\glt}[1]{{\rm GL_2}(#1)}

\newcommand{\iut}{\cite{mochizuki-iut1, mochizuki-iut2, mochizuki-iut3,mochizuki-iut4}}
\newcommand{\topics}{\cite{mochizuki-topics1,mochizuki-topics2,mochizuki-topics3}}

\newcommand{\linv}{\mathfrak{L}}
\newcommand{\bedef}{\begin{defn}}
\newcommand{\eedef}{\end{defn}}
\renewcommand{\act}[1][]{\overset{#1}{\curvearrowright}}
\newcommand{\bfx}{\overline{F(X)}}
\newcommand{\anabelmap}{\leftrightsquigarrow}
\newcommand{\ban}[1][G]{\mathcal{B}({#1})}
\newcommand{\pit}{\Pi^{temp}}
 
 \newcommand{\bL}{\overline{L}}
 \newcommand{\bkm}{\bK_M}
 \newcommand{\vbk}{v_{\bK}}
 \newcommand{\vbkm}{v_{\bkm}}
\newcommand{\ocs}{\O^\circledast}
\newcommand{\ot}{\O^\triangleright}
\newcommand{\ocsk}{\ocs_K}
\newcommand{\otk}{\ot_K}
\newcommand{\ok}{\O_K}
\newcommand{\oko}{\O_K^1}
\newcommand{\oks}{\ok^*}
\newcommand{\Qpb}{\overline{\Q}_p}
\newcommand{\Qpbh}{\widehat{\overline{\Q}}_p}
\newcommand{\tr}{\triangleright}
\newcommand{\ocpt}{\O_{\C_p}^\tr}
\newcommand{\ocpf}{\O_{\C_p}^\flat}
\newcommand{\sG}{\mathscr{G}}
\newcommand{\sY}{\mathscr{Y}}
\newcommand{\sxqp}{\mathscr{X}_{\cpt,\Q_p}}
\newcommand{\syqp}{\mathscr{Y}_{\cpt,\Q_p}}
\newcommand{\sxfe}{\mathscr{X}_{F,E}}
\newcommand{\sxfep}{\mathscr{X}_{F,E'}}
\newcommand{\syfe}{\mathscr{Y}_{F,E}}
\newcommand{\syfep}{\mathscr{Y}_{F,E'}}
\newcommand{\loglt}{\log_{\sG}}
\newcommand{\fc}{\mathfrak{t}}
\newcommand{\ku}{K_u}
\newcommand{\kup}{\ku'}
\newcommand{\kt}{\tilde{K}}
\newcommand{\sGpf}{\sG(\O_K)^{pf}}
\newcommand{\hgm}{\widehat{\mathbb{G}}_m}
\newcommand{\hgmp}{\widehat{\mathbb{G}}_{m;p}}
\newcommand{\bE}{\overline{E}}

\newcommand{\bPi}{\overline{\Pi}}
\newcommand{\bPit}{\bPi^{\rm{\scriptscriptstyle temp}}}
\newcommand{\Pit}{\Pi^{\rm{\scriptscriptstyle temp}}}
\renewcommand{\pit}[1]{\Pi^{\scriptscriptstyle temp}_{#1}}
\newcommand{\pitk}[2]{\Pi^{\scriptscriptstyle temp}_{#1;#2}}
\newcommand{\pio}[1]{\pi_1({#1})}

\newcommand{\xan}{X^{an}}
\newcommand{\yan}{Y^{an}}

\newcommand{\sM}{\mathscr{M}}

\togglefalse{draft}
\newcommand{\FF}{\cite{fargues-fontaine}}
\iftoggle{draft}{\pagewiselinenumbers}{\relax}

\newcommand{\onto}{\twoheadrightarrow}

\begin{abstract}
This paper lays the foundation of the \textit{Theory of Arithmetic Teichmuller Spaces of Number Fields} by  explicitly constructing many arithmetically inequivalent avatars  of a fixed number field. This paper also constructs a topological space of  such avatars and  describes its symmetries. Notably, among these symmetries is a global Frobenius morphism which changes the avatar of the number field! The existence of such avatars  has been asserted (and used) by Shinichi Mochizuki in his work on the arithmetic Vojta and Szpiro conjectures (however the existence of distinct avatars has proved difficult for me to ascertain by his methods). On the other hand, I prove that there is a metrisable space of such avatars  and so one can  quantify the difference between two inequivalent avatars and hence  the aforementioned difficulties disappear in my approach. This renders my theory fundamentally and quantitatively more precise than  Mochizuki's Theory (and my results can also be applied to Mochizuki's Theory). With a view towards global applications of my theory (e.g.  the Arithmetic Szpiro Inequality),  in the Appendix (\S~\ref{se:appendix}), I provide a discussion of the proofs of the Geometric Szpiro Inequality due \cite{bogomolov00} and \cite{zhang01} from the point of view of this paper and  prove the geometric case of Mochizuki's Corollary 3.12.
\end{abstract}

\newcommand{\vlnon}{\mathbb{V}_L^{non}}
\newcommand{\vlarch}{\mathbb{V}_L^{arc}}

\newcommand{\fA}{\mathfrak{A}}

\newcommand{\V}{\mathbb{V}}
\newcommand{\vl}{\V_L}
\newcommand{\lvbh}{\widehat{\overline{L}}_v}
\newcommand{\lvbht}{{\widehat{\overline{L}}_v}^\flat}
\newcommand{\lvbhmax}{\widehat{\overline{L}}_v^{max}}
\newcommand{\lvbhtmax}{{\widehat{\overline{L}}_v}^{max\flat}}
\newcommand{\lvbhreal}{\widehat{\overline{L}}_v^{\R}}
\newcommand{\lvbhtreal}{{\widehat{\overline{L}}_v}^{\R\flat}}

\newcommand{\sX}{\mathscr{X}}
\newcommand{\syflv}{\sY_{\lvbht,L_v}}
\newcommand{\syflvmax}{\sY_{\lvbhtmax,L_v}}
\newcommand{\sxflv}{\sX_{\lvbht,L_v}}
\newcommand{\sxflvmax}{\sX_{\lvbhtmax,L_v}}
\newcommand{\bsY}{\mathscrbf{Y}}
\newcommand{\bsX}{\mathscrbf{X}}
\newcommand{\yadl}{\bsY_L}
\newcommand{\yadlmax}{\bsY_L^{max}}
\newcommand{\yadq}{\bsY_{\Q}}
\newcommand{\yadqmax}{\bsY_{\Q}^{max}}
\newcommand{\xadl}{\bsX_L}
\newcommand{\xadq}{\bsX_{\Q}}
\newcommand{\xadlmax}{\bsX_L^{max}}
\newcommand{\xadqmax}{\bsX_{\Q}^{max}}
\newcommand{\yadlpoint}{\{(L_\wp\into K_{\wp}, K_\wp\isom \cpt)\}_{\wp\in\vlnon}}
\newcommand{\xadlpoint}{\{(L_\wp\into K_{\wp}, K_\wp\isom \cpt)\}_{\wp\in\vlnon}}

\newcommand{\bv}{\bar{v}}
\newcommand{\sGal}[1]{\mathbf{G}_{#1}}

\newcommand{\arith}[1]{\mathfrak{arith}(#1)}
\newcommand{\arithl}{\arith{L}}
\newcommand{\adel}[1]{\mathfrak{adel}(#1)}
\newcommand{\adell}{\adel{L}}
\newcommand{\by}{{\bf y}}

\newcommand{\fet}{\mathcal{F\hat{e}t}}
\newcommand{\sI}{\mathcal{I}}
\newcommand{\sA}{\mathcal{A}}

\newcommand{\logp}{\log^+}
\newcommand{\sFrob}{\mathcal{F\!r\!o\!b}}
\newcommand{\sF}{\mathscr{F}}
\newcommand{\pow}[2]{#1\llbracket#2\rrbracket}

\newcommand{\ainf}{A_{\textrm{inf}}}

\setcounter{tocdepth}{3}
\tableofcontents

\newcommand{\ells}{{\ell^*}}
\newcommand{\glqp}{{\rm GL}_2^+(\Q)}
\newcommand{\tslq}{\widetilde{{\rm SL}_2(\Q)}}
\newcommand{\tslql}{{\widetilde{{\rm SL}_2(\Q)}}^\ells}
\newcommand{\slq}{{\rm SL_2}(\Q)}
\newcommand{\tglqp}{\widetilde{{\rm GL}_2^+(\Q)}}
\newcommand{\tglqpl}{{\widetilde{{\rm GL}_2^+(\Q)}}^\ells}
\newcommand{\slr}{{\rm SL}_2(\R)}
\newcommand{\tslr}{\widetilde{{\rm SL}_2(\R)}}
\newcommand{\slz}{{\rm SL}_2(\Z)}
\newcommand{\tslz}{\widetilde{{\rm SL}_2(\Z)}}
\newcommand{\fH}{\mathfrak{H}}
\newcommand{\syi}{\sY_{\infty}}
\newcommand{\syis}{\sY_{\infty,s}}
\newcommand{\syils}{\sY_{\infty,s}^\ells}
\newcommand{\sxi}{\sX_\infty}
\newcommand{\vphi}{\varphi}
\newcommand{\flog}{\mathfrak{log}}
\newcommand{\vnl}{\V_L^{non}}
\newcommand{\Bl}{\mathbb{B}_L}
\newcommand{\Blp}{\mathbb{B}^+_L}
\newcommand{\val}{\V_L^{arc}}
\newcommand{\ssep}{\S\,}
\newcommand{\bvphi}{\boldsymbol{\varphi}}

\nocite{joshi-teich-summary-comments}\nocite{joshi-teich-quest}
\section{Introduction}
\subsection{The Main Results} For the past hundred and fifty years, since the work of E.~Kummer and R.~Dedekind, the analogy  (see the Table of Analogies \ssep\ref{sss:classical-analogies}) between Riemann surfaces and Number Fields (see \cite{dedekind-weber}) has played a central role in the evolution of ideas and theorems in Number Theory. This paper lays the foundation of what may be called  an Arithmetic Teichmuller Theory of Number Fields by constructing explicitly many inequivalent avatars of a fixed Number Field (each avatar of a number field may be thought of as providing the quasi-conformal equivalence class of the number field).  This theory is reminiscent of Teichmuller's Theory of Riemann surfaces, and further deepens this  fecund analogy between Riemann surfaces and Number Fields. This theory is the zero dimensional special case of the General Theory of Adelic Arithmetic Teichmuller Spaces considered in \cite{joshi-teich,joshi-untilts}. 

In this zero dimensional special case, one glues together, into a continuous family, the many arithmetically inequivalent avatars of a number field, just as the classical Teichmuller Theory glues together, into a continuous family, the many (quasi-conformally) inequivalent avatars of  Riemann surfaces. [The two tables \ssep\ref{sss:classical-analogies}{\bf(1, 2)} provide a quick reference for the classical analogies between Riemann surfaces and number fields and new the constructions of this paper.] The isomorphism class of the chosen number field, of course, remains unchanged, but each of its avatars, termed  \textit{arithmeticoids} in \ssep\ref{se:arithmeticoid-adeloid-frobenioid}, represents a quantifiably distinct (and even topologically distinguishable) version of its field arithmetic--i.e. one deforms the very arithmetic of the number field while keeping its isomorphism class fixed!  

The topological space of such avatars (\ssep\ref{se:adelic-ff-curves}) is a connected metrisable space (\Cref{th:distance-bet-arith}). In particular, it makes perfect sense to talk of the distance between two arithmeticoids. This space  comes equipped with many natural symmetries  and notably it is equipped with a (global) Frobenius morphism  (see \Cref{th:galois-action-on-adelic-ff}, \Cref{cor:global-frobenius}) and provides a global Frobenius morphism of a number field (\Cref{def:global-Frobenius-def}).

Moreover, by \Cref{th:fundamental-property-frobenius}, for each non-archimedean prime of the chosen number field, this global Frobenius morphism is in fact a genuine Frobenius morphism in positive characteristic! These global symmetries operate by (quantifiably) changing the field arithmetic structure of the arithmeticoid (\Cref{cor:galois-action-on-arithmeticoids} and also see \Cref{th:galois-action-on-adelic-ff}).

On the other hand, since the isomorphism class of the number field does not change, neither does the topological isomorphism class of its absolute Galois group and hence, by classical Anabelian results of Neukirch-Uchida (also Jarden-Ritter) Theory, neither does the isomorphism class of the multiplicative monoid of non-zero elements of the number field and the multiplicative monoid of its completions at all of its primes. In other words,  the multiplicative structure of the number field, and all its completions, is preserved under such deformations (\Cref{th:anabelomorphic-deforms}). 

\Cref{th:consequence-inequivalent-arithmeticoids} demonstrates that the choice of an arithmeticoid controls local Berkovich analytic structures of any smooth projective variety, at  any arbitrary non-empty subset of primes of the Number Field. But, for any smooth quasi-projective variety, such a choice keeps the isomorphism class of the variety, its \'etale fundamental group and the isomorphism class of the tempered fundamental groups at all primes fixed. This means the anabelian properties of the variety, the number field and all its $v$-adic completions remain unaffected by the choice of the arithmeticoid.

The existence of such avatars of a number field also implies  that the theory of heights  has dependence on the chosen avatar of the number field (\ssep\ref{se:heights}). Hence for Diophantine applications, one may even contemplate averaging over the height functions arising from distinct arithmeticoids. To substantiate my claim that the existence of such avatars should be expected to have Diophantine consequences (such as \cite{mochizuki-iut4}),  I provide a discussion, from my point of view, of the proofs geometric Szpiro inequality due to \cite{bogomolov00}, \cite{zhang01}, exhibiting, in these proofs the existence of such avatars of a Riemann surface and related symmetries--especially the role of a global Frobenius morphism. Notably, in \Cref{th:geometric-case-of-moccor}, I also prove the geometric case  of \cite[Corollary 3.12]{mochizuki-iut3}  in the context of Geometric Szpiro Inequality as established in \cite{bogomolov00} and \cite{zhang01}.

\newtcolorbox[auto counter,crefname={classical analogies}{Analogies Table}]{rosetta}[2][]{colback=white,coltitle=black,colframe=white!25!brown,fonttitle=\bfseries,
	title=Classical Analogies \thetcbcounter:    #2,#1}
\begin{minipage}{\textwidth}
\subsubsection{Table of Analogies}\label{sss:classical-analogies}
The celebrated analogy between  Riemann surfaces and  Number Fields, introduced in  [\citeauthor{dedekind-weber}],  has played a foundational role in Arithmetic Algebraic Geometry ever since. However, as the first table illustrates, there are gaps in this classically well-understood picture. 
\vskip0.20in 
	\numberwithin{table}{subsection}
\begin{rosetta}[label=tab:classical-analogies,grow to right by=1cm,grow to left by=2cm]{The Classical Case (\cite{dedekind-weber})}{}
\centering
		\begin{tabular}{|p{2in}|p{2in}|p{2.5in}|}\hline 
			Riemann Surfaces	& Algebraic Curves & Number Fields \\ 
			\hline 
			Riemann surface $X$	& connected, smooth curve $Y/\C$ & Number Field $L$  \\ 
			\hline 
			Holomorphic mappings $$X'\to X$$ (and Riemann's moduli Theory)	& morphisms of $\C$-schemes $$Y'\to Y$$  & Field homomorphisms $$L\to L'$$ \\ 
			\hline 
			quasi-conformal mappings $$X'\to X$$ (and Teichmuller's Theory)	& $$\mathbf{?}$$ & $$\mathbf{??}$$ \\ 
			\hline 
		\end{tabular} 
\tcblower
$$\small{\text{Dedekind-Weber Analogy}}$$
			$$\small{\text{The field }\C(X)\text{ of meromorphic functions on } X \leftrightsquigarrow \text{The Number Field }L.}$$ 
	\end{rosetta}
\vskip0.25in
\noindent The gap in the above table shown as $\mathbf{?}$ and $\mathbf{??}$ is  filled in by the theory of this paper and is shown in the table below.  [\iut\ asserts that the gap $??$ can be filled in.]
	\vskip0.25in
	\begin{rosetta}[label=tab:classical-analogies2,grow to right by=1cm,grow to left by=2cm]{The theory of the present series papers}{}
\centering
		\begin{tabular}{|p{1.75in}|p{2.45in}|p{2.45in}|}\hline 
			Riemann Surfaces	& Algebraic Curves  $+ \cdots$ & Number Fields $+\cdots$ \\ 
			\hline 
			Riemann surface $X$	& $\text{connected, smooth curve }Y/\C + \text{ Riemann Surface } Y(\C)$ & Arithmeticoid $\arithl$ of a number field $L$  \\ 
			\hline 
			Holomorphic mappings $$X'\to X$$	& morphisms $(Y',Y'(\C))\to (Y,Y(\C))$ are morphisms of $\C$-schemes  $$Y'\to Y$$ &  field homomorphisms $$L\to L'$$ \\ 
			\hline 
			quasi-conformal mappings $$X'\to X$$	& morphisms $(Y',Y'(\C))\to (Y,Y(\C))$ are quasi-conformal mappings $$Y'(\C)\to Y(\C)$$  & morphisms between arithmeticoids $$\arith{L}\to \arith{L'}$$ \\ 
			\hline 
		\end{tabular} 
\end{rosetta}
	\clearpage
\end{minipage}
\newpage
\newcommand{\iutthr}{\cite{mochizuki-iut1,mochizuki-iut2,mochizuki-iut3}}
\subsection{A global Frobenius morphism of a number field}
At this juncture let me say that it is a folklore (lament) in Diophantine Geometry that if (only) one had a Frobenius morphism for a Number Field then many Diophantine conjectures would admit proofs. 

The remarkable idea of Shinichi Mochizuki's \iut\ is that while a Number Field does not have a Frobenius morphism, there are transmutations of a Number Field which behave like a Frobenius morphism. However, Mochizuki attempts to capture the existence of such transmutations through his use of multiplicative monoids of a Number Field, and its $p$-adic completions  (and their algebraic closures) and through his theory of $\flog$-links (as is discussed in \ssep\ref{ss:relation-with-mochizuki-theory}). This makes it  difficult to discern the existence of such transmutations using methods of \iut.  [For more on the relationship between Mochizuki's monoid theoretic approach and my approach via formal groups see \cite{joshi-formal-groups}.]

My discovery  is that there exists, in fact, a global Frobenius morphism (\Cref{def:global-Frobenius-def}):
$$ 
L\mapsto L^{(1)}
$$
which transmutes the given avatar (i.e. an arithmeticoid $\arithl_{\by}$) of a number field $L$ into another avatar $L^{(1)}=\arithl_{\bvphi(\by)}$ of $L$.  This operation arises from the genuine Frobenius morphism (\Cref{th:galois-action-on-adelic-ff}, \Cref{cor:global-frobenius}), $\bvphi:\yadl\to\yadl$ which for each prime $p$, gives the Frobenius morphism, it induces the $p^{th}$-power mapping on the various multiplicative monoids associated to $L$ (\Cref{cor:mult-strs-arith}), and it changes (local and global) arithmetic  and geometry (as asserted by Mochizuki and as established in my work). These properties and \constrtwo{Theorem }{th:flog-links}, \constrtwo{Theorem }{th:flog-kummer-correspondence}, imply that  $L\mapsto L^{(1)}$ gives a rather precise version of  Mochizuki's $\flog$-link (central to \cite{mochizuki-iut3}).

A variant of the global Frobenius morphism of a number field $L$ is the action  of $L^*\act\yadl$ given by \Cref{th:galois-action-on-adelic-ff}. This provides a global action (distinct from the global Frobenius) of $L^*$ on arithmeticoids of a number field which takes an $x\in L^*$ and transforms an arithmeticoid $\arithl_\by$ by the rule $\arithl_\by \mapsto \arithl_{x\cdot\by}$. This action modifies the global arithmetic avatar (i.e. the arithmeticoid) of $L$ by means of the (principal) divisor ${\rm div}(x)$ of the fractional ideal $(x)$ of $L$.  In  \Cref{th:fundamental-property-frobenius}, I show that $x$ acts trivially if and only if $x$ is a root of unity in $L$. [The corresponding result in Mochizuki's theory is \cite[Proposition 3.10(ii)]{mochizuki-iut3}.] This global action of $L^*$ also changes height functions (\ssep \ref{se:heights}).

As I demonstrate (see \ssep\ref{se:appendix}), a (geometric) version of the global Frobenius morphism  is also a central ingredient of the proofs of  the Geometric Szpiro Inequality in \cite{bogomolov00} and \cite{zhang01}. So application to the Arithmetic Szpiro Conjecture is not unexpected at all!

\subsection{The Product Formula as a global Arithmetic Period Mapping}
An important global aspect of arithmetic of number fields is the Product Formula \cite{artin45}. So it should be no surprise that the product formula is integral to the results of this paper. One important point is that each normalized arithmeticoid $\arithl_\by^{nor}$ (for $\by\in\yadl$) comes equipped with its own product formula. In \Cref{th:hyperplane}, I show that each normalized arithmeticoid $\arithl_\by^{nor}$ defines a hyperplane $H_\by$ in the infinite dimensional $\R$-vector space $\oplus_{v\in\vl}\R_v$ (with $\R_v=\R$ for all ${v\in\vl}$). The hyperplane $H_\by$ is given, quite literally, by the logarithm of the product formula for $\arithl_\by^{nor}$ i.e. $H_\by$ is given by the equation $H_\by:\log(L_\by^*)=0$. This hyperplane is stable under multiplication by $x\in L^*$ precisely because of the product formula holds in $\arithl_\by^{nor}$. Hence  passing to the (projective) space of hyperplanes $\P(\oplus_{v\in\vl}\R_v)$ one obtains   a point $[H_\by]\in \P(\oplus_{v\in\vl}\R_v)$. This point moves with action of various symmetries (\Cref{th:galois-action-on-adelic-ff}) on the space of arithmeticoids $\yadl$. This leads to the global arithmetic period mapping (\Cref{th:hyperplane}) $$\yadl\to \P(\oplus_{v\in\vl}\R_v)\qquad \by\mapsto [H_\by].$$ While the next paragraph \ssep\ref{ss:relation-with-mochizuki-theory} details the relationship between my theory with \iut, let me say here that the existence of such a ``moving hyperplanes'' is tacitly claimed in the proof of \cite[Corollary 3.12]{mochizuki-iut3}; on the other hand, their existence has been denied in \cite[Section 2.2]{scholze-stix} (for additional details on this specific point see \cite[\ssep 4.6]{joshi-teich-rosetta}).

\subsection{Relationship with Shinichi Mochizuki's Theory}\label{ss:relation-with-mochizuki-theory} Let me provide a detailed relationship of my theory with Mochizuki's Theory. 
That such an arithmetic Teichmuller Theory of number fields exists, is a remarkable claim of Shinichi Mochizuki in \iut, and in his parlance, this is \textit{``disentangling/dismantling of various structures associated to a number field''}  (see \cite[\ssep I1 and \ssep I4]{mochizuki-iut1}), and to be sure,  his proof of the arithmetic Vojta and Szpiro Conjectures is underpinned on the existence of such avatars of a number field!  [For my proofs of this disentangling/dismantling of arithmetic (and geometric) structures see \ssep\ref{se:multiplicative-structures}, \ssep\ref{se:arith-and-geom}.]

However,  using the methods of \topics, \iut, the existence of such inequivalent avatars of a number field has proved exceedingly difficult to discern for many readers. \textit{Importantly}, it is my understanding that, to date, no one  claiming expertise with \iut\ has offered a convincing proof of precisely how an avatar of a number field changes under the fundamental operations of Mochizuki's theory (this point is central even for the formulation of \cite[Corollary 3.12]{mochizuki-iut3}).  This has led to significant challenges in  recognizing and (and even ascertaining) the non-triviality of the said theory. [Mochizuki's discussion of these concerns is  in \cite{mochizuki-essential-logic}; my discussion is in \cite{joshi-teich-summary-comments,joshi-untilts,joshi-teich-quest}.] 

Second important point which emerges in my analysis of \iutthr, is that the theory of perfectoid fields \textit{inescapably} enters \iutthr\ because of its requirement (see \cite[\ssep I3]{mochizuki-iut1}) of working with arbitrary geometric base-points  in the theory of tempered fundamental groups (my discussion of this, with proofs, is in \cite{joshi-untilts}; the theory of tempered fundamental groups is documented extensively in \cite{andre-book,andre03}, \cite{lepage-thesis}). Once this second point is recognized, \textit{as far as I see}, any \textit{natural} formulation of \iutthr\ (and hence \iut) should \textit{necessarily} proceed in a manner similar to the approach taken in my papers. This recognition also resolves the first mentioned issue:  my work provides a direct, and a far more precise way of arriving at the existence of such avatars and notably, even allows me to construct  the space of such avatars. Importantly, I am able to demonstrate the  precise mechanism by which an avatar changes   and  the claimed  dismantling (and transformation) of the number field structure takes place by means of the symmetries on the space of arithmeticoids (\Cref{cor:galois-action-on-arithmeticoids} and also see \Cref{th:galois-action-on-adelic-ff}). [Mochizuki considers ``Indeterminacies'' as opposed to symmetries. My discussion of this is in \cite{joshi-teich-summary-comments}.] Moreover each arithmeticoid provides a  \textit{Frobenioid}, which is naturally isomorphic to the perfection of the Frobenioid of a number field considered in \iut. One can also obtain realified Frobenioids by the mechanism detailed in \ssep\ref{se:arithmeticoid-adeloid-frobenioid}. Notably, my theory detailed in this paper provides distinguishable arithmetic labels for Frobenioids of number fields (and hence for Hodge Theaters) considered in \iut\ and required in the proof of \cite[Corollary 3.12]{mochizuki-iut3} and consequently my results may be applied  to \iut. [For additional remarks on this see  \ssep\ref{ss:frobenioids}.]

Thirdly, there exists Teichmuller Theoretic proofs, due to \cite{bogomolov00} and \cite{zhang01}, of the Geometric Szpiro Inequality (I provide a discussion of the geometric case from the \textit{optik} of this paper in \ssep\ref{se:appendix}). Notably, I provide a formulation (and proof) of \cite[Corollary 3.12]{mochizuki-iut3}  in the geometric context (my proof of the arithmetic case of this corollary appears in \cite{joshi-teich-rosetta}). As I have argued in \cite{joshi-teich-quest}, in the history of Diophantine Geometry, an understanding of the proof of the geometric case has usually led the way to the arithmetic case. Hence a Teichmuller Theoretic proof of  the arithmetic Szpiro Inequality, asserted in \iut, is not completely unexpected!  The proofs in the geometric case can be viewed as  averaging procedures.  \textit{Averaging over the Teichmuller Theory of Number Fields using a global Frobenius morphism   is indeed the central, remarkable idea in \iut}.
The results in this paper unequivocally establish the existence  of the Teichmuller Theory of Number Fields. This sets the stage for the proof of arithmetic Vojta and Szpiro Conjectures.

\subsection{Outline of this paper}
\textcolor{red}{The table of contents provides a finer outline of the paper.}
In \ssep\ref{se:multiplicative-structures}, I demonstrate the existence of deformations of a perfectoid field keeping its multiplicative structure fixed  (\Cref{thm:mult-structure}). In \ssep\ref{se:adelic-ff-curves}, I introduce adelic Fargues-Fontaine curves (properties are established in \Cref{def:adelic-ff-curves}, \Cref{th:galois-action-on-adelic-ff}). In \ssep\ref{se:arithmeticoid-adeloid-frobenioid}, the key notion of an arithmeticoid of a number field is introduced (\Cref{def:arithmeticoid}, \Cref{th:existence-inequivalent-arithmeticoids}, \Cref{th:existence-inequivalent-arithmeticoids}). This allows one to talk about deformations of a number field (\Cref{def:number-field}) and the constructions of \ssep\ref{se:adelic-ff-curves}  provide a moduli space or rather a Teichmuller Space of deformations of a fixed number field (\Cref{th:parameter-space-of-deformations}). In \ssep\ref{se:arith-and-geom}, I demonstrate how the choice of a deformation of a fixed number field also deforms local analytic geometries at any chosen collection of primes (\Cref{th:consequence-inequivalent-arithmeticoids}). Applying the (tempered) fundamental group functor, and some of Mochizuki's geometric constructions  provides a passage to  Mochizuki's group theoretic view of the geometric objects of my theory.  Arithmetic loops and the arithmetic loop space point of view is discussed in \ssep\ref{ss:arithmetic-loop-spaces} and existence of arithmetic loops is established in \Cref{thdef:arith-loops}, \Cref{cor:arith-loops-at-primes}. 

Let me say that the theory of the present series of papers has applications to other areas of arithmetic geometry. For applications to the Langlands Correspondence see \cite[\ssep 18]{joshi-anabelomorphy}. In \Cref{ss:arithmetic-loop-spaces}, I detail two important applications of the  theory of this paper. The first is to the theory of loop spaces in arithmetic geometry and the second is the
relationship with the  analogy between primes and knots \cite{morishita2002}, \cite{mazur2012} is discussed in \ssep\ref{ss:arithmetic-loop-spaces}. In \Cref{th:distinct-knots}, I establish that each arithmeticoid or avatar of a number field provides a knot associated to each prime and the existence of a parameter space for such knots. Notably, a new insight which my results offer beyond \cite{morishita2002}, \cite{mazur2012}, is that, an arithmetic knot, like its topological counter part,  may not necessarily distinguished by the fundamental group of the knot complement. With a view to applications to \iut, in \ssep\ref{se:arith-and-galois-cohomology}, I discuss Galois cohomology from the point of view of deformations of a number field. In \ssep\ref{se:heights}, the formalism of deformations is applied to study heights of points on a variety over a fixed number field. Notably, heights are now a function of the choice of the deformation or the chosen avatar of the number field. Finally, with a view towards applications of this theory of Arithmetic Teichmuller Spaces to  Mochizuki's proof of the Arithmetic Szpiro Inequality, in \ssep\ref{se:appendix}, I provide my view of the  (classical) Teichmuller Theoretic proofs of the Geometric Szpiro Inequality due to \cite{bogomolov00}, \cite{zhang01} (Mochizuki's view of these proofs is documented in \cite{mochizuki-bogomolov}). In particular, in \Cref{th:geometric-case-of-moccor}, I prove the geometric case of Mochizuki's Corollary 3.12.

\subsection{Acknowledgments}
In the Spring of 2018, I was on a sabbatical visit  to RIMS, Kyoto,visiting Shinichi Mochizuki\footnote{My visit overlapped with the week-long visit of Peter Scholze and Jakob Stix, but I was not part of the conversations between  Hoshi, Mochiuzki, Scholze and Stix.  I received a copy of \cite{scholze-stix} from Fucheng Tan, when that report was put in circulation.}. During my visit, Mochizuki gave me some private lectures on his approach to the Teichmuller theoretic proofs of the Geometric Szpiro Inequality due to \cite{zhang01,bogomolov00} (Mochizuki's exposition followed his own commentary \cite{mochizuki-bogomolov} on these proofs)\footnote{My work recorded in the present series of papers has emerged from my attempts in understanding these three papers (\cite{bogomolov00}, \cite{zhang01} and \cite{mochizuki-bogomolov}) on the geometric Szpiro Inequality and their relationship to \iut--for more on this see \ssep\ref{se:appendix}.}. My own reflections on the relationship between these geometric proofs and Mochizuki's papers is recorded here\footnote{I am an algebraic geometer and an algebraist of the function theory sort, and I found Mochizuki's purely group theoretic approach (taken in \iut) extremely unsatisfactory. I have preferred to develop Arithmetic Teichmuller Theory from a function theoretic vantage point (as is recorded in my papers).} in \ssep\ref{se:appendix}. This paper (and all the papers in this series) grew out of my attempts to independently verify the remarkable claims of \iut.

Unfortunately, Mochizuki stopped corresponding with me after the release of \cite{joshi-untilts-2020}\footnote{My Untilts paper \cite{joshi-untilts-2020} (current version is \cite{joshi-untilts}) provided an independent demonstration of the existence of arithmetic holomorphic structures and signaled the existence of an Arithmetic Teichmuller Theory bearing a close parallel to Teichmuller's Theory of Riemann surfaces and retaining all the anabelian features of \iut.}.  I have made numerous attempts to reach out to Mochizuki (between Nov. 2020-June 2024) but without much success. This is extremely unfortunate as both of us are describing the same theory (my version provides  a more precise way to establish this theory in greater generality  than Mochizuki's  pure group theoretic approach).

At any rate, it is my pleasure to thank Shinichi Mochizuki for conversations  regarding the geometric Szpiro inequality and matters related to \iut\ during 2018-2020. The influence of Mochizuki's ideas on the present series of papers should be clear to all.

I would also like to thanks Peter Scholze for correspondence and some comments.

\section{Preliminaries}\label{se:preliminaries}
\subparat{Anabelomorphisms and amphoricity} 
I will freely use some terminology established in \cite{joshi-anabelomorphy}. I will say that two schemes $X,Y$ (resp. Berkovich analytic spaces $X,Y$ over valued fields) are \textit{anabelomorphic schemes} (resp. \textit{tempered anabelomorphic analytic spaces}) if there exists a topological isomorphism $\pi_1^{et}(X)\isom \pi_1^{et}(Y)$ (resp. $\pi_1^{temp}(X)\isom \pi_1^{temp}(Y)$) between their \'etale (resp. tempered fundamental groups) and I will refer to any topological isomorphism $\pi_1^{et}(X)\isom \pi_1^{et}(Y)$ as an \textit{anabelomorphism} between $X$ and $Y$ (resp. a \textit{tempered anabelomorphism} between $X$ and $Y$). Anabelomorphy of schemes (resp. anabelomorphy of analytic spaces) is an equivalence relation on schemes (resp. on analytic spaces). By \cite{andre03}, for the Berkovich analytic space $\xan$ arising from geometrically connected, smooth quasi-projective variety $X$ over $p$-adic fields, one has a continuous homomorphism $\pi_1^{temp}(\xan)\to \pi_1^{et}(X)$ and the latter is the profinite completion of the image of this homomorphism.  Hence for such analytic spaces, the existence of a tempered anabelomorphism implies the existence of an anabelomorphism.

The second notion from \cite{joshi-anabelomorphy} which I recall here is the following. A quantity $Q_X$ (resp. a property $P$ or an algebraic structure $S_X$) associated to a scheme $X$ (resp. to an analytic space $X$) is said to be \textit{amphoric} if whenever $X,Y$ are anabelomorphic (resp. tempered anabelomorphic) one has $Q_X=Q_Y$,  both $X,Y$ have property $P$ or neither one does, and one has an isomorphism $S_X\isom S_Y$ of algebraic structures. In other words, the quantity $Q_X$, the possession of property $P$ by $X$, and the isomorphism class of $S_X$ are determined by the anabelomorphism (resp. tempered anabelomorphism) class of $X$.

For non-trivial examples of anabelomorphisms and associated amphoric objects see \cite{joshi-anabelomorphy,joshi-gconj}.

\subparat{Some basic assumptions} \label{pa:restrict-num-fild} Let $L$ be a number field, $\V_L$ be the set of inequivalent non-trivial valuations of $L$, let $\vlnon$ be the set of inequivalent non-archimedean valuations of $L$, let $\vlarch$ be the set of inequivalent archimedean valuations of $L$. \textit{I will assume for simplicity that $L$ has no real embeddings. For instance $L\ni \sqrt{-1}$ is adequate (but not equivalent to the former condition).}

\parat{Valued fields and Archimedean Algebraically closed perfectoid fields}\label{pa:valued-fields}
To put the theory of archimedean valued fields on par with $p$-adic valued fields, I introduce the notion of archimedean algebraically closed perfectoid fields and corresponding Fargues-Fontaine curves. For this purpose it is important to work with valued fields $(K,\abs{-}_K)$  in the sense of \cite[Chapter 6, \ssep 6.1 Corollary 2]{bourbaki-alg-com} i.e. $K$ is a field equipped with a function $$\abs{-}_K:K\to \R^{\geq0}$$ satisfying
\benumlab
\item $\abs{x}_K=0$ if and only if $x=0$,
\item $\abs{x\cdot y}_K=\abs{x}_K\cdot \abs{y}_K$ for all $x,y\in K$;
\item there exists a real number $A>0$ (possibly depending on $\abs{-}_K$) such that 
$$\abs{x+y}_K\leq A\cdot\max(\abs{x}_K,\abs{y}_K)$$ 
for all $x,y\in K$.
\eenum
Observe that $\abs{-}_K$ provides a homomorphism $K^*\to \R^{>0}$ and one says that $\abs{-}_K$ is trivial if the image of $K^*$ in $\R^{>0}$ is $\{1\}$. Throughout this paper, I will only work with non-trivially valued fields.

By \cite[Chapter 1, Theorem 10]{artin-algebraic-book}, one may always take $A=\max(\abs{1}_K,\abs{2}_K)$.
A valued field $(K,\abs{-}_K)$ is said to be a non-archimedean valued field if and only if $A=1$ otherwise $(K,\abs{-}_K)$  is said to be an archimedean valued field. 

Let $\abs{-}_{\C}$ be the usual complex absolute value, then $(\C,\abs{-}_{\C})$ is an archimedean valued field in the above sense on $\C$ (with $A=2$) and the pair $(\C,\abs{-}_\C)$ is a  Banach field i.e. $\abs{-}_\C$ is a norm on $\C$ with respect to which $\C$ is complete. Recall that

\bpro[Ostrowski's Theorem] Any archimedean valuation on $\C$ (in the above sense) is of the form $$\abs{-}_{\C}^s\text{ for    } 0<s\in \R.$$
\epro 

For treating archimedean valued fields  essentially on par with non-archimedean valuations I make the following definition.

\begin{defn}\label{def:archimedean-perfectoid}
An \textit{archimedean algebraically closed perfectoid field} is a valued field $(K,\abs{-}_K)$ which is isomorphic  with $(\C,\abs{-}_{\C}^s)$ for some $s\in \R^{>0}$. 
\end{defn} 

\begin{defn}\label{def:archimedean-perfectoid-tilt}
Let $(K,\abs{-}_K)$ be an algebraically closed, archimedean perfectoid field. 
\benumlab
\item We define the tilt $(K^\flat,\abs{-}_{K^\flat})$ to be $(K^\flat,\abs{-}_{K^\flat})=(K,\abs{-}_K)$. 
\item In particular, one has the  \textit{tilt of complex numbers} $\C$ given by  $(\C^\flat,\abs{-}_{\C^\flat})=(\C,\abs{-}_{\C})$. 
\item An \textit{untilt of $\C^\flat$} is an archimedean, algebraically closed perfectoid field $(K,\abs{-}_K)$ and an isometry $K^\flat\isom \C^\flat$.
\eenum
\end{defn}

By this definition $(\C,\abs{-}_{\C}^s)$ is an untilt of $\C^\flat$ for any $s\in\R^{>0}$.

\parat{Archimedean Fargues-Fontaine curve}\label{ss:archimedean-ff-curve}
One also needs to work with archimedean analog of Fargues-Fontaine curves. Here, there is no natural choice (\cite{fargues-classfield} works with a different choice from the one discussed here). The choice made here is adapted for what one needs in  Diophantine situations.

I will take \be\sY_{\C^\flat,\C}=\R^{>0}.\ee   

One also equips $\sY_{\C^\flat,\C}$ with an action of $\C^*$ given by:
\be\C^*\times \R^{>0}\ni(z,r)\mapsto \abs{z}\cdot r.\ee

The following is now  elementary:

\bpro\label{pro:archimedean-curve} 
The mapping $\theta\ni\R^{>0}\mapsto (\C,\abs{-}_\C^{\theta})$ sets up a bijection between untilts of $\C^\flat$ and the archimedean Fargues-Fontaine curve $\sY_{\C^\flat,\C}$. The curve $\sY_{\C^\flat,\C}$ is equipped with the above described action of $\C^*$. Notably $z\in \C^*$ operates trivially if and only if $z\in S^1\subset \C^*$.
\epro

Now suppose that $v\in\V_L$ is an archimedean prime of the chosen number field $L$. Then $L\into L_v=\lvbh\isom \C$ and  by using the above considerations one can talk of untilts of $\lvbh$ for archimedean primes of $L$. Notably I will write $$\sY_{\lvbh^\flat,L_v}=\sY_{\C^\flat,\C}$$ for the archimedean Fargues-Fontaine curve considered above.  

\newcommand{\sL}{\mathscr{L}^1}
\newcommand{\sV}{\mathscr{V}}
\renewcommand{\Rp}{\R^{>0}}
\subparat{The Gelfand Spectrum and the Archimedean Fargues-Fontaine curve}
To construct the space of functions on $\sY_{\C^\flat,\C}$ from the Banach space point of view, similar to the non-archimedean case discussed in \cite{fargues-fontaine}, in 
this section I will recall some basic facts from archimedean spectral theory from \cite{gelfand-book}. 

Let $\abs{-}_{\C}$ be the usual absolute value. Consider the topological, locally compact abelian groups $(\R,+)$ and $(\Rp,\times)$. Fix a Haar measure,  denoted by $dh$, on each of these groups. Note that these two groups are topologically isomorphic via  $\exp: \R\to \Rp$ and its inverse $\log:\Rp\to \R$. 

Let $\sL(\R)$ (resp. $\sL(\R^{>0})$ be the space of all measurable complex valued functions $f:\R\to \C$ (resp. $f:\R^{>0}\to \C$) which are absolutely integrable on $\R$, and introduce the norm
$$\abs{f}=\int_{\R} \abs{f(h)}_{\C}dh <\infty,$$
and  a similar integral for $\Rp$,
which makes $\sL(\R)$ (resp. $\sL(\R^{>0})$) into a Banach space (see \cite[Chap IV, \ssep20, Theorem 1]{gelfand-book}).

So one has an isomorphism of the Banach rings 
$$\sL(\R)\isom \sL(\Rp).$$

The convolution of two such functions $f,g\in\sL(\R)$ (resp. $\sL(\Rp)$) is given by
$$(f*g)(t)=\int_{\R}f(t-h)g(h)dh$$
and a similar definition for $\Rp$ (this integral exists for almost all $t$).

Convolution satisfies $$\abs{f*g}\leq \abs{f}\abs{g}$$ 
turns the Banach spaces $\sL(\R)$ and $\sL(\Rp)$ into commutative Banach algebras \cite[\ssep 20, 3, Theorem 1]{gelfand-book}. \textit{Note that these rings do not contain a unit element so these are non-unital Banach rings. } 

To deal with the lack of unit in these rings, let  $\sV(\R)$ (resp. $\sV(\Rp)$) be the group ring of $(\R,+)$ (resp. $(\Rp,\times)$). These are rings obtained by formally adjoining a unit to $\sL(\R)$ and  $\sL(\Rp)$ respectively \cite[\ssep 20, 8, Definition 2]{gelfand-book}. 

Consider the Banach algebra $\sV(\Rp)$. Let $\fm\subset \sV(\Rp)$ be a maximal ideal. Then $\fm$ is a closed maximal ideal in the Banach topology on $\sV(\Rp)$ and the quotient $\sV(\Rp)/\fm$ is a Banach ring and by the Gelfand-Mazur Theorem one has an isomorphism of Banach algebras $$\sV(\Rp)/\fm \isom \C$$ where $\C$ is equipped with its standard Banach norm i.e. isomorphic to $(\C,\abs{-}_\C)$. 

Let $\mathfrak{max}(\sV(\Rp))$ be the maximal spectrum of $\sV(\Rp)$. By \cite[Chapter IV]{gelfand-book} the maximal spectrum $\mathfrak{max}(\sV(\Rp))$ can be equipped with a natural topology and the spectrum contains a distinguished maximal ideal $\fm_{\infty}$,  such that $\{\fm_\infty\}$ is closed in the topology and one has an homeomorphism 
$$  \mathfrak{max}(\sV(\Rp)) -\{\fm_\infty\} \isom \mathfrak{max}(\sL(\Rp)) \isom \Rp.$$ [Moreover $\mathfrak{max}(\sV(\Rp))$ is compact and is the one point compactification of $\mathfrak{max}(\sV(\Rp)) -\{\fm_\infty\}$.]

In particular, one observes from this that, one should define
$$B_{\C^\flat,\C}=\sL(\Rp).$$
Then take $\sY_{\C^\flat,\C}$ to be the Gelfand Spectrum:
$$\mathfrak{max}(\sL(\Rp))=\mathfrak{max}(B_{\C^\flat,\C})=\sY_{\C^\flat,\C}.$$
If $\fm_s\in\sL(\Rp)$ is the maximal ideal corresponding to the point $s\in\Rp=\sY_{\C^\flat,\C}$, then I will take the archimedean  semi-valuation $$\abs{-}_{\fm_s}:\sL(\Rp)\to \sL(\Rp)/\fm_s=\C\to \R^{\geq0}$$ given, for   any $f\in \sL(\Rp)$, by $$\abs{f}_{\fm_s}=\abs{f\bmod{\fm_s}}_\C^s.$$

Thus the Banach algebra $\sL(\Rp)$  provides a natural parameter space of archimedean valuations  on $\C$ given by Ostrowski's Theorem.

Now define the Frobenius morphism 
\be\label{eq:complete-arch-ff-frobenius}\vphi_\infty:\sY_{\C^\flat,\C} \to \sY_{\C^\flat,\C}\ee in this context to be the identity morphism on the group $\Rp$ i.e. \be\vphi_\infty=1:\sY_{\C^\flat,\C}\to \sY_{\C^\flat,\C}.\ee 

Then one defines the ``complete'' i.e. compact archimedean Fargues-Fontaine curve to be the topological quotient space  
\be\label{eq:complete-arch-ff-curve}
\sX_{\C^\flat,\C}=\sY_{\C^\flat,\C}/\vphi^\Z=\Rp/\vphi_\infty^\Z=\Rp.\ee

\subparat{Some recollections of  non-archimedean Fargues-Fontaine curves}
Let  $p>0$ be a prime number. Let $v\in\vnl$ be a non-archimedean prime of $L$.  Let $p_v$ be the residue characteristic of $L_v$.  For the basic theory of non-archimedean perfectoid fields, readers are referred to \cite{scholze12-perfectoid-ihes}. In this paper $F$ will  always refer to an algebraically closed perfectoid field of characteristic $p>0$; subscripted versions of $F$ such as $F_v$ etc. will refer to an algebraically closed perfectoid field of characteristic $p_v>0$.   Let $E$ be a $p$-adic field. The basic theory of  Fargues-Fontaine curves  $\syfe$ and $\sxfe$ is detailed in \cite{fargues-fontaine}. Notably one has a quotient of adic spaces
$$\syfe\to\syfe/\vphi^\Z=\sxfe$$
by the powers of the Frobenius morphism $\vphi:\syfe\to \syfe$.

 In the presence of a number field $L$, and a place $v\in\V_L$, one will apply the theory of \cite{fargues-fontaine}, to the data $F=\lvbh^\flat, E=L_v$ and obtain the curve $\syflv$, and its Frobenius morphism $\vphi_v:\syflv\to\syflv$ and its quotient (as an adic space) $$\sxflv=\syflv/\vphi^\Z_v.$$ 
 
Let $\C_p$ be the completion of an algebraic closure of $\Q_p$. Let $\cpt$ be the tilt of $\C_p$ (see \cite{scholze12-perfectoid-ihes}). Let $\fm_{\cpt}\subset \O_{\cpt}\subset \cpt$ be the maximal ideal of $\cpt$. I will also work with the Fargues-Fontaine curves $\syqp, \sxqp$ (see \cite{fargues-fontaine}). 
 
If $v\in\V_L$ is a prime of $L$ such that $v|p$ then $\lvbh\subset \C_p$ is an isometric embedding which is compatible with the respective Galois actions on both the sides. Note that the residue fields of $\lvbh$ may not coincide with that of $\cpt$ in general.

\subparat{Canonical valuations and normalized valuations}
Let me remark that for Diophantine considerations involving heights one needs to work with suitably normalized valuations so that the product formula holds \eqref{eq:prod-formula1}. As is remarked in the axiomatic characterization of the product formula given in \cite{artin45}, one can work with valued fields in the sense of \ssep\ref{pa:valued-fields} and the product formula holds for a number field with suitable normalizations i.e. at the archimedean primes one can work archimedean valued fields (in the sense of \ssep\ref{pa:valued-fields}) instead of normed archimedean fields.

Let $y\in\syflv$ be a closed classical point. Then the data $(L_v\into K_v,K_v^\flat \isom F_v)$ comes equipped with a natural absolute value given by \cite[Proposition 2.2.17]{fargues-fontaine}. This will be referred to as the \textit{canonical valuation on $K_v$} for this datum of $y$.

On the other hand one can also normalize the valuation of $K_v$ by requiring that $\abs{p}_{K_v}=\frac{1}{p^{f_v}}$ where $f_v$ is the degree of the residue field of $L_v/\Q_p$. This will be called the \textit{$L_v$-standard normalized valuation on $K_v$}.

\brem\label{re:normalization1} 
Regardless of the choice of a  normalization,  an important aspect of the existence of Fargues-Fontaine curves \cite{fargues-fontaine}  is that one cannot  normalize valuations of all untilts of $\lvbht$ simultaneously. More precisely the function $$\syflv\to \R$$ given by 
\be\label{eq:p-adic-val-of-p} y\mapsto \abs{p}_{K_y},\ee
is not a constant function on $\syflv$. \textit{This function should be considered as the non-archimedean analog of Beltrami parameter in classical Teichmuller Theory  \cite[Chapter 1, 1.4.1]{imayoshi-book}.}  For readers unfamiliar  with classical Teichmuller Theory, let me remark that the classical Beltrami parameter \cite[Chapter 1, 1.4.2]{imayoshi-book} is a measure of local distortions of the metric data of a Riemann surface. The function \eqref{eq:p-adic-val-of-p} is a measure of  local (i.e. at $p$) distortions of $p$-adic metrics. Mochizuki's discussion of this phenomenon is in \cite[Remark 3.9.3]{mochizuki-iut1}.
\erem

\section{Deformation of a field keeping its multiplicative structure fixed}\label{se:multiplicative-structures} 
Mochizuki's approach in \iut\ may be described as working with deformations of fields while keeping their multiplicative structure fixed \cite[\ssep I4, Page 27]{mochizuki-iut1} and \citep[\ssep I5, Page 34,35]{mochizuki-iut1}. In my attempt to algebraize this idea, I was led in \cite{joshi-formal-groups} to the notion of monoid formal group laws and the existence of \textit{universal monoid formal group laws.} An important observation which emerged from \cite{joshi-formal-groups} is that both Mochizuki's ``deformations of arithmetic keeping the multiplicative structure fixed'' and the theory of perfectoid fields with fixed tilts,  appear unified from the perspective of \cite{joshi-formal-groups}. This recognition was central to my subsequent work. 

\subparat{Multiplicative monoids of perfectoid fields}\label{pa:mult-monoids} The following result is an important observation which emerged during the writing of \cite{joshi-formal-groups}, \cite{joshi-teich}.

\bdefn 
Fix a prime $p$ and for any algebra $R$. I will write $R$ for the multiplicative monoid of $R$ and define a new multiplicative monoid $\widetilde{R}$ defined by 
\be 
\widetilde{R}=\invlim_{x\mapsto x^p} R.
\ee
\edefn

Let $K$ be a non-archimedean perfectoid field with residue characteristic $p$, let $(K,\times)$ be the multiplicative monoid of the field $K$ (I will write the multiplicative monoid $(K,\times)$ simply as $K$). Then one has the multiplicative monoid $\widetilde{K}$ by
\be 
\widetilde{K}=\invlim_{x\mapsto x^p} K.
\ee

For an Archimedean algebraically closed perfectoid field $K$, define
$$\tilde{K}=K^*$$
(This is consistent with the choice of $\vphi=1$ in the archimedean case.)
\textcolor{red}{Note the asymmetry in the non-archimedean and the archimedean case.}

Elements of $\widetilde{K}$ may be thought of as collections $(x_n)_{n\in\Z}$ of sequences of elements $x_n$ of $K$ such that $x_{n+1}^p=x_n$ for all $n\in\Z$. 
The following is a consequence of \cite[Lemma 3.4]{scholze12-perfectoid-ihes}.
\blem\label{le:conseq-of-scholze} 
There is a continuous isomorphism of multiplicative monoids
$$\widetilde{K}\isom K^\flat.$$ In particular, the isomorphism class of the monoid $\widetilde{K}$ is independent of perfectoid field $K$ and depends only on $K^\flat$.
\elem

\subparat{Multiplicative structure of a perfectoid field}
Note that $0\in \widetilde{K}$ and this makes $\widetilde{K}$ less useful in capturing the additive structure of $K$ using the $p$-adic logarithm. So following \cite{fargues-fontaine} one works with the submonoid
$$\widetilde{1+\fm_{\O_K}}\subset \widetilde{K}$$ obtained from the submonoid $1+\fm_{\O_K}\subset \O_K$ of $1$-units of $\O_K$  for capturing the additive structure of $K$ via the multiplicative monoid 
$\widetilde{1+\fm_{\O_K}}$. First of all, the above description of $\widetilde{K}$ (given by \Cref{le:conseq-of-scholze}) has an analog for $\widetilde{1+\fm_{\O_K}}$ due to \cite[Chapitre 4]{fargues-fontaine}: 
\bpro 
The multiplicative monoid $\widetilde{1+\fm_{\O_K}}$ is topologically isomorphic to $1+\fm_{\O_{K^\flat}}$ and hence In particular, its topological isomorphism class depends only on $K^\flat$. 
\epro
To recover the additive structure of the field one can now use $p$-adic logarithms.
One has a continuous surjective homomorphism $$\log:\widetilde{1+\fm_{\O_K}}\to K\text{ given by }(x_n)_{n\in\Z}\mapsto \log(x_0),$$ where $\log$ is the $p$-adic logarithm, and the  kernel of this homomorphism is a one dimensional $\Q_p$-vector space. The above homomorphism is even a surjection of Banach spaces \cite[Proposition 4.5.14, Exemple 4.5.15]{fargues-fontaine}. So one should think of the multiplicative monoid $$\widetilde{1+\fm_{\O_K}}$$ as capturing the multiplicative structure of the field $K$ and providing the additive structure of the field $K$ (from the multiplicative structure) via the $p$-adic logarithm.

The following theorem is an observation of \cite{joshi-teich}:
\bthm\label{thm:mult-structure}
Let $F=\cpt$ be the tilt of $\C_p$. Let $(K,K^\flat\isom \cpt)$ be a characteristic zero untilt of $\cpt$. Then 
\benumlab
\item One has an isomorphism of topological multiplicative monoids 
$$\widetilde{K} \isom \cpt.$$
\item One has an isomorphism of topological multiplicative monoids (this is in fact an isomorphism of $\Q_p$-Banach spaces)
$$\widetilde{1+\fm_{\O_K}}\isom 1+\fm_{\O_{\cpt}}.$$
\item In particular, $\widetilde{K}$ and $\widetilde{1+\fm_{\O_K}}$  is independent of the choice of $K$ (i.e. the multiplicative structure of $K$ is independent of the choice of $K$).
\item Moreover, here exists untilts $(K_1,K_1^\flat\isom \cpt)$, $(K_2,K_2^\flat\isom \cpt)$ such that the topological fields $K_1$ and $K_2$ are not topologically isomorphic but one has an isomorphism of multiplicative structures $\widetilde{K_1}\isom \widetilde{K_2}$.
\item (there is also a variant for arbitrary algebraically closed perfectoid field $F$ of characteristic $p>0$ with a choice of a Lubin-Tate formal group $\sG$ which shall be added later).
\eenum
\ethm
\bp 
The first assertion is \cite[Lemma 3.4]{scholze12-perfectoid-ihes} and the second is \cite{fargues-fontaine}. The third follows from the fact that topological isomorphism class of $\widetilde{K}\isom F$ is independent of $K$ as $F$ is independent of $K$. The fourth is a consequence of the first three and  of my central observations in \cite{joshi-teich} regarding the consequences of the fundamental theorem of \cite{kedlaya18} and \cite{joshi-formal-groups} (where I prove the existence of universal monoid formal group law which allows me to algebraically demonstrate that the multiplicative structure of rings may be fixed but the additive laws can be made to vary).
\ep

\subparat{Continuous families of fields with a fixed multiplicative structure} The following is an important and immediate consequence which provides a precise version of Mochizuki's idea \cite[Page 28]{mochizuki-iut1} of keeping the multiplicative structure of a field fixed and allowing its topological structure to vary. 
\bcor
The Fargues-Fontaine curve $\sY_{\cpt,\Q_p}$ provides a continuous family of algebraically closed perfectoid fields whose multiplicative structure is fixed (in the above sense) while the topological structure of the fields vary. 
\ecor

\subparat{Formal group approach of \cite{joshi-formal-groups} to fixing multiplicative structures}\label{ss:formal-group-approach}  My paper \cite{joshi-formal-groups} played a crucial role in my recognition of the fact that Mochizuki's assertion regarding additive and multiplicative structures can be algebraized and are related to \cite{joshi-formal-groups}. The relationship to Fargues-Fontaine theory emerges from this algebraization.  In the theory of the Fargues-Fontaine curve $\sY_{F,\Q_p}$, the multiplicative monoid $\hgm(\O_F)$ is held fixed. The precise relationship between Mochizuki's idea of keeping multiplicative structures fixed \iut\ and \cite{fargues-fontaine} is given by the existence of universal monoid formal group law, in which a multiplicative (or additive) monoid of the field is held fixed but the field structure varies.  This is made precise in the following consequence of \cite{joshi-formal-groups}.

\bthm\label{th:formal-group-law}
Let $F$ be an algebraically closed, perfectoid field of characteristic $p>0$ (for example $F=\cpt$). Let $\hgm/\Z_p$ be the multiplicative formal group over $\Z_p$. Consider the multiplicative monoid $\hgm(\O_F)=1+\fm_{F}$. Then there  exists a universal one-dimensional $\hgm(\O_F)$-formal group law $\sF^{univ}$ over a certain ring $L_{\hgm(\O_F)}$ and a ring homomorphism $$L_{\hgm(\O_F)}\to W(\O_F)$$ such that the tautological ${\hgm(\O_F)}$-monoid formal group law $$\sF(x,y)=x+y\in \pow {W(\O_F)} {x,y}$$ arises, by pull-back,  via this ring homomorphism, from the universal formal group law  $\sF^{univ}$ over $L_{{\hgm(\O_F)}}$.
\ethm

\bp 
Observe that the tautological formal group law $\sF(x,y)=x+y\in \pow {W(\O_F)} {x,y}$ is a monoid formal group law for the multiplicative monoid $\hgm(\O_F)$ by homomorphism of multiplicative monoids $$\hgm(\O_F)=1+\fm_{\O_F}\to \End_{W(\O_F)}(\sF)$$ given by multiplication by a Teichmuller representative $$a\mapsto [a]\cdot x\in \pow {W(\O_F)} {x}.$$
Now the proof is immediate from the existence of the universal monoid formal group law $\sF^{univ}_M$  for the monoid $M=\hgm(\O_F)$ over a ring $L_M$ constructed in \cite{joshi-formal-groups}. By the universal property of this universal monoid formal group law, the $\hgm(\O_F)$-monoid formal group law $\sF(x,y)$ over $W(\O_F)$ arises from the universal monoid formal group law by means of a homomorphism $L_{\hgm(\O_F)}\to W(\O_F)$ given by the universal monoid formal group law of \cite{joshi-formal-groups}.
\ep

\brem
For $F=\cpt$, ring  $W(\O_F)=W(\O_{\cpt})=\ainf$. The tautological formal group law $\sF$ referred to above is essentially the addition law of Witt vectors (considered as a formal group law).  
\erem

\subparat{Interchangeability of addition and multiplication}
Mochizuki has asserted in \cite[\ssep I4]{mochizuki-iut1} that addition and multiplication are interchangeable in his theory and this plays an important role in \iut.

Let me illustrate how this appears in my theory. 
\bthm\label{th:interchange} Let $F$ be an algebraically closed perfectoid field of characteristic $p>0$ (e.g. $F=\cpt$) and let $\fm_F$ be the maximal ideal of $\O_F$. Then one has an isomorphism of the $\Z_p$-modules
$$(\fm_F,+) \mapright{AH} (1+\fm_F,\times),$$
where $AH$ is the Artin-Hasse exponential function 
$$AH(T)=\exp\left(\sum_{n\geq 0}\frac{T^{p^n}}{p^n}\right)\in\pow{\Z_p}{T}.$$ Especially one has an additive and a multiplicative (but equivalent) descriptions of closed classical points on $\sY_{F,\Q_p}$:
$$\left((\fm_F,+)-\{0\}\right)/\Z_p^* \isom \abs{\sY_{F,\Q_p}}\isom \left((1+\fm_F,\times)-\{1\}\right)/\Z_p^*.$$
\ethm
\bp 
This is a consequence of \cite[Proposition 2.3.10]{fargues-fontaine}, \cite[Example 2.3.11]{fargues-fontaine} and \cite[Example 4.4.7]{fargues-fontaine}. The action of $\Z_p^*$ is the Lubin-Tate action of $\Z_p$ on the left and action of $\Z_p^*$ is the corresponding action given explicitly in \cite[Example  2.3.11]{fargues-fontaine}.
\ep
\brem
One should view this theorem as asserting that the curve $\syqp$ has an additive avatar as well as a multiplicative avatar and one may switch from the additive avatar to the multiplicative avatar and this switching is facilitated by the above theorem. In some proofs in \cite{fargues-fontaine} this switching is quite advantageous.
\erem
\brem
An important way of thinking about \Cref{thm:mult-structure} and \Cref{th:interchange} is that as decoupled from field structure, addition and multiplication are in some sense providing equivalent information but when viewed  in the (topological) field structure the relationship between addition and multiplication is intertwined in such a   complicated way that the fields may not be topologically isomorphic. In \iut, Mochizuki works with $\bQ_p$ and wants to  make a similar assertion--but using the methods of \iut, it is difficult to exhibit the difference between two versions of $\bQ_p$ provided by two different arithmetic holomorphic structures.
\erem

\section{Adelic Fargues-Fontaine curves}\label{se:adelic-ff-curves}
I want to introduce adelic Fargues-Fontaine curves over number fields. 
\subparat{Definition of $\yadl$} The following definition is central to the constructions of this paper:
\begin{defn}\label{def:adelic-ff-curves} Let $L$ be a number field.  For each $v\in\vl$,  let $L_v$ be the completion of $L$ at $v$. Let $\lvbh$ be a fixed algebraic closure of $L_v$ and let $\lvbht$ be its tilt \cite{scholze12-perfectoid-ihes}. Let $\abs{\syflv}$ be the topological space of closed classical points of the Fargues-Fontaine curve $\syflv$. Then the adelic Fargues-Fontaine curves is the product 
$$\yadl=\prod_{v\in\V_L} \abs{\syflv}$$
and
$$\xadl=\prod_{v\in\V_L} \abs{\sxflv}.$$
\end{defn}
Note that the adelic curves are written in boldfaced script letter. The following lemma will be useful later on.
\blem\label{le:yadl-is-metric}
The topological space  $\yadl$ is a metrisable space (hence Hausdorff and paracompact).
\elem
\bp 
The set $\V_L$ is countable. For each archimedean $v$, $\syflv$ is obviously a metric space. For each non-archimedean $v$, each factor $\syflv,\sxflv$ is equipped with the topology given by a metric  \cite[Proposition 2.3.2]{fargues-fontaine} and hence
$\yadl$ and $\xadl$ are equipped with the product topology and are countable products of metric spaces and hence the product topology  can be given by a metric. 
\ep

\brem\label{rem:yadl-variant} 
Let me remark that the following variant of the construction of $\yadl$ is also useful. Let $L$ be a number field.  For each $v\in\vl$,  let $L_v$ be the completion of $L$ at $v$. Let $p_v$ be the residue characteristic of $L_v$. Let $\C_{p_v}$ be the completion of some algebraic closure of $\Q_{p_v}$ and let $\C_{p_v}^\flat$ be the tilt of $\C_{p_v}$ (\cite{scholze12-perfectoid-ihes}). This makes sense for all $v$. 

Let $\abs{\sY_{\C_{p_v}^\flat,L_v}}$ (resp. $\abs{\sX_{\C_{p_v}^\flat,L_v}}$) be the topological space of closed classical points of the Fargues-Fontaine curve $\sY_{\C_{p_v}^\flat,L_v}$ (resp. $\sX_{\C_{p_v}^\flat,L_v}$). Then the following variant of the adelic Fargues-Fontaine curve $\yadl$ will also be considered: 
$$\yadl'=\prod_{v\in\V_L} \sY_{\C_{p_v}^\flat,L_v}$$
and
$$\xadl'=\prod_{v\in\V_L} \abs{\sX_{\C_{p_v}^\flat,L_v}}.$$
The properties of $\yadl'$ and $\xadl'$ are established in a manner entirely analogous to that of $\yadl$ and $\xadl$ respectively.
\erem

\subparat{Galois and $L^*$ actions on $\yadl$ and a global Frobenius morphism on $\yadl$}\label{ss:GalL-defined}
Fix an algebraic closure $\bL$ of $L$ and let $G_L=\gal(\bL/L)$ be the absolute Galois group of $L$ computed using $\bL/L$. Let $v\in\vl$ be a prime of $L$, let $L_v$ be the completion of $L$ at $v$ and $\bL_v$ a choice of an algebraic closure of $L_v$, and $G_v=\gal(\bL_v/L_v)$ be the absolute Galois group of $\bL_v/L_v$. Let 
\be\label{def:prod-gal-group}\sGal L=\prod_{v\in\vl}G_v\ee  be the product of all the local galois groups equipped with the product topology, making $\sGal{L}$ into a topological group. The following proposition, while not needed in this paper,  is a classical result about anabelian properties of absolute Galois groups of number fields:

\bpro\label{pr:GL-amphoric}
The topological group $\sGal{L}$ is $G_L$-amphoric i.e. the isomorphism class of the topological group $\sGal{L}$ is determined by the isomorphism class of $G_L$.
\epro
\bp 
This is immediate from the fact that the set of toplogical groups  $\{G_{L_v}:v\in\V_L\}$ is characterized as the set of maximal closed subgroups of $G_L$ which are of ``Mixed-Characteristic Local Field Type'' modulo the action of $G_L$ by conjugation (see \cite[Proposition 3.5]{hoshi-number}).
\ep

Let $x\in L^*$ be an element of $L^*$, then one has a natural embedding of $L^*\into L_v^*$ given in the natural way. Then one has the following
\bthm\label{th:galois-action-on-adelic-ff}\
\benumlab
\item There is a natural action of $\sGal{L}$ on $\yadl$ and $\xadl$. 
\item There is a global Frobenius morphism $\bvphi:\yadl\to \yadl$, which for each non-archimedean $v\in\V_L$, is given by the Frobenius morphism of $\varphi_v:\syflv\to \syflv$ (\cite[Chapitre 2]{fargues-fontaine} (which changes local arithmetic holomorphic structures in the sense of \cite{joshi-untilts}).
\item There is a natural action of $L^*$  on $\yadl$.
\item For each non-archimedean $v\in\V_L$, the $L^*$ action on each factor $\syflv$ is through powers of the Frobenius morphism of $\varphi:\syflv\to \syflv$ and this action changes local arithmetic holomorphic structures in the sense of \cite{joshi-untilts}. 
\item Notably $L^*$ operates as a global Frobenius morphism on $\yadl$ and moves arithmetic holomorphic structures.
\item Let $\sG_v$ be a one dimensional Lubin-Tate group over $\O_{L_v}$ for each non-archimedean prime $v\in\V_L^{non}$. Then one has a natural action of the group $$\prod_{v\in\V_L^{non}} {\rm Aut}_{\O_{L_v}}(\sG(\O_{F_v}))$$ on $\yadl$ which moves arithmetic holomorphic structures (\cite{joshi-teich}). 
\item This action given by \cite{joshi-teich} is explicitly given by the action of the group considered in {(\bf5)} by its action  on primitive elements of degree one \cite[Chapitre 1,2]{fargues-fontaine} at all $v\in\vlnon$ as follows.  Let $$\sigma=(\sigma_v)\in \prod_{v\in\V_L^{non}} {\rm Aut}_{\O_{L_v}}(\sG(\O_{F_v}))$$ and let$\pi_v$ be a uniformizer in $\O_{L_v}$ then $$\sigma_v([a_v]-\pi_v)=[\sigma_v(a_v)]-\pi_v$$ for any $a_v\in F_v$ with $0<\abs{a_v}_{F_v}<1$.
\item For $L=\Q$, one can take $\sG_p=\hgm/\Z_p$, $\hgm(\O_{\cpt})=1+\fm_{\O_{\cpt}}$ for each prime $p<\infty$ and then the assertion of {\bf(5)} provides an action of 
$$\prod_{p<\infty} {\rm Aut}_{\Z_p}(\hgm(\O_{\cpt}))=\prod_{p<\infty} {\rm Aut}_{\Z_p}\left(\widetilde{1+\fm_{\C_p}}\right)=\prod_{p<\infty} {\rm Aut}_{\Z_p}\left(1+\fm_{\cpt}\right).$$
on $\yadq$.
\eenum
\ethm
\brem\   
\benumlab
\item The action given in \Cref{th:galois-action-on-adelic-ff}({\bf 5}, {\bf 6}) should be considered as being the arithmetic analog of the action of the Virasoro algebra on moduli spaces of curves \cite{beilinson00b}, \cite{frenkel01-book}. 
\item Notably in the point of view espoused by \cite{fargues-fontaine} of \textit{``functions of the coordinate p'' (for a fixed prime number $p$)} the action given in ({\bf 5,6,7}) changes the \textit{``local coordinate p (for any fixed prime number $p$).''}
\item In \cite[Theorem 3.11 and Corollary 3.12]{mochizuki-iut3}, Mochizuki considers the  action of the group $$\prod_p{\rm Aut}_{\Z_p}(\O_{\bQ_p}^{*}/\mu(\bQ_p))$$ (here $\mu(\bQ_p)$ is the group of roots of unity in $\bQ_p$ and the product is over all rational primes $p$) and asserts  that the arithmetic holomorphic structures change (under the action of these automorphisms) in \iut. [Technical term used for this in loc. cit. is application of (Mochizuki's) Indeterminacy Ind2.] But, \textit{personally},  I have found this important claim of \cite{mochizuki-iut3}  difficult to verify. Hence the discovery in \cite{joshi-teich} of the correct version  of Mochizuki's action (provided here in {\Cref{th:galois-action-on-adelic-ff}(\bf 7)}) is important as it provides a key ingredient required in formations of \cite[Theorem 3.11 and Corollary 3.12]{mochizuki-iut3}. 
\item  Important point demonstrated in \cite{joshi-teich} (and here in \Cref{th:consequence-inequivalent-arithmeticoids}), is that such changes of the variable $p$ (given by ({\bf 5,6,7}) also change geometry and topology of algebraic varieties over $p$-adic fields, but keep the (tempered) fundamental group (and hence the \'etale fundamental group) fixed and this signals the existence of an Arithmetic Teichmuller Theory.
\eenum
\erem
\bp 
The group $\sGal{L}$ operates on $\yadl$ and $\xadl$ through the action of its factor $G_v$ on the factor $\syflv$ and $\sxflv$. This proves {\bf(1)}.

The assertions {\bf(2)} is clear from the construction of $\yadl$ and the properties of Fargues-Fontaine curves established in \cite{fargues-fontaine}. 

The assertions {\bf(3), (4)} may be proved together. To define the action of $L^*$ on $\yadl$ it is enough to use the diagonal embedding $x\mapsto (x)_{v\in\V_L}$ of $$ L^*\into \prod_{v\in\V_L} L_v^*.$$ Let $v\in\V_L$ be a non-archimedean prime of $L$. Let me remind the reader that in my notational convention, $\syflv$ is the topological space of closed classical points of the adic curve $\syflv$ (one can also work with the adic curve, but it is not essential to do so at the moment).  Observe that for each non-archimedean $v\in\V_L$,  there is a natural Lubin-Tate action of $L^*_v$ on  $\fm_{\O_{F_v}}$ in which any uniformizer of the ring $\O_{L_v}$ operates by Frobenius on $\fm_{\O_{F_v}}$ and this extends to an $L_v^*$-action on $\fm_{\O_{F_v}}$). This Lubin-Tate action is explicated in \cite[Chap 2]{fargues-fontaine}. 
Now to arrive at the action of $L_v^*$ on $\syflv$, one notes, from \cite[Proposition 2.3.9]{fargues-fontaine}, that $$\syflv=(\fm_{F_v}-\{0\})/\O_{L_v}^*.$$ So the $L_v^*$ action is such that the action of the subgroup $\O_{L_v}^*\subset L^*_v$ is trivial on $\syflv$. 

For archimedean primes this can be explicated using \Cref{pro:archimedean-curve}: Suppose $v\in\V_L$ is an archimedean prime of $L$. Then $L_v\isom \C$ (as an archimedean valued field) and  by the construction given in \Cref{pro:archimedean-curve}, $\sY_{\lvbh^\flat,L_v}$ is equipped with an action of $L_v^*\isom \C^*$ and $L^*$ operates on $\sY_{\lvbh^\flat,L_v}$ via the inclusion $L^*\into \sY_{\lvbh^\flat,L_v}$.
 This proves {\bf(2), (3)} and {\bf(4)}.

An important discovery of \cite{joshi-teich} is  the action of ${\rm Aut}_{\O_{L_v}}(\sG(\O_{F_v}))$ on $\sY_{F,L_v}$. This is explicated in detail in \cite[Theorem 7.29.1]{joshi-teich} which proves {\bf(5,6,7)}.
\ep

The following corollary to the above constructions (\Cref{th:galois-action-on-adelic-ff}) provides a variant of the global Frobenius morphism constructed above. 
\bcor\label{cor:global-frobenius}
For any number field $L$ as above, the topological space $\yadl$ is equipped with a continuous global Frobenius mapping
$$\bvphi:\yadl\to \yadl$$
given, for any $\by\in\yadl$ by $$\bvphi(\by) =(\varphi_{v}(y_v))_{v\in\vl}.$$
\ecor
\bp 
The proof is clear.
\ep

\brem\label{rm:global-frob-remark}\
\benumlab
\item \textit{In particular, $\bvphi$ (and the multiplicative action of $L^*$ on $\yadl$) is a proxy for a global Frobenius morphism on number fields because locally, at each $v$, an algebraic number operates by the powers of the local Frobenius morphism of $\syflv$.}
\item   Both $\bvphi$ and $L^*\act\yadl$ move points on $\syflv$ (for each non-archimedean $v$) and  change arithmetic holomorphic structures by powers of the local Frobenius morphism!
\item As was shown in \cite{joshi-teich-estimates}, moving in a Frobenius orbit on $\syflv$ corresponds to moving in the vertical column of $\mathfrak{log}$-links in \cite{mochizuki-iut3}. This action $L^*\act \yadl$ provides a global way to do this. 
\item As is discussed in the appendix (\ssep\ref{se:appendix}), in the proof of the geometric Szpiro inequality given in \cite{bogomolov00}, \cite{zhang01} there is a similar ``global Frobenius morphism! (see \ssep\ref{app:global-frobenius})'' [This is my central observation, beyond \cite{mochizuki-bogomolov}, regarding  \cite{bogomolov00}, \cite{zhang01}, and \cite{mochizuki-bogomolov} which  shapes my approach to the arithmetic case.] 
\item As pointed out in \cite{joshi-teich}, \Cref{th:galois-action-on-adelic-ff}(5) is arithmetic analog, at all primes of $L$, of the Virasoro action in \cite{beilinson88}, \cite{kontsevich87} and especially in the Geometric Langlands Program of \cite{beilinson00b}, and should be understood as arising from continuous changes of the ``variable $p$'' (at each fixed prime $p$)--enhancing the view of \cite[Chapitre 1]{fargues-fontaine} of ``holomorphic functions of the variable $p$.''
\eenum
\erem

\brem Let $F$ be an algebraically closed perfectoid field of characteristic $p>0$ and let $E,E'$ be two $p$-adic fields with an anabelomorphism $G_E\isom G_{E'}$. In \cite{joshi-gconj} I have shown that the absolute Grothendieck Conjecture fails to hold for Fargues-Fontaine curves i.e. the curves $\sxfe,\sxfep$ are anabelomorphic but not necessarily isomorphic. \textit{This result of \cite{joshi-gconj} has important consequences for the number field setting.} Let $L$ be a number field. The main theorem of \cite{joshi-gconj} implies that the topological isomorphism class of $\sGal{L}$ does not necessarily determine the isomorphism class of the data $\{G_{L_v}\act\sxflv\}_{v\in\V_L}$ uniquely. Equivalently the data $\{(G_{L_v}\act\syflv,\vphi_v)\}_{v\in\V_L}$ is not determined uniquely by the isomorphism class of $\sGal{L}$). Equivalently the data $\sGal{L}\act\yadl$ is not uniquely determined by the topological isomorphism class of $\sGal{L}$. As is observed in \cite{joshi-teich-rosetta}, the $\{(G_{L_v}\act\syflv,\vphi_v)\}_{v\in\V_L}$ (resp. $\{(G_{L_v}\act\syflv,\vphi_v)\}_{v\in\V_L}$ equivalently $\sGal{L}\act\yadl$) can be viewed as the ``prime-strips'' of my theory. Notably one may consider the action of ${\rm Aut}(\sGal{L})$ on $\yadl$.
This is the origin of Mochizuki's first Indeterminacy Ind1 (see \cite{mochizuki-iut3}) from my point of view.
\erem

\subparat{Correspondences on $\yadl$}\label{ss:correspondences}
Let $L$ be a number field. In this section I want to elaborate on \Cref{th:galois-action-on-adelic-ff}{\bf(6)} and establish the existence of some natural correspondences on the adelic Fargues-Fontaine curve $\yadl$.
The theorem is the following. 

\bthm\label{th:adelic-correspondences}
Let $L$ be a number field satisfying \ssep\ref{pa:restrict-num-fild}. Let $\yadl$ be the adelic Fargues-Fontaine curve associated to $L$. Let $$\sigma=(\sigma_p)_{p\in\V_{\Q}}\in \left(\prod_{p<\infty}{\rm Aut}_{\Z_p}((1+\fm_{\O_{\cpt}}))\right)\times {\rm Aut}(\R^{>0}).$$ Let $\by\in\yadl$. Then $\sigma$ provides a divisorial correspondence on $\yadl$ i.e. $\sigma$ provides one-to-many mapping (written as a finite formal sum over points) on $\yadl$:
$$\sigma(\by)=\sum_{j=1}^n(\by_{j})$$
with $n\leq [L:\Q]$.
\ethm
\bp 
Let $\by=(y_v)_{v\in\V_L}$. Let $v\in\V_L$ be a non-archimedean prime lying over some prime number $p$ (i.e. $v|p$). By \cite{fargues-fontaine}, one has a natural finite \'etale morphism $f_v:\syflv\to\syqp$ of degree $[L_v:\Q_p]\leq [L:\Q]$. Let $y_v'=f(y_v)\in \syqp$. Now as was observed in \cite{joshi-teich}, the $p$-component, $\sigma_p$ of $\sigma$ is a $\Z_p$-linear automorphism and provides an action on    $$\syqp=\left((1+\fm_{\O_{\cpt}})-\{1\}\right)/\Z_p^*.$$ Let $y_v''=\sigma_p(y_v')\in \syqp$. Since $f_v:\syflv\to\syqp$ is a finite \'etale morphism, its fiber $f_v^{-1}(y''_v)\subset\syflv$ is finite and let $$f_v^{-1}(y''_v)=\{y_{v,1}'',\ldots,y_{v,n}''\},$$ be the schematic fiber. Then one has $n\leq [L:\Q]$. The archimedean case is dealt with similarly.  This gives the one-to-many mapping $$\by\mapsto f^{-1}((y''_v)_{v\in\V_L})=\{\by_1:=(y_{v,1}'')_{v\in\V_L},\ldots,\by_n:=(y_{v,n}'')_{v\in\V_L}\}.$$
Then one has the correspondence 
$$\by\mapsto\sum_{j=1}^n(\by_j),$$
with $n\leq [L:\Q]$.
This proves the theorem.
\ep

\subparat{A fundamental property of the $L^*$-action} The following result is important despite its proof being elementary.
\bthm\label{th:fundamental-property-frobenius} 
Let $\alpha\in L^*$. Then $\alpha$ acts trivially on $\yadl$ if and only if $\alpha\in L$ is a root of unity.
\ethm
\bp 
First let $v\in\vnl$. If $\alpha$ acts trivially on $\syflv$ then $\alpha\in \O_{L_v}^*$ and hence by hypothesis, $\alpha$ is an $\O_{L_v}$-unit for all $v\in\vnl$. Thus $\alpha\in\O_L^*$. 

Now suppose $v$ is archimedean prime of $L$. As $\alpha$ acts on $\syflv=\R^{>0}$ through multiplication by $\abs{\alpha}_{\lvbh}$ on $\R^{>0}$ (see \Cref{pro:archimedean-curve}), so $\alpha$ acts trivially implies that at all archimedean primes $v$, one has, by \Cref{pro:archimedean-curve}, $$\abs{\alpha}_{\lvbh}=1.$$ 

Hence $\alpha$ is an algebraic integer which has absolute value one at all archimedean primes and thus by an old result of Kronecker (see \cite[Chapter 2, Exercise 11]{marcus-book}), $\alpha$ is a root of unity contained in $L$. This proves the assertion.
\ep

\bcor 
Let $L$ be a number field with no real embeddings and assume that $L$ does not contain any non-trivial roots of unity, then the $L^*$-action on $\yadl$ has trivial stabilizers.
\ecor

\brem 
Mochizuki makes a similar but less precise assertion in \cite[Remark 1.2.3(ii)]{mochizuki-iut3}.
\erem

\brem 
Theorem asserts that the action of $L^*$ on the topological space $\yadl$ is stacky in general.
Using the formalism of topological stacks, one can nevertheless construct the quotient topological stack $[\yadl/L^*]$ which should play the role of the correct $\xadl$ i.e. the correct definition of $\xadl$ should be
$$\xadl=[\yadl/L^*].$$ But I will not pursue this for now.
\erem

\subparat{The global topological rings $\Bl\supset\Blp$} Now let me outline two important corollaries of \Cref{th:galois-action-on-adelic-ff}.
Let $\O_L^\triangleright=\O_L-\{ 0\}\subset L^*$ be the multiplicative monoid of non-zero elements of $\O_L$.  For each $v\in\vnl$, let $\pi_v\in\O_{L_v}$ be a uniformizer. Let $\vphi_v$ be the Frobenius of $\syflv$. Let $B_{L_v}$ denote the ring constructed in \cite[Chapitre 2]{fargues-fontaine}. 

 Let $$B_{L_v}^{\vphi=\pi_v}\subset B_{L_v}$$ be the Banach subspace on which Frobenius operates by multiplication by $\pi_v$.

I will now extend the definition of these rings to archimedean primes of $L$.
For $v\in\val$, let $B_{\C_v^\flat,L}=B_{\C^\flat,\C}$ be the Banach algebra constructed in the archimedean setting and let $$B_{\C_v^\flat,L}^+=B_{\C_v^\flat,L},$$ and one takes $\pi_v=1$ and $\vphi_v=1$. This means $$B_{\C_v^\flat,L}^{\vphi_v=\pi_v}=B_{\C_v^\flat,L} \text{ for all } v\in\val.$$

Hence $B$-rings (and related objects) are defined for all $v\in\V_L$.
\bdefn\label{def:adelic-B-ring}
The adelic $B$-ring for a number field $L$ with no real embeddings is the ring, equipped with the product topology:
$$\Bl=\prod_{v\in\V_L} B_{L_v},$$
and  let $$\Blp=\left\{(x_v)\in \Bl: x_v\in B^+_{L_v} \text{ for all but finitely many }v \right\}.$$
Let $$\Bl(1)=\prod_{v\in\V_L} B_{L_v}^{\vphi=\pi_v}.$$
\edefn

\brem
One should think of $\Bl$ as functions on $\yadl$ and $\Blp$ as functions on $\yadl$ which are bounded at all but finitely many primes of $L$.
\erem

\subparat{Fundamental actions on $\Bl$} The following corollary is now clear from the preceding results and definitions:
\bcor\label{cor:fun-actions} 
Let $L$ be a number field with no real embeddings.
Then
\benumlab
\item there is a natural action of $\sGal{L}$ on $\Bl$, which preserves $\Bl(1)$, and
\item a natural action of $L^*$ on $\Bl$ which corresponds to the $L^*$ action on $\yadl$ constructed above.
\item The natural action of $\O_L^\triangleright\subset L^*$ on $\Bl$ preserves $\Bl(1)$.
\item There is a natural action of $\prod_{v\in\V_L^{non}} {\rm Aut}_{\O_{L_v}}(\sG(\O_{F_v}))$ on $\Bl(1)$. 
\eenum
\ecor

\subparat{Definition of $\yadl^{max}$ and $\xadlmax$} The following ``maximally complete'' variant of $\yadl$ and $\xadl$ are also useful in Diophantine considerations.
\begin{defn}\label{def:adelic-ff-curves-maximally-complete} Let $L$ be a number field.  For each $v\in\vl$,  let $L_v$ be the completion of $L$ at $v$.  Let $\lvbhmax$ be a maximally complete field containing $L_v$ isometrically (such fields exist by \cite{kaplansky42}, \cite{poonen93}) and let $\lvbhtmax$ be the tilt of $\lvbhmax$ (see \cite{scholze12-perfectoid-ihes}). The isomorphism class of $\lvbhmax$ is uniquely determined by the fact that it is maximal among all the immediate extensions of $\bL_v$. But there are still many inequivalent untilts of $\lvbhtmax$ as witnessed by the existence of the Fargues-Fontaine curve $\sY_{\lvbhtmax,L_v}$. 

Let $\abs{\syflvmax}$ (resp. $\abs{\sxflvmax}$) be the topological space of closed classical points of the Fargues-Fontaine curve $\syflvmax$ (resp. $\abs{\sxflvmax}$). Then the adelic , maximally complete Fargues-Fontaine curves are the products of topological spaces 
	$$\yadlmax=\prod_{v\in\V_L} \abs{\syflvmax}$$
	and
	$$\xadlmax=\prod_{v\in\V_L} \abs{\sxflvmax}.$$
\end{defn}

The following property of $\yadl^{max}$ will be useful in understanding how one can use it in Diophantine contexts.
\bpro\label{pro:working-ymax} 
Let $\by=(y_v)_{v\in\vl}\in\yadl^{max}$. Then for each $v\in\vl$
\benumlab
\item  the residue field $K_{y_v}$ of $y_v\in  \abs{\sY_{F_v,L_v}}$ is a maximally completed (hence algebraically closed) perfectoid field with an isometric embedding of the $p$-adic field $L_v$, and 
\item the value group of $K_{y_v}$ is equal to the value group of $\lvbh$. 
\item If $v\in\vlnon$ then the residue field of $K_{y_v}$ is the residue field of $\lvbh$. 
\item In particular, for each $v\in\vl$, $K_{y_v}$ contains the algebraic closure of $\bL_v$ of $L_v$ and also its completion $\lvbh$ (in $K_{y_v}$).
\eenum  \epro
\bp 
This is immediate from \cite{fargues-fontaine} and \cite{scholze12-perfectoid-ihes}.
\ep

\subparat{Definition of ``realified'' curves $\yadl^{\R}$ and $\xadl^{\R}$}\label{ss:realified-curves} The following realified versions of  variant of $\yadl$ and $\xadl$ are also useful in Diophantine considerations. These realified curves are constructed as follows. 

Let $\overline{\F}_p\supset \F_p=\Z/p$ be an algebraic closure of the field $\F_p$ with $p$ elements. One would like to work with algebraically closed perfectoid fields with value group equal to $\R$ and with residue field $\overline{\F}_p$. Let  $$\boldsymbol{C_p^\flat}=\overline{\F}_p((t^\R))$$ of Han-Malcev power series. By \cite{kaplansky42}, \cite{poonen93} this field has the following  properties
\benumlab
\item the residue field of $\boldsymbol{C_p^\flat}$ is  $\overline{\F}_p$, and 
\item the value group of $\boldsymbol{C_p^\flat}$ is $\R$ (hence the valuation is of rank one), 
\item and $\boldsymbol{C_p^\flat}$ is a maximally complete field of characteristic $p>0$.
\eenum  
Moreover  $\boldsymbol{C_p^\flat}$ is also algebraically closed and hence perfect of characteristic $p>0$ and hence perfectoid (by \cite{scholze12-perfectoid-ihes}) and uniquely determined up to an  isomorphism of valued fields, by the three listed properties \cite{poonen93}.

By \cite{fargues-fontaine}, one may construct Fargues-Fontaine curves with the pair $(\boldsymbol{C^\flat_p},E)$ consisting of the field $\boldsymbol{C^\flat_p}$ and any $p$-adic field $E$ as the input datum.

\begin{defn}\label{def:adelic-ff-curves-realified} Let $L$ be a number field with no real embeddings.  For each $v\in\vl$,  let $L_v$ be the completion of $L$ at $v$. If $v\in\vlnon$, let the residue characteristic of $v$ be denoted by $p_v$.  Let us define algebraically closed perfectoid fields $F_v$ as follows:
\be 
F_v	=
\begin{cases}
\boldsymbol{C^\flat_{p_v}} & \text{ if } v\in \vlnon,\\
L_v\isom \C^\flat=\C & \text{ if }v\in\vlarch.
\end{cases}
\ee  
The \textit{adelic, realified  Fargues-Fontaine curves $\yadl^\R$ and $\xadl^\R$}  are the products of topological spaces 
	$$\yadl^\R=\prod_{v\in\V_L} \abs{\sY_{F_v,L_v}}$$
	and
	$$\xadl^\R=\prod_{v\in\V_L} \abs{\sX_{F_v,L_v}}.$$
\end{defn}
The following property of $\yadl^\R$ will be useful in understanding how one can use it in Diophantine contexts.
\bpro\label{pro:working-yR} 
Let $\by=(y_v)_{v\in\vl}\in\yadl^\R$. Then for each $v\in\vl$
\benumlab
\item  the residue field $K_{y_v}$ of $y_v\in  \abs{\sY_{F_v,L_v}}$ is a maximally completed (hence algebraically closed) perfectoid field with an isometric embedding of the $p$-adic field $L_v$, and 
\item the value group of $K_{y_v}$ is equal to $\R$. 
\item If $v\in\vlnon$ then the residue field of $K_{y_v}$ is $\overline{\F}_{p_v}$. 
\item In particular, for each $v\in\vl$, $K_{y_v}$ contains the algebraic closure of $\bL_v$ of $L_v$ and also its completion $\lvbh$ (in $K_{y_v}$).
\eenum  \epro
\bp 
This is immediate from \cite{fargues-fontaine} and \cite{scholze12-perfectoid-ihes}.
\ep

\section{Arithmeticoids, Adeloids, Ideloids and Frobenioids}\label{se:arithmeticoid-adeloid-frobenioid}
I want to introduce the notion of arithmetic deformations of  number fields. Now let $L$ be a number field, $\V_L$ be the set of inequivalent non-trivial valuations of $L$, let $\vlnon$ be the set of inequivalent non-archimedean valuations of $L$, let $\vlarch$ be the set of inequivalent archimedean valutions of $L$. 

\subsection{Definition of an arithmeticoid} Let $v\in\V_L$ be a valuation of $L$. Let $L_v$ be the completion of $L$ at $v$, let $\lvbh$ be the completion of a fixed algebraic closure of $L_v$. Let $\lvbh^\flat$ be its tilt (for $v\in\vlarch$ this is \Cref{def:archimedean-perfectoid-tilt} and for $v\in\vlnon$ this is \cite{scholze12-perfectoid-ihes}). Thus tilting makes sense for every $v\in\V_L$ i.e. for every valuation of $L$.

\begin{defn}\label{def:arithmeticoid}
An \textit{arithmeticoid of $L$} or an \textit{arithmetic deformation of $L$} is a point $\by\in\yadl$. I will write $\arithl_\by$ for the arithmeticoid provided by $\by$. Giving an arithmeticoid $\arithl_\by$ means that for each $v\in\V_L$ one is given an untilt $(L_v\into K_v,K_v^\flat\isom \lvbh^\flat)$. Thus $$\arithl_\by=((L_v\into K_v,K_v^\flat\isom \lvbh^\flat))_{v\in\V_L}\in\yadl.$$ I will refer to $$(L_v\into K_v,K_v^\flat\isom \lvbh^\flat)$$ as  \textit{the local data at $v$ provided by $\arith{L}_{\by}$}. The topological ring (with the product topology) $$\prod_{v\in\V_L} K_v$$ will be called the \textit{arithmetic ring of the arithmeticoid $\arith{L}_{\by}$}. 
If the point $\by$ is not required in a given context, then I will simply write $\arith{L}$ instead of $\arith{L}_\by$. Since an arithmeticoid corresponds to a point $\by\in\yadl$, it makes sense to talk about \textit{ the equality of arithmeticoids}. I will use the term \textit{equivalence  of arithmeticoids} to be synonymous with the term \textit{equality of arithmeticoids}.
\edefn

The following is useful in understanding the definition of equality (or equivalence) of arithmeticoids.

\blem
Two arithmeticoids of $L$ are  equal if and only if for every $v\in\V_L$, one has an isomorphism of untilts $(L_v\into K_v,K_v^\flat\isom \lvbh^\flat)\isom (L_v \into K_v',K_v^{'\flat}\isom \lvbh^\flat)$.
\elem
\bp 
This is clear from the definition.
\ep

The following is elementary:

\blem\label{le:subfield-lemma} 
Let $\arith{L}$ be an arithmeticoid of $L$. Then its arithmetic ring  $\prod_{v\in\V_L} K_v$ contains an algebraic closure $\bL$ of the embedded subfield $L$ and for each $v\in\V_L$ one has an algebraic closure  $\bL_v\subset K_v$ of $L_v$ contained in $K_v$.
\elem
\bp 
Let  $R=\prod_{v\in\V_L} K_v$ be the arithmetic ring of the arithmeticoid $\arith{L}$. Clearly $R$ is of characteristic zero. The embedding $L\into L_v$ (given by the completion of $L$ at $v$) provides an embedding $$L\into R=\prod_{v\in\V_L} K_v$$ given by $x\mapsto (x)_v$. Elementary considerations show that every non-constant polynomial $f(T)\in L[T]\subset R[T]$ has a root in $R$. Thus $R$ contains $\bL$ as an embedded subring.
\ep

I will denote the embedding of $L$ in an arithmeticoid $$\iota_L:L\into \arith{L}$$ and the embedding of $\bL$ by $$\iota_{\bL}:\bL\into \arith{L}.$$
From now on these embeddings are to be implicitly understood--in other words I will write $$\arith{L}$$ instead of $$(\arith{L}, \iota_L:L\into\arith{L}).$$

\brem 
One should think of an arithmeticoid of $\arith{L}$ as providing a copy of the arithmetic of $L$ and $\bL$.
\erem

\bthm\label{th:distance-bet-arith} Let $L$ be a number field with no real embeddings. Then $\yadl$ is a metrisable  space.
\ethm
\bp
Since $\yadl$ is a metrisable space by \Cref{le:yadl-is-metric} and hence for a choice of standard metric $d$ on a countable product of metric spaces, one defines the distance between two arithmeticoids $\arith{L}_{\by_1}$ and $\arith{L}_{\by_2}$ by using the chosen metric on $\yadl$ i.e. one may take
$$d(\arith{L}_{\by_1},\arith{L}_{\by_2})=d(\by_1,\by_2).$$ This completes the proof.
\ep

\brem
In particular, one can quantify quite naturally, how the arithmetic of $L$ presented by $\arith{L}_{\by_1}$ is different from that presented by $\arith{L}_{\by_2}$. Ultimately the idea of \iut\ is that one wants to average over different ways of doing arithmetic (and consequently geometry by \cite{joshi-teich}) itself!
\erem

\bdefn Let $R_\by$ be the arithmetic ring of an arithmeticoid $\arith{L}_\by$. Then the \textit{tilt of the arithmeticoid $\arith{L}_\by$} is the topological ring given by
$$R^\flat_\by=\prod_{v\in\V_L} K_v^\flat.$$
\edefn

\textcolor{red}{Note that:} 
\benumlab
\item The nomenclature `Tilt of an arithmeticoid' should not be confused with tilting as defined in \cite{scholze12-perfectoid-ihes}, but the nomenclature is natural and hence I use it.
\item The tilt of an arithmeticoid is a topological ring of characteristic zero because  it has factors which are of positive characteristic $p>0$ for  any rational primes $p$.
\eenum

\bpro\label{pr:isom-tilts-of-arithmeticoids} 
Let $L$ be a number field. Then the topological isomorphism class of the tilt of any arithmeticoid $\arith{L}_\by$ is independent of the arithmeticoid $\by$.
\epro
\bp 
Indeed, the arithmetic ring $R_\by$ of $\arith{L}_\by$ has tilt $R_\by^\flat$ and using the tilting data $(K^\flat\isom \lvbht)_{v\in\V_L}$ of the arithmeticoid one obtains an isomorphism
$$R^\flat_\by\isom \prod_{v\in\V_L} \lvbht.$$
Since the topological ring on the right is independent of the choice of $\by$, the assertion follows.
\ep

\subsection{Topological equivalence of arithmeticoids} Now let me introduce the notion of topological equivalence of arithmeticoids. This will help understand the relevance of \Cref{pr:isom-tilts-of-arithmeticoids}.
\bdefn I will say that two arithmeticoids $\arithl_{\by_1}, \arithl_{\by_2}$ are \textit{topologically equivalent} if and only if their arithmetic rings are  topologically isomorphic. 
\edefn

\brem\
\benumlab
\item Evidently two topologically inequivalent arithmeticoids are  inequivalent. 
\item \textcolor{red}{Let me remark that topological equivalence of arithmeticoids ignores the tilting information completely and therefore is a weak notion from the point of view of this paper. Secondly,  it is also a difficult notion to work with.} 
\item \Cref{pr:isom-tilts-of-arithmeticoids} asserts that the tilts of arithmeticoids are always topologically isomorphic.
\item On the other hand, in \Cref{th:existence-inequivalent-arithmeticoids} I show that topologically inequivalent arithmeticoids of $L$ also exist. 
\eenum
\erem

\subsection{Normalizing valuations}\label{ss:normalization}
\textit{The fundamental global constraint in the arithmetic of a number field is the validity of the product formula  for $L$} which asserts that for the standard choice of normalizations for valuations $(L_v,\abs{-}_v)$ (given in \cite{artin45}),  the global product formula holds for $L$ i.e. for all $x\in L^*$ one has:
\be\label{eq-prod-form} 
\prod_{v\in\vl}\abs{x}_{v}=1,
\ee
or equivalently
\be\label{eq-log-prod-form}
\sum_{v\in\vl}\log\abs{x}_{v}=0.
\ee
Now suppose one is given a $\by=(\xi_v)_{v\in\vl}\in\yadl$, then one has the associated collection of algebraically closed, complete valued fields $(K_{\xi_v})_{v\in\vl}$. Since the two valuations $\abs{-}_v)$ and $\abs{-}_{K_{\xi_v}}$ on $L_v$ are equivalent, there exist, for each $v\in\vl$, real numbers $\alpha_v\in\R$ such that
\be
(L_v,\abs{-}_v)=(L_v,\abs{-}_{K_{\xi_v}}^{\alpha_v}).
\ee

Hence the product formula \eqref{eq-prod-form} reads, for all $x\in L^*$ as
\be\label{eq:prod-formula1} \prod_{v\in\vl}\abs{x}_{K_{\xi_v}}^{\alpha_v}=1
\ee

Hence one says that 
\be\alpha_\by=(\alpha_v)_{v\in\vl}\in \prod_{v\in\vl} \R_{>0}\ee is the \textit{normalization coordinate of $\by$}.
When valuations induced on $\{L_v \}_{v\in\vl}$ by the local data of $\arith{L}_{\by}$ are normalized so that \eqref{eq:prod-formula1} holds, I will write this as $\arith{L}_{\by}^{nor}$ and call it a normalized arithmeticoid of $L$ corresponding to $\by$.

Now suppose that $\yadl\ni\by'=(\xi_v')_{v\in\vl}\neq (\xi_v)_{v\in\vl}=\by\in\yadl$ then one obtains the normalization coordinates $\alpha_{\by'}=({\alpha_v'})_{v\in\vl}$ of $\by'$ i.e. for each $v\in\vl$, real numbers $\alpha_v'\in\R$ such that
\be
(L_v,\abs{-}_v)=(L_v,\abs{-}_{K_{\xi_v'}}^{\alpha_v'}).
\ee
Because the absolute values $(\abs{-}_{K_{\xi_v}})_{v\in\vl}$ depend on $\by\in\yadl$, one has
\be \alpha_{\by'}=(\alpha_v')_{v\in\vl}\neq (\alpha_v)_{v\in\vl}=\alpha_{\by}\ee
in general and hence one says that  the \textit{product formula normalization  of $\by$ moves with $\by\in\yadl$}. Especially, one sees that there is no uniform choice of the normalization coordinate $\alpha_{\by}$ for all $\by\in\yadl$.

\newcommand{\adeloid}{\mathcal{A\!d\!e\!l\!o\!i\!d}}
\newcommand{\ideloid}{\mathcal{I\!d\!e\!l\!o\!i\!d}}
\subparat{Adeloids and ideloids provided by arithmeticoids}
\begin{defn}
	Let $\arith{L}_{\by}$ be an  arithmeticoid of $L$. Then the \textit{adeloid of $L$} associated to $\arith{L}_{\by}$ is the subring
	$$\adeloid_{\arith{L}}=\left\{(x_v)_{v\in \V} \in \prod_{v\in\V_L} K_v: x_v\in\O_{K_v} \text{ for all but finitely many }v\in\V \right\},$$
	and the ideloid $\ideloid_{\arith{L}}$ of associated to $\arith{L}_{\by}$  is the subgroup given by
	$$\ideloid_{\arith{L}}=\left\{(x_v)_{v\in \V}\in \prod_{v\in\V_L} K_v^*: x_v\in \O^*_{K_v} \text{ for all but finitely many }v\in\V \right\}.$$
\end{defn}	

\blem\label{le:mult-act-on-ideloid}
The adeloid $\adeloid_{\arith{L}}$ is a topological ring and the ideloid $\ideloid_{\arith{L}}$ is a topological group. The ideloid $\ideloid_{\arith{L}}$ is equipped with a natural multiplicative action of $L^*$:
$$L^*\times \ideloid_{\arith{L}} \to  \ideloid_{\arith{L}}$$ which is given by
$$ (x,(z_v)_{v\in\V_L})\mapsto (x\cdot z_v)_{v\in\V_L}.$$
\elem
\bp 
This is clear from the definitions.
\ep

\bdefn\label{def:arith-deg-ideloid}
Let $\arith{L}_\by$ be an arithmeticoid, let $\alpha_{\by}=(\alpha_v)_{v\in\vl}$ be its normalization coordinate; let $\ideloid_L$ be the ideloid of $\arith{L}$. Then one has a natural homomorphism of groups
$$\deg_{\by}:\ideloid_{\arith{L}}\to \R$$
given by 
$$\deg_{\by}((x_v)_{v\in\V_L})= \sum_{v\in\V_L} \log\abs{x_v}_{K_v}^{\alpha_v}.$$
This is a  finite sum by the definition of an ideloid. This homomorphism shall be referred to as \textit{(normalized) arithmetic degree homomorphism of an ideloid}.
\edefn

The following property of the normalized arithmetic degree is now immediate from \Cref{def:arith-deg-ideloid} and \eqref{eq:prod-formula1}:
\blem 
In the notation of \Cref{def:arith-deg-ideloid}, the normalized arithmetic degree has the following property:
for all $x\in L^*$ and for all $(z_v)_{v\in\vl}\in \ideloid_{\arith{L}}$ one has
$$\deg_{\by}(x\cdot (z_v)_{v\in\V_L}) = \deg_{\by}( (z_v)_{v\in\V_L}).$$
In particular the normalized arithmetic degree descends to a homomorphism
$$ \deg_{\by}:\frac{\ideloid_{\arith{L}}}{L^*}\to \R.$$
\elem
\bp 
This is immediate from the product formula \eqref{eq:prod-formula1} and \Cref{def:arith-deg-ideloid}.
\ep

\subsection{Basic properties of arithmeticoids} 
\bpro 
Let $\arith{L}$ be an arithmeticoid. Then the $\arith{L}$ provides an adeloid $\sA_{\arith{L}}$ as well as an ideloid $\sI_{\arith{L}}$ of $L$.  Let $L'\supset L$ be a  finite extension of $L$. Then 
\benumlab
\item the adeloid $\adeloid_{\arith{L}}$ provides the ring   $\mathfrak{adeles}({L'})\subset \adeloid_{\arith{L}}$ of adeles of $L'$, and
\item the ideloid  $\ideloid_{\arith{L}}$ provides  the group of ideles $\mathfrak{ideles}({L'})\subset \ideloid_{\arith{L}}$ of $L'$.
\item The rings of adeles  of $L$ provided by two arithmeticoids may not be comparable as objects embedded in their respective arithmeticoids.
\eenum
\epro
\bp 
This is clear from the respective definitions.
\ep

The following  will be used in the rest of the paper:
\bthm\label{th:existence-inequivalent-arithmeticoids}
Let $L$ be a number field with no real embeddings.  
\benumlab
\item Each arithmeticoid $\arith{L}_\by$ of $L$ provides  a  Galois category, denoted  $\fet(L)_\by$, of finite extensions of $L$ contained in its arithmetic ring $R_\by=\prod_{v\in\V_L} K_v$, and
\item  hence provides an isomorph $G_{L;\by}$ of the absolute Galois group $G_L$ of $L$.
\item Each arithmeticoid provides, for each $v\in\V_L$,  a  Galois category, denoted  $\fet(L_v)_{\by_{v}}$, of finite extensions of $L_v$ contained in $K_v$ and,
\item  hence a preferred isomorph $G_{L_v;\by_v}=G_{L_v;K_v}$ of the absolute Galois group $G_{L_v}$ of $L_v$, and
\item hence a preferred isomorph of the product group ${\sGal{L}}_{;\by}$.
\item The class of arithmeticoids of $L$ is non-empty; and
\item topologically inequivalent arithmeticoids of $L$ also exist.
\eenum
\ethm
\bp Let $\arith{L}$ be an arithmeticoid of $L$. By \Cref{le:subfield-lemma}, the arithmeticoid provides an algebraic closure $\bL$ of $L$, hence it also provides the category of finite extensions contained in $\bL$ and hence a it also provides a preferred isomorph $G_L=\gal(\bL/L)$. 

Similarly, for each $v$, the local data of $\arith{L}$ at $v$ provides an algebraically closed perfectoid field $K_v$ and hence a preferred algebraic closure of $L_v$ namely the algebraic closure  $\bL_v\subset K_v$ of $L_v\subset K_v$. Hence one obtains the group $\gal(\bL_v/L_v)$ for each $v$. Hence an isomorph $\sGal{L}$. 

To prove that the class of arithmeticoids of $L$ is non-empty one proceeds as follows. For each prime $v$ of $L$ choose an algebraic closure $\bL_v$ of $L_v$ and
let $\lvbh$ be the completion of $\bL_v$ as $\lvbh$ is an algebraically closed, perfectoid field let $\lvbh\isom \lvbht$ constructed by \cite{scholze12-perfectoid-ihes}. Then for each $v\in\V_L$, one has an  untilt $(L_v\into \lvbh,\lvbh^{\flat}=\lvbh^{\flat})$. Write $y_v\in \syflv$ for the point of the Fargues-Fontaine curve corresponding to this choice of an untilt at $v$. Then write  $\by=(y_v)_{v\in\V_L}$. This provides an  arithmeticoid of $L$. 

One could of course have arrived at this by choosing (by axiom of choice) any point $\by\in\yadl$ and take the arithmeticoid $\arith{L}_{\by}$ which this choice provides.

The last assertion is more subtle. This comes from an observation made in \cite{joshi-teich}. As was pointed out in \cite{joshi-teich}, by \cite{kedlaya18}, there exists, for any $v\in \vnl$, a pair of closed classical points $y_{1,v},y_{2,v}\in \syflv$ such that the residue fields $K_{1,v}$ and $K_{2,v}$ are not topologically isomorphic. Choose  a $v\in\vnl$ and $\by_1,\by_2\in\yadl$ such that $\by_{1,v}=y_{1,v},\by_{2,v}=y_{2,v}$. Then the arithmeticoids $\by_1,\by_2$ are not topologically inequivalent by construction (and hence not equivalent).
\ep

\subsection{Multiplicative structures provided by an arithmeticoid}
The notation of \ssep\ref{pa:mult-monoids} will be in force. The following important result is an immediate consequence of the above discussion and \Cref{thm:mult-structure}.
\bthm\label{th:mult-strs-arith} 
Let $L$ be a number field and let $\arithl_{\by_1}, \arithl_{\by_2}$ be two arithmeticoids of $L$. 
Then for each valuation $v\in\V_L$ one has an isomorphism of the respective multiplicative structures of the algebraically closed perfectoid fields at $v$ in the sense of \Cref{thm:mult-structure}. In particular, one has an isomorphism of their respective topological multiplicative monoids (defined in \Cref{pa:mult-monoids}):
$$\prod_{v\in\V_L} \tilde{K}_{1,v} \isom \prod_{v\in\V_L}\tilde{K}_{2,v}.$$
On the other hand the arithmeticoids $\arithl_{\by_1}$ and $\arithl_{\by_2}$ need not be equal or even be topologically equivalent.
\ethm
\bp 
This is clear from \Cref{le:conseq-of-scholze}, \Cref{thm:mult-structure} and \Cref{th:existence-inequivalent-arithmeticoids}.
\ep

\bcor\label{cor:mult-strs-arith} 
In the notation of \Cref{th:mult-strs-arith}, suppose additionally that $\by_1=\by$ and $\by_2=\bvphi(\by)$. For each 
$v\in\vlnon$, let $p_v$ be the residue characteristic of $v$, and for $v\in\vlarch$, let $p_v=1$. Then the topological isomorphism of \Cref{th:mult-strs-arith} is given explicitly by:
$$\prod_{v\in\V_L} \tilde{K}_{1,v} \isom \prod_{v\in\V_L}\tilde{K}_{2,v}$$
given by $$(\tilde{x}_v)_{v\in\vl}\mapsto (\tilde{x}_v^{p_v})_{v\in\vl},$$
i.e. $\bvphi$ is the $p_v^{th}$-power operation on the local multiplicative monoids corresponding to $\by$ and $\bvphi(\by)$.
\ecor
\bp 
This is clear from \Cref{th:mult-strs-arith}, \Cref{th:galois-action-on-adelic-ff} and  \Cref{cor:global-frobenius}.
\ep

\subsection{Three fundamental actions} 
 The following is a corollary of \Cref{th:galois-action-on-adelic-ff}
\bcor\label{cor:galois-action-on-arithmeticoids}
There is a  natural action of each of the groups, $L^*$, $\sGal{L}$ and $\prod_{v\in\V_L^{non}} {\rm Aut}(\sG(\O_{F_v}))$, on  the collection of arithmeticoids $$\left\{\arith{L}_{\by}:\by\in\yadl \right\}$$ of $L$ given by the action of these groups on $\yadl$ given by \Cref{th:galois-action-on-adelic-ff}.
\ecor

\subsection{Deformations  of a number field}
\bdefn\label{def:number-field}
	A \textit{arithmetic deformation of a number field $L$} is the choice of an arithmeticoid $\arith{L}_\by$ for some point $\by\in\yadl$. This comes equipped with
	\benumlab
	\item an embedding $\iota_L: L \into \arith{L}_\by$, and 
	\item an embedding $\iota_{\bL}: \bL \into \arith{L}_\by$, and 
	\item for every finite extension $L'$ of $L$, an adeloid $\sA_{\arith{L'}_{\by}}$, and 
	\item an ideloid $\sI_{\arith{L'}_\by}$. 
	\item Each $\arith{L}$ is also equipped with an intrinsic copy of the absolute Galois group $G_L$ and the group $\sGal{L}$.
	\eenum
	Using the normalization coordinate $\alpha_{\by}=(\alpha_v)_{v\in\vl}$ of $\by$,
	valuations on the local data of an arithmeticoid $\arith{L}$ can be normalized so that following global arithmetic requirement, namely the validity of the product formula:  for all  $x\in L^*$ one has
	\be\label{eq:prod-formula2} 
	 \prod_{v\in\V_L}\abs{x}_{K_v}^{\alpha_v}=1.
	\ee
\edefn

\brem
Thanks to the last assertion of \Cref{th:consequence-inequivalent-arithmeticoids} one can say that arithmeticoids of $L$ provide isomorphs of a number field $L$ and its algebraic closure $\bL$ which may even be topologically distinguished from each other (in general). 
\erem

\subsection{The global Frobenius morphism of a number field}
One of the central consequences of \Cref{cor:global-frobenius} is that a fixed number field $L$ comes equipped with a global Frobenius morphism. This has been asserted in \iut\ and \cite{mochizuki-frobenioid1} but my construction, while conforming to Mochizuki's monoidal philosophy espoused in his construction of Frobenioids in \cite{mochizuki-frobenioid1}, is quantitatively  far more precise.

Let $L$ be a number field (assumed to have no real embeddings for simplicity) and let $\arithl_\by$ be an arithmeticoid of $L$ with $\by\in\yadl$ (resp. $\by\in\yadlmax$, $\by\in\yadl^\R$). 
\bdefn\label{def:global-Frobenius-def}
Let $\bvphi:\yadl\to\yadl$ be the global Frobenius morphism \Cref{cor:global-frobenius}. Let $\arithl_{\bvphi(\by)}$ be the arithmeticoid of $L$ corresponding to $\bvphi(\by)$. Then 
$$\arithl_\by \mapsto \arithl_{\bvphi(\by)}$$
is called the \textit{global Frobenius morphism of the number field $L$}. Using the notations prevalent in the context of the Frobenius morphism of characteristic $p>0$  (for example in \cite{Illusie1979a}), let me write 
$$ L_{\by}:=\arithl_\by $$
and 
$$L_{\by}^{(1)}:=L_{\bvphi(\by)}=\arithl_{\bvphi(\by)}=:\arithl_{\by}^{(1)}.$$
Then the global Frobenius morphism of the number field is the operation
$$L_{\by}\mapsto L^{(1)}_{\by} \text{ (equivalently } \arithl_\by\mapsto \arithl_{\by}^{(1)}).$$
If $\by$ is clear from the context then simply write this as  $$L \mapsto L^{(1)}.$$ The mapping  $L\mapsto L^{(1)}$ the global Frobenius morphism of the number field $L$. 
\edefn
\brem 
The properties of $L\mapsto L^{(1)}$ are easily established using \Cref{th:galois-action-on-adelic-ff}, \Cref{cor:global-frobenius}, \Cref{cor:mult-strs-arith} and \Cref{rm:global-frob-remark}. Notably, by \Cref{cor:global-frobenius}, the relationship between the local multiplicative structures  corresponding to $L$ and $L^{(1)}$ is through the $p_v^{th}$-power mapping for each $v\in\vl$.
\erem

\subsection{The product formula as a global arithmetic period mapping}\label{ss:period-map-prod-form}
\newcommand{\sH}{\mathscr{H}}
\newcommand{\grV}{\mathcal{G\!r}(\sV)}
An important aspect of the theory of classical Teichmuller Spaces is the existence of period mappings. In this section I  describe a natural period mapping in the theory of Arithmetic Teichmuller Spaces of number fields. 
\bthm\label{th:hyperplane} 
Let $L$ be a number field, fixed as above. Consider the vector space $$\sV_L=\bigoplus_{v\in\vl} \R_v$$
where each $\R_v=\R$ is given its usual absolute value. Then 
\benumlab
\item Given an arithmeticoid $\arithl_\by$ of $L$ with $\by=(y_v)_{v\in\vl}$ with arithmetic ring $R_\by=\oplus_{v\in\vl}K_{y_v}$ (\Cref{def:arithmeticoid}) and normalization coordinate $\alpha_{\by}=(\alpha_v)_{v\in\vl}$ then the ideloid $\ideloid_{\arithl_\by}$ provides a homomorphism of groups
$$\log_{\arithl_\by}: \ideloid_{\arithl_\by}  \to \sV_{L}$$ given by
$$ \ideloid_{\arith{L}_\by} \ni (x_v)_{v\in \vl} \mapsto  (\log\abs{x_v}_{K_{y_v}}^{\alpha_v})_{v\in\V_L}.$$
I will simply write  $\log_{\by}$ instead of $\log_{\arithl_\by}$ if there is no possibility for confusion.
\item Each normalized arithmeticoid $\arithl_\by$ provides a degree homomorphism (\Cref{def:arith-deg-ideloid})
$$ \begin{tikzcd}\ideloid_{\arith{L}_\by} \arrow[r,"\log_{\by}"]\arrow[rr, bend left, "\deg_{\by}"] &  \sV_L \arrow[r] & \R\end{tikzcd}$$
given by  the formula
$$ (x_v)_{v\in \vl}\mapsto (\log\abs{x_v}_{K_{y_v}}^{\alpha_v})_{v\in\V_L}\mapsto  \sum_{v\in\vl} \log(\abs{x_v}_{K_{y_v}}^{\alpha_v})=\deg_{\by}((x_v)_{v\in\vl}).$$
\item As a consequence, each normalized arithmeticoid $\arith{L}_{\by}^{nor}$ (\cref{ss:normalization}) provides a hyperplane $H_{\by}\subset \sV_{L}$ given by the logarithm of product formula equation \eqref{eq:prod-formula1} for  $L\into\arithl_\by^{nor}$:
$$H_{\by}:\deg_{\by}((x)_{v\in\V_L})=\sum_{v\in\V_L}\log\abs{x}_{K_{y_v}}^{\alpha_v}=0\qquad (\text{for all }x\in L^*).$$ One can think of $H_{\by}$ as being given by the equation
$$H_{\by}:\log_{\by}(L^*)=0.$$
\item Explicitly, the product formula hyperplane $H_\by\subset \bigoplus_{v\in\vl}\R$ is given by the equation
$$ 
H_\by=\left\{ z=(z_v)\in \bigoplus_{v\in\vl} \R: \sum_{v\in\vl} \alpha_v\cdot z_v=0 \right\}.
$$
[The sum in the formula  is well-defined because any $z=(z_v)_{v\in\vl}\in \bigoplus_{v\in\vl}\R$ has all but finitely many coordinates $z_v=0$.]
\item As the product formula holds for the normalized arithmeticoid $\arithl_\by^{nor}$, the hyperplane $H_\by$ is stable under multiplicative action of  $L^*\act\ideloid_{\arith{L}}$ given in \Cref{le:mult-act-on-ideloid}.
\item On the other hand, the action $L^*\act\yadl$ constructed in \Cref{th:galois-action-on-adelic-ff} moves $H_\by$ to $H_{x\cdot\by}$ for any $x\in L^*$.
\item Let $\P(\sV_L)$ be the projective space of hyperplanes in $\sV_{L}$.  Then one has a natural and non-constant mapping:
$$\begin{tikzcd} \yadl\arrow[r, "\by\mapsto H_{\by}"] &\P(\sV_L),\end{tikzcd}$$
called the period mapping of the  Arithmetic Teichmuller Theory of the number field $L$ which associates to an arithmeticoid $\arith{L}_{\by}$ the hyperplane $H_{\by}$.
\item Let $\P(\yadl)\subset \P(\sV_L)$ be the image of $\yadl$ under this function. Then $\P(\yadl)$ is equipped with a natural action of all the groups and symmetries (for example \Cref{cor:galois-action-on-arithmeticoids} and \Cref{th:galois-action-on-adelic-ff}) which act up on $\yadl$.
\eenum
\ethm
\bp 
This is clear from the properties of $\yadl$ and of arithmeticoids established so far. The fact that mapping $\by\to H_\by$ is non-constant is immediate from the fact that normalizations of arithmeticoids do not carry over because locally at each prime $v\in\vlnon$, valuations cannot be simultaneously normalized on $\syqp$.
\ep
\brem\label{rem:hyperplane}\ 
\benumlab
\item The mapping $\yadl\to \P(\sV)$ should be considered as the fundamental period mapping of Arithmetic Teichmuller Theory of Number Fields. 
\item Mochizuki's discussion of the role of product formulas is in \cite[Remark 3.9.6]{mochizuki-iut3}, \cite[Remark 5.2.1]{mochizuki-iut1}, \cite[Page 71]{mochizuki-gaussian}.  In the context of  the proof of  \cite[Corollary 3.12]{mochizuki-iut3} (where in the above proposition, one takes $L$ to be $L_{mod}$) given in \cite{joshi-teich-rosetta}, one works with tuples of arithmeticoid $(\by_1,\ldots,\by_{\ells})\in\Sigma_L\subset \yadl^\ells$ and this gives a tuple of hyperplane $$(H_{\by_1},\ldots,H_{\by_\ells})$$ arising from the $\ells$-tuple of ideloids $$\left(\ideloid_{\arith{L}_{\by_1}},\ldots,\ideloid_{\arith{L}_{\by_\ells}}\right).$$ The above proposition remains valid even for this case of $\ells$-tuples of hyperplanes. 
\item The relationship between the coordinates of $(H_{\by_1},\ldots,H_{\by_\ells})$ is given by the scaling relationship \cite[Theorem 4.2.2.1]{joshi-teich-rosetta}.
\item It is worth noting that in the context of \cite[Corollary 3.12]{mochizuki-iut3} Mochizuki considers (without a transparent proof)  a version of this construction by means of the theory of realified Frobenioid $\sFrob(L)^\R$ of $L$ (see my discussion of this Frobenioid is in \ssep\ref{ss:realified-frobenioids}). Specifically Mochizuki works with the tuple of realified Frobenioids $(\sFrob(L)_1^\R,\ldots, \sFrob(L)_\ells^\R)$ which provides, by the same degree mapping, the hyperplanes considered here.
\item The existence of the phenomenon and properties established by \Cref{th:hyperplane} is vital to establishing  \cite[Corollary 3.12]{mochizuki-iut3}. But this type of proposition is claimed without proof in \cite{mochizuki-iut3}, \cite{fucheng}, \cite{yamashita}. Notably, because of the scaling property \cite[Theorem 4.2.2.1]{joshi-teich-rosetta},  the hyperplane $H_{\by_{j}}$ comes equipped with a scaling factor of $j^2$ (for appropriate coordinate $v\in\V_L$ depending on the semi-stable elliptic curve as in \cite{joshi-teich-rosetta}). The lack of proof of this property in \cite{mochizuki-iut3} leads to the assertions of impossibility of this phenomenon in \cite[Section 2.2, Page 10]{scholze-stix}.
\eenum
\erem

\subsection{$\yadl$, $\yadl^{max}$ and $\yadl^\R$ are  parameter spaces for deformations of a number field}
\bthm\label{th:parameter-space-of-deformations} Let $L$ be a number field with no real embeddings. Then $\yadl$,   $\yadl^{max}$ and $\yadl^\R$, are all metrisable parameter spaces of arithmeticoids of $L$ i.e. notably $\yadl$ is a metrisable space parameterizing deformations of $L$.
\ethm
\bp
Evidently, by \Cref{def:arithmeticoid}, arithmeticoids of $L$ are parameterized by points of $\yadl$ and its metrisability was proved in \Cref{th:distance-bet-arith}. The assertion for $\yadl^{max}$ and $\yadl^\R$ is also clear from this.
\ep

\brem 
This should be considered as the number field analog of the well-known fact that the classical Teichmuller space is a metrisable parameter space for Riemann surfaces \cite{imayoshi-book}.
\erem

\subsection{Arithmeticoids provide anabelomorphic deformations of a number field}
An arithmeticoid of $L$ is an anabelomorphic deformation of the arithmetic of the number field $L$ i.e. one has the following.
\bthm\label{th:anabelomorphic-deforms} 
Let $L$ be a number field and let $\arith{L}_{\by_1},\arith{L}_{\by_2}$ be two arithmeticoids of $L$. For $i=1,2$ write $L_i\subset \arith{L}_{\by_i}$ for the embedded number subfield $L$ and similarly for any $v\in\V_L$, write $L_{i,v}$ for the embedded $L_v$. Then
\benumlab
\item  the isomorphism of $G_{L_1}\isom G_L\isom G_{L_2}$  provides an isomorphism of the multiplicative monoids
$$L_1^*\isom L^*\isom L_2^*$$
provided by the isomorphism $\gal(\bL_1/L_1)\isom \gal(\bL_2/L_2)$;
\item similarly one has for each $v\in\V_L$ the isomorphism $G_{L_{1,v}}\isom G_{L_{2,v}}$ provides an isomorphism of the topological monoids
$$L_{1,v}^*\isom L_{2,v}^*.$$
\item however there may be no topological isomorphism of $K_{1,v}\isom K_{2,v}$ and hence no isomorphism between $L_1\into L_{1,v}$ and $L_2\into L_{2,v}$ compatible with their embeddings in $L_1\into L_{1,v}\into K_{1,v}$ and $L_2\into L_{2,v}\into K_{2,v}$. 
\item The two arithmeticoids $\arith{L}_{\by_1},\arith{L}_{\by_2}$ need not be equivalent and need not even be topologically equivalent;
\item  In particular, the difference between the two  avatars of $L_1, L_2$ of $L$ given by  $\arith{L}_{\by_1},\arith{L}_{\by_2}$ can be quantified by means of \Cref{th:distance-bet-arith}.
\eenum
\ethm
\bp 
The first assertion is clear from the definition. The second assertion is a consequence of the fact that the monoid $L^*$ is amphoric (see \cite{joshi-anabelomorphy}) i.e. an invariant of the isomorphism class of the topological group $\gal(\bL/L)$ (and by classical anabelian reconstruction of number fields  $L^*$ may be reconstructed from the topological $\gal(\bL/L)$--see \cite{hoshi-number}). 

The third assertion is an important observation of \cite{joshi-teich,joshi-untilts}.  The fourth is immediate from \Cref{th:existence-inequivalent-arithmeticoids} and the last assertion is clear from  \Cref{th:distance-bet-arith}.
\ep

\brem\ 
\benumlab
\item Thanks to \Cref{th:anabelomorphic-deforms}{\bf(4)} one can assert that  arithmeticoids provide topological deformations of a number field in general and that each arithmeticoid provides an avatar of the number field in which its field arithmetic is slightly different from the corresponding structure in another inequivalent arithmeticoid.
\item The theory of arithmeticoids is compatible with Mochizuki's viewpoint (and my viewpoint) of the fluidity of arithmetic structures and in fact provides a precise version of it in which differences between two such avatars of $L$ can be quantified by means of \Cref{le:yadl-is-metric}. 
\item \Cref{th:anabelomorphic-deforms}{\bf(2)} allows us to view $\arith{L}$ as providing an anabelomorphic (i.e. keeping the isomorphism class of the absolute Galois group $G_L$ fixed) deformation of the arithmetic of a number field. To put it simply: \textit{An arithmeticoid of $L$ is an anabelomorphic deformation of the arithmetic of the number field $L$.}
\eenum
\erem

\subsection{Deformations of a number field are linked by $\Bl$}
The precise assertion is the following. [For a similar point of view see \cite{joshi-formal-groups}.]

\bthm 
Let $\arith{L}_{\by}$ be an arithmeticoid. Let $R_{\by}$ be its arithmetic ring. Then there exists a closed ideal $\mathbb{I_{\by}}\subset \Bl$ such that 
$$\Bl/\mathbb{I_{\by}}=R_{\by}.$$
\ethm
\bp 
For each $v\in \V_L$, let $\fm_{\by,v}\subset B_{\lvbh,L_v}$ be the closed maximal ideal corresponding to the (closed classical) point $\by_v\in\syflv$. Then the product ideal $$\mathbb{I}_{\by}=\prod_{v\in\V_L} \fm_{\by_v}$$ is closed in product topology on $\Bl$. One has a natural surjection
$$\Bl=\prod_{v\in\V_L} B_{\lvbh,L_v}/\fm_{\by_v}\to \prod_{v\in\V_L} K_{\by,v},$$
which is given on the factor corresponding to $v$ by quotienting by the maximal ideal $\fm_{\by,v}$. The kernel of this surjection is $\mathbb{I}_{\by}$. 
\ep

\subsection{Frobenioid associated to a number field and to an  arithmeticoid}\label{ss:frobenioids}
For a domain $R$ let $R^\tr$ be the multiplicative monoid of non-zero elements of $R$.
Let me demonstrate that there is a natural Frobenioid associated to an arithmeticoid. For Mochizuki's definition of Frobenioid of a number field it is enough to understand \cite[Example 6.3]{mochizuki-frobenioid1}.

\bdefn\label{def:frobenioid-number-field}
The Frobenioid associated to the number field $L$ \cite[Example 6.3]{mochizuki-frobenioid1} is the triple $$\sFrob(L)=(L^*,\Phi(L), L^*\to \Phi(L)^{gp}) $$
where $$\Phi({L})=\bigoplus_{v\in\V_L} (\O_{L_v}^\tr)/(\O_{L_v}^*)$$
and $$L^*\to \Phi(L)^{gp}=\bigoplus_{v\in\V_L} {L_v}^*/\O_{L_v}^*$$ is the natural homomorphism.
\edefn

Now one can define a Frobenioid associated to any arithmeticoid as follows:
\bdefn\label{def:frobenioid-arithmeticoid}
Let $(\arith{L}, \iota_L:L\to \arith{L})$ be an arithmeticoid. Write $\iota_{L_v}:L_v\to K_v$.
Let $\sFrob(\arith{L})$ be the Frobenioid defined as follows. Let $$\Phi(\arith{L})=\bigoplus_{v\in\V_L} (\O_{K_v}^\tr)/\O_{K_v}^*,$$
and let $$L^*\to \Phi(\arith{(L)})^{gp}$$
be the natural homomorphism given in the $v$ component by $$x\mapsto \iota_{L_v}(x)\bmod{\O_{K_v}^*}\in K_v^*/\O_{K_v}^*$$
Then the Frobenioid associated to $\arith{L}$ is given as follows:
\be \sFrob(\arith{L})=\left(L^*,\Phi(\arith{L}), L^*\to \Phi(\arith{L})^{gp}\right)
\ee
\edefn
Note that $K_v^*/\O_{K_v}^*$ is the value group of $K_v$ and by \cite{scholze12-perfectoid-ihes} this value group is naturally isomorphic to the value group of its tilt $K_v^\flat\isom \lvbh$.

\brem
The data of a Frobenioid of a number field is much weaker than the data of an arithmeticoid. The next proposition shows this Frobenioid of an arithmeticoid is naturally isomorphic to the perfection of the Frobenioid associated to a number field in \cite{mochizuki-frobenioid1}. 
\erem

\brem 
Let me remark that if $X$ is a reasonable adic space--for example $X$ arises as the adic space associated to any quasi-projective variety, then one can define a Frobenioid associated to $X$. This remark will not be used in this paper.
\erem

\subsection{The isomorphism $\sFrob(\arith{L})\isom \sFrob(L)^{pf}$} 
The assertion is the following:
\bpro\label{pr:isom-frobenioids1} 
Let  $\arith{L}$ be an arithmeticoid of a number field $L$. Then there is a natural isomorphism of $$\sFrob(\arith{L})\isom \sFrob(L)^{pf}$$ of Frobenioids, where
$$\sFrob(\arith{L})=\left(L^*,\Phi(\arith{L}), L^*\to \Phi(\arith{L})^{gp}\right)$$ is the Frobenioids defined above and the perfection of the Frobenioid of the number field $L$ defined by Mochizuki in \cite[Example 6.3]{mochizuki-frobenioid1} and given by
$$\sFrob(L)^{pf}=\left(L^*,\Phi(L)=\bigoplus_{v\in\V_L} \O_{L_v}^\tr/\O_{L_v}^*, L^*\to (\Phi(L)^{pf})^{gp}\right).$$
 On the other hand many non-isomorphic arithmeticoids of $L$ provide isomorphic Frobenioids isomorphic to $\sFrob(L)^{pf}$.
\epro
\bp 
The proof is immediate from the definitions. 
The perfection, $\sFrob(L)^{pf}$ replaces $\Phi(L)$ by its perfection $\Phi(L)^{pf}$.
\ep

\subsection{Realified Frobenioids}\label{ss:realified-frobenioids}  \cite{mochizuki-frobenioid1} and \iut\ also work with  realified Frobenioid $\sFrob(L)^{\R}$. This is possible in the approach discussed above on noting that the value group of an algebraically closed perfectoid field is rank one i.e. naturally given by specifying a valuation taking values in $\R$. So one can view the Frobenioid of $\sFrob(\arith{L})$ as equipped with a natural realification $$\sFrob(\arith{L})^\R$$ and then extend the above isomorphism to an isomorphism of realified Frobenioids $$\sFrob(L)^\R\isom \sFrob(\arith{L})^\R.$$

By the construction of $\yadl^\R$ in \ssep\ref{ss:realified-curves}, if $\by\in \yadl^\R$ then $\sFrob(\arithl_\by)$ natural equals its realification in the sense of \cite{mochizuki-frobenioid1} and \ssep\ref{ss:frobenioids}: 
$$\sFrob(\arithl_\by)=\sFrob(\arithl_\by)^\R.$$
\brem 
\textit{Thanks to \Cref{pr:isom-frobenioids1}, one could say that perfections of the Frobenioids of  $p$-adic fields in \iut\ are a proxy for  (algebraically closed) perfectoid fields in \iut} i.e. my central observation is that in \iut\ the passage to perfections and realifications of multiplicative (attempt to) provide functionality similar to that provided by perfectoid fields (and their tilts) in an arithmeticoids.  

More precisely one needs the functionality provided by arithmeticoids in \iut. Notably Mochizuki's Hodge-Theaters requires the notion of  arithmeticoids of number fields rather than the realification of Frobenioids of number fields. By the above proposition and the preceding theory of arithmeticoids, greater precision is gained  by doing so (even for \iut), and at the same time one still may invoke results of Mochizuki's \cite{mochizuki-frobenioid1}.
\erem

\section{Arithmeticoids and Arithmetic Teichmuller spaces}\label{se:arith-and-geom}
\newcommand{\fjxl}{\mathfrak{J}(X/L)}
\subpara Let $L$ be a number field. Let $X/L$ be geometrically connected, smooth, hyperbolic curve over $L$.  Let $\fjxl$ be the adelic arithmetic Teichmuller space considered in \cite{joshi-teich}. A point of the adelic Teichmuller space $\fjxl$  may be understood as a collection of local arithmetic holomorphic structures for each $v\in \V_L$:
\be 
(Y/L_v, L_v\into K_v, K_v^\flat\isom \lvbht, *_{K_v}:\sM(K_v)\to \yan_{L_v})_{v\in\V_L}
\ee
such that for each $v\in\V_L$, $Y/L_v$ is tempered anabelomorphic with $X/L_v$ (this condition is independent of the choice of arithmetic holomorphic structures on $X$ and $Y$).

\bpro\label{pr:arith-datum-frob} 
Let $(Y/L_v, L_v\into K_v, K_v^\flat\isom F_v, \sM(K_v)\to \xan_{L_v})_{v\in\V_L}\in\fjxl$ be any point of the adelic arithmetic Teichmuller space. Then this provides a natural arithmeticoid 
$\arith{L}$ given by  $$\arith{L}=(L_v\into K_v, K_v^\flat\isom F_v)_{v\in\V_L},$$
and hence a natural isomorphism of perfect Frobenioids (\Cref{pr:isom-frobenioids1}) 
$$\sFrob(\arith{L})\isom \sFrob(L)^{pf}$$
and an isomorphism of realified Frobenioids $$\sFrob(\arith{L})^\R\isom \sFrob(L)^{\R}$$
given  in \ssep\ref{ss:realified-frobenioids}.
\epro
\bp 
The proof is immediate from the definition of arithmeticoid and from \Cref{pr:isom-frobenioids1} and \ssep\ref{ss:realified-frobenioids}.
\ep
\brem 
In \cite{joshi-teich-rosetta} it is shown that one may canonically construct Mochizuki's Hodge Theater (and its variants) using the datum of any adelic point of $\fjxl$.
\erem

\subparat{Deforming the  number field also deforms local analytic geometries}
The theorem given below is an important consequence of \cite{joshi-teich} and demonstrates how the choice  of a deformation of a number field (i.e. of an arithmeticoid of $L$) changes local analytic geometries at an arbitrary set of primes of the number field--while keeping the \'etale fundamental group of $X/L$ and the tempered fundamental group of $X/L_v$ for each prime $v$ fixed! The point to take away from this theorem proved below is that deforming the number field also deforms geometry of varieties over that number field (by deforming $p$-adic analytic geometries).

\bthm\label{th:consequence-inequivalent-arithmeticoids}
Let $L$ be a number field and let $\emptyset\neq  S\subset \vnl$ be a non-empty set of non-archimedean primes of $L$ ($S$ is not be assumed to be finite). Let $\arith{L}_{\by_1}$ be an arithmeticoid of $L$. Let $X/L$ be a geometrically connected smooth, projective variety over $L$. Then there exists an arithmeticoid $\arith{L}_{\by_2}$ of $L$ such that 
\benumlab
\item the arithmeticoids $\arith{L}_{\by_1}$ and $\arith{L}_{\by_2}$ are not topologically equivalent (and hence not equivalent), and
\item for each $v\in S$, the $\Q_v$-analytic spaces $\xan_{K_{1,v}}$ and $\xan_{K_{2,v}}$ are not isomorphic.
\item The pairs $(X,\arith{L}_{\by_1}),(X,\arith{L}_{\by_2})$ each provide distinct members of the isomorphism class of $X/L$.
\item Each pair $(X,\arith{L})$  provides an isomorph of the \'etale fundamental group $$\pi_1^{et}(X/L,*_{\bL}:\Spec(\bL)\to X)$$ for some choice of an $\bL$-geometric base-point with $\bL$ provided by the given arithmeticoid $\arith{L}$, and for each $v\in \V_L$, an isomorph of the tempered fundamental group $$\pi_1^{temp}(\xan_{L_v},*_{K_v}:\sM(K_v)\to \xan_{L_v})$$
for some choice of a $K_v$-geometric base-point with $K_v$ provided by $\arith{L}$; 
and for each $v\in \V_L$, an isomorph of the geometric tempered fundamental subgroup $$\pi_1^{temp}(\xan_{K_v},*_{K_v}:\sM(K_v)\to \xan_{K_v}) \subset \pi_1^{temp}(\xan_{L_v},*_{K_v}:\sM(K_v)\to \xan_{L_v}) $$
(for some choice of a $K_v$-geometric base-point with $K_v$ provided by $\arith{L}$)
\eenum
\ethm
\bp 
Let me choose an arithmeticoid $\arith{L}_{\by_2}$ as follows. For each $v\in S$, let me choose $\by_{2,v}\in\syflv$ such that  the perfectoid residue fields  $K_{1,v}$ and $K_{2,v}$ are not topologically isomorphic (these choices are possible by the axiom of choice).  The existence of such fields is established by \cite{kedlaya18}. For $v\not\in\ S$, choose $\by_{2,v}=\by_{1,v}$. This gives a $\by_2\in\yadl$. By construction, for all $v\in S$, one has $\by_{1,v}\neq\by_{2,v}\in \syflv$. Now that one has chosen $\by_2$, one obtains, by definition, an  arithmeticoid $\arith{L}_{\by_2}$ such that $\arith{L}_{\by_1}$ and $\arith{L}_{\by_2}$ are not equivalent and also not topologically equivalent by construction. This proves {\bf(1)}.

The second assertion is an immediate consequence of the choice of $\arith{L}_{\by_2}$ and one of the main theorems of \cite[Section 3]{joshi-teich}.

The third assertion is clear. The fourth assertion is a consequence of \cite[Section 3]{joshi-teich}. This proves the theorem.
\ep

\newcommand{\hol}[3]{\mathfrak{hol}_{#1}(#3)_{#2}}
\newcommand{\holt}[2]{\mathfrak{hol}(#1)_{#2}}

\brem 
It is precisely because of \Cref{th:consequence-inequivalent-arithmeticoids} that one can say the theory developed in the present series of papers is the adelic analog of classical Teichmuller Theory and \Cref{th:consequence-inequivalent-arithmeticoids} leads me to the following natural definition.
\erem
\begin{defn}\label{def:global-holomorphoids}
Let $L$ be a number field (assumed to have no real embeddings) and let $X/L$ be a quasi-projective variety over $L$. Then a \textit{global holomorphoid of $X/L$}, or more simply an \textit{arithmetic holomorphoid of $X/L$}, or even more simply a \textit{holomorphoid of $X/L$} is the datum $(X, L\into \arithl_\by)$ consisting of an arithmeticoid $\arithl_\by$  with  $\by=(y_v)_{v\in\vl}$ and, for each $v\in\vl$, a morphism 
\be\label{eq:geom-base-point-holomorphoid} \sM(K_{y_v})\to \xan_{L_v}\ee
of analytic spaces. I will write the holomorphoid $(X,\arith{L})$ as $(X/\arith{L})$, or even more suggestively as $\holt{X/L}{}$, or  $\holt{X/L}{\by}$ (if I want to indicate the arithmeticoid datum $\arith{L}_{\by}$ (with $\by\in\yadl$) explicitly). Evidently, if $(X,\arithl_\by)$ is a holomorphoid then 
$$(X/L_v, L_v\into K_{y_v}, K_{y_v}^\flat\isom F_{y_v}, \sM(K_{y_v})\to \xan_{L_v})_{v\in\V_L}$$ is a point i.e. an object of the adelic Arithmetic Teichmuller Space $\fjxl$ of \cite{joshi-teich,joshi-untilts}.
\end{defn}

\brem 
To place this in the context of classical Teichmuller Theory, let $\Sigma$ be a connected, compact Riemann surface and let $\abs{\Sigma}$ be the underlying topological space. Then any Riemann surface $\Sigma'\in T_{\Sigma}$ in the Teichmuller space $T_\Sigma$ of $\Sigma$ may be called an \textit{holomorphoid} of $\abs{\Sigma}$. In particular,  $\Sigma$ is the tautological holomorphoid of $\abs{\Sigma}$. Giving the data of a holomorphoid of $\abs{\Sigma}$ is equivalent to giving the pair $(\Sigma', f:\Sigma'\to \Sigma)$ consisting of a Riemann surface $\Sigma'$ and a quasi-conformal mapping $f:\Sigma'\to \Sigma$.
\erem

\numberwithin{equation}{subsection}
\section{Arithmeticoids and Galois cohomology}\label{se:arith-and-galois-cohomology}
In this section I want to describe Galois cohomology provided by an arithmeticoid. An important point which emerges is that Galois cohomology classes of geometric or arithmetic interest have dependence on the chosen arithmeticoid and by means of the results established in \cite[Section 7]{joshi-untilts}, one may compare cohomology classes provided by distinct arithmeticoids in one common location such as the Galois cohomology of a fixed arithmeticoid. This leads to the key idea of collation of cohomology classes arising from distinct avatars of a number field which is central to the formulation of \cite[Theorem 3.11, Corollary 3.12]{mochizuki-iut3}. 

\subparat{Standard arithmeticoids}\label{pa:stand-arithmeticoid} Let $L$ be a number field, let $\arith{L}_\by$ be an arithmeticoid of $L$. I will chose an arithmeticoid $\arith{L}_{\by_0}$ which I will call the \textit{standard arithmeticoid of $L$}. This is chosen as follows. For each prime $v\in \vnl$ choose a closed classical point $y_{v}\in\syflv$ which lies in the fiber of over the canonical point of $\sxflv$. For $v\in\vlarch$, let $y_v$ be any point of $\syflv$ (remember that for $v\in\vlarch$, $\vphi_v=1$). The point $\by_0=(y_v)_{v\in\vl}\in \yadl$ thus provides by definition an arithmeticoid $\arith{L}_{\by_0}$. I will refer to $\arith{L}_{\by_0}$ as the \textit{standard arithmeticoid of $L$} and $\by_0\in\yadl$ will be called the \textit{standard point} of $\yadl$.

Note that $\arith{L}_{\by_0}$ is not unique! There are many possible arithmeticoids which may be used as standard arithmeticoids of $L$.
\subparat{Galois cohomology}\label{ss:galois-cohomology-facts} Let $L\into \bL\into \arith{L}_{\by}$ be the algebraic closure of $L$ provided by $\arith{L}_{\by}$ and let $G_{L}=\gal(\bL/L)$ be the absolute Galois group of $L$ provided by $\arith{L}_{\by}$. For each $v\in \vl$ let $K_v$ be the algebraically closed perfectoid field provided by $\by$. Let $L_v\into \bL_v\in K_v$ be the algebraic closure of $L_v$ provided by $K_v$ and let $G_{L_v;K_v}=\gal(\bL_v/L_v)$ be the absolute Galois group of $L_v$ given using the algebraic closure $\bL_v\subset K_v$ of $L_v$ in $K_v$. 

For non-archimedean primes $v\in \vl$ let $\Z_{p_v}(1)_{K_v}$ be the $G_{L_v;K_v}$-module of  $p_v$-power roots of unity in $K_v$ where $p_v$ is the residue characteristic of $L_v$. 

Then for each $i\geq 0$ and any prime $v\in \vl$ (archimedean or non-archimedean), one takes Galois cohomology $H^i(G_{L_v},\Z_{p_v}(1))$ to be defined as
\be
H^i(G_{L_v},\Z_{p_v}(1))=
\begin{cases}
	Ext_{MHS}^i(\Z(0),\Z(1)) & \text{if } v\in\vlarch;\\
	 H^i(G_{L_v,K_v},\Z_{p_v}(1)) & \text{if } v\in\vlnon.
\end{cases}
\ee
The following proposition summarizes the key properties of $H^i(G_{L_v},\Z_{p_v}(1))$ and $H^i(G_{L_v},\Q_{p_v}(1))$ which will be used throughout this paper.
\bpro Let $L$ be a number field as above. Then one has the following multiplicative description of $H^i(G_{L_v},\Z_{p_v}(1))$ at all primes $v\in\vl$.
\benumlab
\item For $v\in\vlarch$, one has $L_v\isom\C$ and $$H^i(G_{L_v},\Z_{p_v}(1))=Ext_{MHS}^1(\Z(0),\Z(1))\isom L_v^*\isom \C^*.$$ 
\item For $v\in\vlnon$ one has (by Kummer theory) an isomorphism of $\Z_{p_v}$-modules $$H^1(G_{L_v,K_v},\Z_{p_v}(1))\isom \invlim_n (L_v^*/L_v^{*p_v^n}),$$ 
and
$$H^1(G_{L_v,K_v},\Q_{p_v}(1))\isom \left(\invlim_n (L_v^*/L_v^{*p_v^n})\right)\tensor_{\Z_{p_v}}\Q_{p_v}.$$
\eenum
\epro
\bp 
For the first assertion one uses \cite{deligne-local} and the assumption that $L$ has no real embeddings. For the second assertion see \cite{perrin-riou1994} or \cite{nekovar93}.
\ep

\brem Note that because of these isomorphisms, the main groups of interest in this paper have a multiplicative description.
\erem

\bdefn The galois cohomology of an arithmeticoid with integral coefficients (resp. with rational coefficients) is the product given by 
$$ H^i({\arith{L}_\by},\Z(1)):=\prod_{v\in\vlnon} H^i(G_{L_v;K_v},\Z_{p_v}(1))\times \prod_{v\in\V_L^{arc}} Ext^i(\Z(0),\Z(1)),$$
and respectively the product
$$H^i({\arith{L}_\by},\Q(1)):=\prod_{v\in\vlnon} H^i(G_{L_v;K_v},\Q_{p_v}(1))\times \prod_{v\in\V_L^{arc}} Ext^i(\Q(0),\Q(1)).$$
 I will refer to $H^i({\arith{L}_\by},\Z(1))$ (resp. $H^i({\arith{L}_\by},\Q(1))$) as the \textit{integral Galois cohomology provided by the arithmeticoid $\arithl_\by$} (resp. \textit{rational Galois cohomology provided by the arithmeticoid $\arithl_\by$}). 
\edefn
Note that $H^i({\arith{L}_\by},\Z(1))$ is a module over the ring $\left(\prod_{v\in\V_L^{non}}\Z_{p_v}\right)\times \Z$ with the archimedean component treated as an abelian group i.e. a $\Z$-module and $H^i({\arith{L}_\by},\Q(1))$ is a module over $\left(\prod_{v\in\V_L^{non}}\Q_{p_v}\right)\times \Q$.

\brem 
One way to think about the above definition is that there are many Galois cohomology theories--each indexed by an arithmeticoid. As discussed in \cite{joshi-formal-groups}, in the case of topological spaces, formal groups play a central role via \cite{quillen69}; likewise, in the present theory formal groups (\cite{joshi-formal-groups}) play a similar role from the point of view of Galois Cohomology theory.
\erem
\subsection{Bloch-Kato subspaces} The Bloch-Kato subpaces of these groups may also be defined (see \cite{bloch90} or \cite{nekovar93}) as follows:
let $$H^i_f({\arith{L}_\by},\Z(1)):=\prod_{v\in\vlnon} H^i_f(G_{L_v;K_v},\Z_{p_v}(1))\times \prod_{v\in\V_L^{arc}} Ext^i(\Z(0),\Z(1)),$$ and also
$$H^i_f({\arith{L}_\by},\Q(1)):=\prod_{v\in\vl} H^i_f(G_{L_v;K_v},\Q_{p_v}(1))\times \prod_{v\in\V_L^{arc}} Ext^i(\Q(0),\Q(1)).$$ I will refer to $H^i_f({\arith{L}_\by},\Z(1))$ (resp. $H^i_f({\arith{L}_\by},\Q(1))$) as the integral Fontaine subspace of $H^i({\arith{L}_\by},\Z(1))$  (resp. the rational Fontaine subspace of $H^i({\arith{L}_\by},\Q(1))$). 
For any $v\in\vlnon$ one has an isomorphism (see \cite{bloch90}, \cite{perrin-riou1994} or \cite{nekovar93})
\be  
H^1_f(G_{L_v},\Q_{p_v}(1))\isom  \left(\invlim_n(\O_{L_v}^*/\O_{L_v}^{*p_v^n})\right)\tensor_{\Z_{p_v}}\Q_{p_v}.
\ee

\subparat{Amphoricity of cohomology of an arithmeticoid} Now suppose $\arithl_{\by_1},\arithl_{\by_2}$ are two arithmeticoids of $L$. Then one has the following (obtained by applying \cite[Theorem 7.7.1]{joshi-untilts} at each prime $v\in\vl$):
\bpro\label{pr:arith-isom-gal-cohom}\
\benumlab 
\item The galois cohomology groups $$H^i(\arith{L}_{\by_1},\Z(1)), H^i_f(\arith{L}_{\by_1},\Z(1)), H^i(\arith{L}_{\by_1},\Q(1)), H^i_f(\arith{L}_{\by_1},\Q(1))$$ are $G_L$-amphoric. 
\item In other words: 
for all integers $i\geq0$ one has  isomorphism of groups
$$H^i(\arith{L}_{\by_1},\Z(1))\isom H^i(\arith{L}_{\by_2},\Z(1)),$$
and
$$H^i(\arith{L}_{\by_1},\Q(1))\isom H^i(\arith{L}_{\by_2},\Q(1)),$$
\item and also of the Fontaine subspaces:
$$H^i_f(\arith{L}_{\by_1},\Z(1))\isom H^i_f(\arith{L}_{\by_2},\Z(1)),$$
and
$$H^i_f(\arith{L}_{\by_1},\Q(1))\isom H^i_f(\arith{L}_{\by_2},\Q(1)).$$
\eenum
\epro
\bp 
The isomorphism class of $G_L$ determines, by \Cref{cor:fun-actions}, the isomorphism class of $\sGal{L}$ and  the isomorphism class of each of its factors $G_{L_v}$, so one reduces to the local case. For proving the local case  it is sufficient to invoke  \constrone{Theorem }{pr:galois-cohomology-and-teich} for each factor. This proves all the assertions.
\ep

\subsection{Collation of cohomology classes} The following proposition is now proved along the lines of \constrone{Corollary }{cor:collation-classes} and provides a way of reading cohomology classes arising from distinct arithmeticoids in the Galois cohomology of a fixed arithmeticoid. 
\bpro\label{pr:collation-principle} 
Let $X/L$ be a geometrically connected, smooth, quasi-projective variety over $L$. Let $\sA\subset\yadl$ be a collection of arithmeticoids and write $$\sA_X=\{X/\arith{L}_\by:\by\in \sA\}$$ for the collection of holomorphoids of $X/L$ provided by the arithmeticoids in $\sA$. For each holomorphoid $X/\arith{L}_\by\in\sA_X$ let $$\Psi_\by\subseteq  H^1({\arith{L}_\by},\Z(1))$$ be a collection of cohomology classes arising from $X/\arithl_\by$ (for example one can take $\Psi_\by=\{{\bf q}_{C/\arithl_\by}\}$ constructed in \eqref{eq:Tate-class-galois-cohom}). Let $\by_0$ be the standard arithmeticoid defined by \ssep\ref{pa:stand-arithmeticoid}. Then there exists a subset 
$$\Psi_{X,\sA}\subseteq  H^1({\arith{L}_{\by_0}},\Z(1))$$
which consists of the union,  over $\by\in\sA$, of the images of all the subsets $\Psi_\by$ under all the  isomorphisms (of topological groups) of each factor of $$H^1(\arith{L}_\by,\Z(1))\isom  H^1(\arith{L}_{\by_0},\Z(1))$$
provided by \Cref{pr:arith-isom-gal-cohom}.
\epro
\bp 
The proof is clear from \constrone{Corollary }{cor:collation-classes}.
\ep
\brem 
One can formulate a similar assertion for coefficients $\Z_p(1)$, $\Q_p(1)$ (for $p\in\V_{\Q}$) etc. 
\erem

\subparat{A worked example: the cohomology class of a Tate elliptic curve} Let $C/\arith{L}_\by$ be an elliptic curve over an arithmeticoid of $L$. Let $v_1,\ldots,v_n$ be the non-archimedean primes of $L$ at which $X$ has semi-stable reduction. Observe that (by \ssep\ref{ss:galois-cohomology-facts}) regardless of whether $v$ is archimedean or non-archimedean, $H^1(G_{L_v},\Q_{p_v}(1))$ has a multiplicative description. This allows one to construct an element for each prime $v\in\vl$ as follows. 

For $j=1,\ldots, v_n$, at each $v_j$ one obtains a Tate elliptic curve and hence a Galois cohomology class in $H^1(G_{L_{v_j}},\Q_{p_{v_j}}(1))$ given by the Tate parameter $q_{v_j}$ as follows. Using the Tate parameter $q_{v_j}$ one can construct an explicit Galois cohomology class which is explicitly described for instance in \cite[II.4]{berger2002}. 

 Take a compatible system $\left\{(q_{v_j})^{1/p_v^n}\right\}_{n\geq0}$ of $p^n$-th roots in $\bL_v$ of $$q_{v_j}\in \bL_v.$$ This provides a Galois cohomology class in $H^1(G_{L_{v_j}},\Q_{p_{v_j}}(1))$ which is explicitly described in \cite[II.4]{berger2002} or \cite[Chapitre 10]{fargues-fontaine} and which I will write by abuse of notation as $$\left\{(q_{v_j})^{1/p_v^n}\right\}_{n\geq0}\in H^1(G_{L_{v_j}},\Q_{p_{v_j}}(1)).$$ 

At all primes of non-semi-stable reduction one  takes the class to be trivial cohomology class $1$ (written multiplicatively). At any archimedean prime $v$ one can take the Schottky parameter $q_v\in Ext^1(\Z(0),\Z(1))$ (see \ssep\ref{ss:schottky-parameterization}). Explicitly

\be  
q^{coh}_{C,y_v}=\begin{cases}
	1 \in H^1(G_{L_v},\Z_{p_v}(1)) & \text{if } v\in\vnl, \text{ and } v\neq v_1,\ldots,v_n,\\
	\left\{(q_{v_j})^{1/p_v^n}\right\}_{n\geq0}  \in H^1(G_{L_v},\Q_{p_v}(1)) & \text{if } v=v_j \text{ for } j=1,\ldots,n;\\
	q_v \in Ext^1_{MHS}(\Z(0),\Z(1)) & \text{if } v\in\vlarch.
\end{cases}
\ee

Thus one obtains a well-define element  
\be\label{eq:Tate-class-galois-cohom} {\bf q}_{C,\by}\in H^1({\arith{L}_\by},\Z(1)).
\ee

[In \cite{joshi-teich-rosetta}, I will provide a variant of this construction which provides a crystalline cohomology class associated to a Tate elliptic curve.]

\brem\ 
\benumlab
\item In \cite{joshi-teich-rosetta} this collation  principle (\Cref{pr:collation-principle}) is applied to certain crystalline cohomology classes constructed from a suitably chosen set of holomorphoids of a Tate elliptic curve at each odd prime of semi-stable reduction. This leads one to the set underlying the assertion of \cite[Corollary 3.12]{mochizuki-iut3}. 
\item Importantly for \iut\ and especially \cite[Corollary 3.12]{mochizuki-iut3}, one requires crystalline cohomology classes constructed from Tate elliptic curves (constructed explicitly in \cite{joshi-teich-rosetta}).
\item In \cite[Corollary 3.12]{mochizuki-iut3} (detailed in \cite{joshi-teich-rosetta}), one applies \Cref{pr:collation-principle} to collate these crystalline classes arising from a suitable set of arithmeticoids (given explicitly in \cite{joshi-teich-rosetta}) to arrive at  Mochizuki's theta-values set i.e. the set on the left hand side of  \cite[Corollary 3.12]{mochizuki-iut3} and the said corollary provides a lower bound on the (suitably defined) volume of this set. 
\eenum
\erem

\section{Heights and arithmeticoids}\label{se:heights}
\subparat{Preliminaries} For defining heights of algebraic numbers (such as elements of $L$) I will follow the treatment of \cite{bombieri-gubler}.  Let $$\log:\R^+ \to \R$$ be the natural logarithm of a positive real number. Following \cite{bombieri-gubler}, let
$$\logp(t)=\max(0,\log(t))$$
and extend $\logp$ from $\R^+\to \R$ to $$\logp:\R^{\geq0} \to \R$$ 
 by setting $$\logp(0)=0.$$
 
\subparat{Height function associated to an arithmeticoid}
\bdefn\label{def:heights} 
Let $\arith{L}_\by$ be an arithmeticoid of $L$ with $\by=(y_v)_{v\in\vl}$. Choose the standard normalization for the valuations on $\arith{L}_\by$ i.e. work with $\arith{L}_{\by}^{nor}$ i.e. work with the normalization coordinate   $\alpha_{\by}=(\alpha_v)_{v\in\vl}$ of $\by$ (\cref{ss:normalization})  (I will suppress the notational distinction between $\arith{L}_\by$ and $\arith{L}_{\by}^{nor}$ for simplicity). Let $\bL$ be the algebraic closure of $L$ provided by this choice of arithmeticoid. Let $P=(x_0,x_1,\ldots,x_n)\in \P^n(L)$ be a point in projective space. Then the height 
$$h_{\arith{L}_{\by}}(P)=h_{\arith{L}_{\by}}((x_0,x_1,\ldots,x_n))=\sum_{v\in\V_L}\max_j(\log(\abs{x_j}_{K_{y_v}}^{\alpha_v})).$$

If $z\in L$ is an algebraic number, then its height is defined to be  $$h_{\arith{L}_\by}(z)=h_{\arith{L}_\by}((1,z))\in \P^1(L)$$ and  is explicitly given by
$$h_{\arith{L}_\by}(z)=\sum_{v\in\V_L} \max(0,\log(\abs{z}_{K_{y_v}}^{\alpha_v}))=\sum_{v\in\V_L} \logp\abs{z}_{K_{y_v}}^{\alpha_v}.$$
\edefn

\brem 
Using \cite[Lemma 1.3.7 and 1.5.1]{bombieri-gubler}, one may also define height $h_{\arith{L}_{\by}}(P)$ of any point $P\in\P^n(\bL)$. This will be used, but is left for the reader.
\erem

\blem 
This definition is independent of the choice of homogeneous coordinate representatives of $P$ i.e. $(x_0,x_1,\ldots,x_n)$ and $(\lambda\cdot x_0,\lambda\cdot x_1,\ldots,\lambda\cdot x_n)$ provide the same height.
\elem

\bp 
This is a consequence of the fact that we work with a normalized arithmeticoid i.e with absolute values normalized so that the product formula \eqref{eq:prod-formula1} holds.
\ep

\subparat{Dependence of heights on the arithmeticoid}\label{re:normalization2} 
The following is essential to this theory:
\bpro\label{th:height-depends-on-arithmeticoid} 
Let $\arith{L}$ be an arithmeticoid of $L$. Then the associated height  function $$h_{\arith{L}}:\P^n(\bL) \to \R$$  depends non-trivially on $\arith{L}$ i.e. $h_{\arith{L}}$ depends on arithmetic holomorphic structures chosen to compute it. Notably the symmetries given by \Cref{cor:galois-action-on-arithmeticoids} which act upon the set of arithmeticoids change the height functions provided by arithmeticoids.
\epro 
\bp
As remarked in \Cref{re:normalization1}, $y\mapsto \abs{p}_{K_y}$ is a non-constant function on $\syflv$. This means that the function $$h_{\arith{L}}:\P^n(\bL) \to \R$$
depends on the choice of the arithmeticoid $\arith{L}$ in a non-trivial way. So the assertion follows.
\ep

\brem 
Regardless of how valuations  are normalized in an arithmeticoid, if one has two arithmeticoids $\arith{L}_{\by_1}$ and $\arith{L}_{\by_2}$, then the height computed in one may not agree with that in the other because of \Cref{re:normalization1}. Especially, the normalization of $\arith{L}_{\by_1}^{nor}$ cannot be carried over to $\arith{L}_{\by_2}^{nor}$.
\erem

\subparat{First instructive example} The following will be useful in understanding the above definition of heights.
\begin{example}\label{ex:example1}	
Let $L=\Q$, so $\V_L=\V_{\Q}=\{p=\infty \text{ or } p \text{ is a prime number} \}$. Let us choose an arithmeticoid of $L$. Let $v\in \V_L$ be a non-archimedean prime. Let $\{y_{v,n}\in\syflv:n\in \Z\}\subset \syflv=\syqp$ be the fiber,  of the morphism $\syqp\to\sxqp$, over the canonical point of $\sxflv=\sxqp$. 

For each non-archimedean prime $v=p$, chose the point $y_{v,0}=y_{p,0}$ on the non-archimedean Fargues-Fontaine curve $\syflv$. This provides the valued field $K_v=\C_p$ for each prime number $p$. 

For an archimedean prime $v=\infty$, chose a point $y_0$ on the archimedean Fargues-Fontaine curve which provides the standard archimedean absolute value on $K_v=\C_\infty=\C$. 

These choices provide an arithmeticoid of $L$ arising from $\by_0=(y_{v,0})_{v\in\V_L}\in\yadl$. 
Let us choose this arithmeticoid of $L$ arising from $\by_0=(y_{v,0})_{v\in\V_L}$ for computing heights.

Let $\arith{L}_{\by_0}$ be the above chosen arithmeticoid of $L=\Q$. Let $z=5$. Suppose the valuation on each local data of $\arith{L}_{\by_0}$, i.e.  $K_p$, for each $p\in\V_{\Q}$ is normalized so that $\abs{p}_{K_p}=\frac{1}{p}$ for $p\neq \infty$, then one has 
$$\logp(\abs{z}_{K_p})=
\begin{cases} 
		0 			& \text{ if $p\neq\infty$,} \\
		\log(z) & \text{ for $p=\infty$.}
\end{cases},$$
and 
$$h_{\arith{\Q}_{\by_0}}(z)=\log(5)$$
and if $z=\frac{1}{5}$ then
$$\logp(\abs{z}_{K_p})=
\begin{cases} 
0 			& \text{ if $p\neq 5$,} \\
\log(5) & \text{ for $p=5$.}
\end{cases}$$
and hence
$$h_{\arith{\Q}_{\by_0}}(z)=\log(5).$$

Note that for any $x\in\Q^*$ the product formula \eqref{eq:prod-formula1} means 
$$\log\abs{x}_{K_\infty}=-\sum_{p<\infty}\log(\abs{x}_{K_p}).$$
\end{example}
\subparat{Second instructive example} Now let me illustrate the consequence of global symmetries on arithmeticoids parameterized by $\yadl$ on height functions given by arithmeticoids. The important point for this discussion is that the action  $L^*\act\yadl$ which as this example will illustrate, modifies local height contributions.
\begin{example}\label{ex:example2}	
I will keep the context and the notation of \Cref{ex:example1}. Specifically this means our arithmeticoid $\arith{L}_{\by_0}$ as above. Now let us consider the action of $L^*\act\yadl$ (for $L=\Q$) i.e the action of multiplying the arithmeticoid $\arith{\Q}_{\by_0}$ considered above by some non-zero rational number. 

I will work out the case of multiplication by $n\in \Z$. Let $$n=p_1^{\alpha_1}\cdot p_2^{\alpha_2}\cdots p_k^{\alpha_k}$$ be the prime factorization of $n$. One wants to consider the effect of ``multiplying'' an arithmeticoid $\arith{L}_{\by_0}$ considered in the previous example by $n$ i.e. one wants to consider the arithmeticoid $$n\cdot \arith{\Q}_{\by_0}$$
given by \Cref{th:galois-action-on-adelic-ff}.

The  action of $n$ on $\arith{\Q}_{\by_0}$ is via the action of $L^*$ on $\yadl$, given  by \Cref{th:galois-action-on-adelic-ff},  is through the Lubin-Tate action on   its local factors.  Let $p$ be an arbitrary prime. Since $L=\Q$, I will write $p$ instead of $v$ and $K_p$ (resp. $F_p$) instead of $K_v$ (resp. $F_v$) etc. This action is obtained from the Lubin-Tate action of $\Z_p$ on $\fm_{F_p}$ and the quotient relation $$\abs{\sY_{F_p,\Q_p}}=\left(\fm_{F_p}-\{0\}\right)/\Z_p^*,$$
and $p\in\Z_p$ operates by Frobenius on $\sY_{F_p,\Q_p}$.  

Now there are two distinct possibilities. If $p$ does not divide $n$ then $n\in\Z_p^*$ and hence the action of $n$ on $\syqp$ is trivial. If $p$ divides $n$ i.e. if $p$ is one of the primes $p_1,\ldots,p_k$, then $n$ acts on  the factors $$\sY_{\C_{p_1}^\flat,\Q_{p_1}}, \sY_{\C_{p_2}^\flat,\Q_{p_2}}, \ldots, \sY_{\C_{p_k}^\flat,\Q_{p_k}}$$ of  $\yadl$ by powers of their respective Frobenius morphisms. The   powers of Frobenius being given by $\alpha_1,\alpha_2,\ldots,\alpha_k$ i.e on $$\sY_{\C_{p_i}^\flat,\Q_{p_i}},$$ the integer $n$ operates by $$\varphi_{p_i}^{\alpha_i}.$$

\brem 
As Mochizuki explains in \cite[\ssep 2.4]{mochizuki-gaussian}, one wants some sort of Frobenius morphism of a number field for bounding heights on number fields (as I show in \ssep\ref{se:appendix}, the proofs of \cite{bogomolov00} and \cite{zhang01} involve such a height function). Such a global Frobenius morphism does not exist on a fixed number field.   On the other hand there is such a global Frobenius morphism on the parameter space  of arithmeticoids $\yadl$ given by  \Cref{th:galois-action-on-adelic-ff}, \Cref{cor:galois-action-on-arithmeticoids}. This is the reason why one is interested in stabilized height functions.
\erem

In \cite[Appendix]{joshi-teich-estimates} I have established the connection between Mochizuki's $\mathfrak{log}$-links and Frobenius morphism in the theory of \cite{joshi-teich,joshi-untilts,joshi-teich-estimates}. Because of the choice of the arithmeticoid $\by_0$ in \Cref{ex:example1},  the above (Frobenius) action is equivalent to traversing the column of $\mathfrak{log}$-links in \cite{mochizuki-iut3} (see \cite[Appendix]{joshi-teich-estimates} for the elucidation of this important point) at all primes dividing $n$.

By \cite[Theorem 10.20.1]{joshi-teich-estimates} one knows that if $y_p\in\syqp$ is a closed classical point,  then Frobenius morphism $\varphi_{p}$ of $\syqp$ provides the following relationship between valuations of $p$ in the respective residue fields $K_{y_p}$ and $K_{\varphi_{p}(y_p)}$:
$$\abs{p}_{K_{y_p}} = \abs{p}_{K_{\varphi_{p}(y_p)}}^{1/p}.$$ 

Write $\vphi_i=\vphi_i$ for the Frobenius morphism at $p_i$. Hence for $n\cdot\arith{L}_{\by_0}$, this calculations shows that 
$$\abs{p_i}_{K_{y_{p_i}}} = \abs{p_i}_{K_{\varphi^{{\alpha_i}}_{i}(y_{p_i})}}^{1/p^{\alpha_i}_i} \text{ for } i=1,\ldots,k.$$

So $n$ acts on $\arith{L}_{\by_0}$, by modifying local arithmetic holomorphic structures, and modifies  the local contributions to height functions of $\arith{L}_{\by_0}$ at  all the primes $p_1,\ldots,p_k$ dividing $n$. 

On the other hand, since one is working with normalized arithmeticoids $\arith{L}^{nor}$,  the normalization constraint i.e. the product formula \eqref{eq:prod-formula1}  implies that changes in $p$-adic valuations forced by the action of $n$ must be compensated  by changes at other primes so that the normalization constraint--i.e. the product formula \eqref{eq:prod-formula1} continues to hold.

Thus the action of $n$ on arithmeticoids provides  a natural variation between the height functions $h_{\arith{L}}$ and $h_{n\cdot\arith{L}}$ of the normalized arithmeticoids $\arith{L}^{nor}$ and $n\cdot\arith{L}^{nor}$. Since the product formula is encoded in terms of the period mapping considered in \ssep\ref{ss:period-map-prod-form} (\Cref{th:hyperplane}), this means that the hyperplanes $H_{\arithl}\neq H_{n\cdot\arithl}$ in general.

This means that in comparing  height functions provided by the arithmeticoids $\arith{L}_{\by}$ and $\arith{L}_{n\cdot \by}$, the inequality
$$h_{n\cdot \arith{\by_0}} \geq h_{\arith{\by_0}},$$
which is strict in general and which can be viewed as an upper semi-continuity property of heights provided by arithmeticoids under  the global action  $L^*\act\yadl$ especially $\O_{L}^\triangleright\act \yadl$ (where $\O_L^\triangleright$ is the multiplicative monoid of non-zero integers in $L$). This upper semi-continuity property creates opportunities in Diophantine problems where one wants upper bounds on heights. 
\end{example}
\subparat{${L}^*$-Stabilized height functions}
Examples~\ref{ex:example1}, \ref{ex:example2}, and \ref{ex:example3} motivate the next definition.
\newcommand{\bfh}{\mathbf{h}}
\bdefn\label{def:stab-height} 
Let $\arith{L}_{\by}$ be an arithmeticoid of $L$. Then the \textit{$L^*$-stabilized height function} (or simply the \textit{stabilized height function}) $\bfh_{\arith{L}_\by}$ of $\arith{L}_{\by}$ is given by
$$\bfh_{\arith{L}_{\by}}(x) =\sup\left\{h_{\alpha\cdot \arith{L}_{\by}}(x)=h_{ \arith{L}_{\alpha\cdot\by}}(x) : \alpha\in L^* \right\}.$$
\edefn
\brem 
Thus $\bfh_{\arith{L}}$ is the stabilization of the height function with respect to the action $L^*\act\yadl$ i.e. with respect to the action of $L^*$ on arithmeticoids of $L$.
\erem

\subparat{The height stabilization principle}
The following is a trivial consequence of \Cref{def:stab-height}, but it will be useful to state it explicitly:
\bthm[The height stablization principle]
Let $\arith{L}_\by$ be an arithmeticoid of $L$. Then one has for all $x\in L^*$
$$\bfh_{\arithl_\by}(x)\geq h_{\arithl_\by}(x).$$
In general, the above inequality is a strict inequality.
\ethm
\brem 
One may think of the height stabilization principle as an upper semi-continuity of heights under $L^*$-action. In \iut\ a conceptually similar (and related) assertion appears in the form of \cite[Proposition 3.5]{mochizuki-iut3}.
\erem 
\subparat{Third instructive example}  The following example provides some insight into why this theorem is interesting (even though its proof is a tautology).
\begin{example}\label{ex:example3}	
	Let $L=\Q$ and let $\by_0\in\yadq$ be chosen as in \Cref{ex:example1,ex:example2}. Suppose one has three primitive integers $a,b,c$ such that
	$$a+b=c.$$
	Then it makes perfect sense to talk about $$h_{\arith{\Q}_{\by_0}}(a\cdot b\cdot c)$$
	this is the height of $a\cdot b\cdot c$ with respect to the arithmeticoid $\arith{\Q}_{\by_0}$  and also talk of $$h_{\arith{\Q}_{(a\cdot b\cdot c)\cdot\by_0}}(a\cdot b\cdot c)$$
	which is the height of $a\cdot b\cdot c$ with respect to the arithmeticoid $\arith{\Q}_{(a\cdot b\cdot c)\cdot\by_0}$ which is an arithmeticoid where the metrics (i.e. arithmetic holomorphic structures) at each prime $p|(a\cdot b\cdot c)$ are modified by the (Frobenius) action of ${(a\cdot b\cdot c)}$ on the arithmeticoid $\arith{\Q}_{\by_0}$. 
\end{example}
\brem\
\benumlab
\item 	\textit{This seems to suggest that there is some relationship between the height stablization principle and some $abc$-type  inequality!}  \textcolor{red}{Note that a precise relationship should be expected to be more complicated than the above discussion suggests.}
\item
Note that as the multiplicative structures provided by arithmeticoids (in the sense of \Cref{thm:mult-structure}) are compatible, it makes sense to think about $a\cdot b \cdot c$ (even in the context of number fields) but not about $a+b=c$ in a uniform way among all arithmeticoids provided by $\yadl$.
\item As is discussed in \ssep\ref{se:appendix} (especially \ssep\ref{ss:height-stabilzation}), a similar height-stabilization principle is at work in the proofs of \cite{bogomolov00} and \cite{zhang01}.
\eenum
\erem

\section{Arithmeticoids and Arithmetic loops and knots in a number field}\label{ss:arithmetic-loop-spaces}
The results of this subsection are not used in the rest of this paper and may be completely skipped if the reader is primarily interested in applications of this paper to \iut. Main reason this material is included here is to illustrate how my work Arithmetic Teichmuller Spaces relates (non-trivially) to other areas of interest in arithmetic geometry. For applications to the Langlands correspondence see \cite[\ssep 18]{joshi-anabelomorphy} (additionally let me remark that  the formulation of the Langlands Correspondence (global and local) requires a choice of an arithmeticoid. This raises the question of how the Langlands Correspondence for one arithmeticoid relates to that in another arithmeticoid; the results of \cite[\ssep 18]{joshi-anabelomorphy} can be considered to be preliminary results in this direction).

\subsection{Arithmetic Loops}\numberwithin{equation}{subsection}
Let $t$ be a variable and let $\overline{\C((t))}$ be an algebraic closure of $\C((t))$. The field $\overline{\C((t))}$ can be described using the Newton-Pusieux Theorem as a field of complex power series in fractional powers $t^\Q$ of $t$. 

Now let $X/\C$ be a complex variety. It is well-known (see \cite{faltings2003}) that one may think of a morphism $$\Spec(\C((t)))\to X$$ and more generally a morphism $$\Spec(\overline{\C((t)))}\to X$$ as an algebraic loop in $X$.

\textit{There is a natural arithmetic analog of this picture  which I arrived at by my work on arithmetic Teichmuller Theory and the role of geometric base points in the theory of tempered fundamental groups and in \iut) as shown in \Cref{thdef:arith-loops} given below.}

\bthmdef\label{thdef:arith-loops}
Let $L$ be a number field and let $p$ be a prime number. Let $\O_L$ be the ring of integers of $L$. Consider $\O_L$ equipped with its standard Banach ring structure. Then one has an analytic space $\Spec(\O_L)^{an}$ associated to $\O_L$. Let $K\supset L\supset \O_L$ be an algebraically closed, perfectoid field, with residue characterstic $p$, containing $L$. Then the morphism of analytic spaces
\be 
\sM(K)=\Spec(K)^{an}\to \Spec(\O_L)^{an}
\ee
is an \textit{arithmetic loop in $L$ at the prime number $p$.}
\ethmdef

\newcommand{\cpmax}{\C_p^{max}}
\bp The easiest way of understanding how to establish this assertion is to first assume that $K=\cpmax$ where $\cpmax$ is the maximally complete field containing $\Q_p$ isometrically.  Then  by \cite{kaplansky42}, \cite{lampert1986} and \cite{poonen93} one knows that $\cpmax$ is algebraically closed perfectoid, maximally complete, of uncountable transcendence degree over $\C_p$, and of the same cardinality as $\C_p$ and that  $\cpmax$ is described as a power series in fractional powers $p^\Q$ of $p$ and with prime-to-p-roots-of-unity as coefficients i.e. any $x\in \cpmax$ is a series of the form
$$x=\sum_j{\zeta_j}p^{a_j},$$
for a suitable choice of roots of unity $\zeta_j$ of orders prime to $p$ and suitable rationals $a_j$ contained in some well-ordered subset of $\Q$. [\cite{lampert1986} provides an explicit algorithm for finding such an expression for any $x\in\bQ_p$.]

\textit{Thus one may consider $\sM(\cpmax)\to \Spec(\O_L)^{an}$ as an arithmetic loop at some prime $\wp\subset \O_L$ lying over the prime $p$ i.e. $\wp|p$.}

Now to establish the general case i.e. $K\supset \Q_p$ is an arbitrary algebraically closed field containing an isometrically embedded $\Q_p$. Let $k$ be the residue field of $K$ and let $\abs{K}$ be the value group of $K$. Then by \cite{kaplansky42,poonen93} there exists a field $K^{max}$ which is maximally complete, algebraically closed, perfectoid and which contains $K$ (isometrically embedded) and such that $K^{max}$ has residue field $k$ and has the same value group as $K$ and $K^{max}$ can be described as a field of power series in $p^{\abs{K}}$ (Han-Malcev power series) in a manner similar  to that of $\cpmax$ (see \cite{kaplansky42,poonen93} for details).  
\ep

In particular,  one sees that if $\arith{L}_{\by}$ is an arithmeticoid of $L$ with the arithmetic ring $R_{\by}$, then the ring homomorphism $\O_L\into L\into R_{\by}$ given by the data of $\arithl_\by$ should be considered as providing an arithmetic  loops at all primes of $\O_L$:
\bcor\label{cor:arith-loops-at-primes}
Let $\arith{L}_{\by}$  be an arithmeticoid of $L$ and let $X/\O_L$ be a quasi-projective and generically smooth scheme over the ring of integers of a number field $L$ (assumed to have no real embeddings). Then the choice of any geometric base-point (i.e. a geometric $K_v$-valued base-point for each analytic space $\xan_{K_v}$) of $\holt{X/L}{\by}$ provides an arithmetic loop in each of the analytic spaces $\xan_{K_v}$ at all primes $v\in\vlnon$.
\ecor

\subsection{Arithmetic Knots}
\numberwithin{equation}{subsection}
Now let me explain the relation between my work and the analogy between knots and primes proposed  in \cite{morishita2002,mazur2012}. Let $L$ be a number field as above. Let $v\in\vlnon$ be a prime of $L$ lying over  a rational prime $p$. Let $\F_v$ be the residue field of $\O_{L_v}$.  Let $K_v\supset L_v\supset  L$ be an algebraically closed, perfectoid field with an embedding of $L$ and with residue field $\overline{\F}_v$ which is an algebraic closure of $\F_v$ and $L_v$ is the completion of $L$ at $v$. Then one has the diagram of ring homomorphisms
\be 
\begin{tikzcd}
	\O_{L} \ar[r, hook]\ar[d, two heads] & \O_{K_v}\ar[d,two heads] \ar[r,hook] & K_v\\
	\F_v \ar[r,hook] & \overline{\F}_v & & 
\end{tikzcd}
\ee
The first column of the diagram determines an embedding  $\Spec(\F_v)\into \Spec(\O_L)$, which is a \textit{knot associated to the prime $v\in\vlnon$} in the terminology of \cite{morishita2002} and \cite{mazur2012}. The datum of the algebraically closed (perfectoid field) $K_v$ provides a geometric point for computing the fundamental group  of this knot. This fundamental group  is identified with a certain quotient of the absolute Galois group $G_{L;K_v}$ of $L$ computed using the algebraic closure of $L$ in $K_v$ (the precise definition of these groups can be gleaned from \cite{morishita2002}, \cite{mazur2012}). \textit{My remark here is that there are in fact many such knots and relative to each other they are generally knotted!}  My observation is as follows:
\bthm\label{th:distinct-knots} 
Let $\arithl_{\by_1}$ and $\arithl_{\by_2}$ be two arithmeticoids of $L$. Let $v\in\vlnon$ be a prime of $L$.
Then 
\benumlab
\item These two arithmeticoids provide, for each prime $v\in\vlnon$ the data $(L,v\in\vlnon,K_{y_{v,1}}\supset L_v\supset L)$ and $(L,v\in\vlnon,K_{Y_{v,2}}\supset L_v\supset L)$ which give rise to two diagrams as above and hence to two such knots at each prime $v\in\vlnon$.
\item In general, and relative to each other, the  knots at $v$ are not necessarily unknotted--i.e. there may not exist any natural isomorphism between the local data $(L,v\in\vlnon,K_{y_{v,1}}\supset L_v\supset L)$ and $(L,v\in\vlnon,K_{Y_{v,2}}\supset L_v\supset L)$ which takes the first knot to the second knot. 
\item Suppose $C/\arithl_{\by_1}$ and $C/\arithl_{\by_2}$ are two holomorphoids of an elliptic curve with split multiplicative reduction at some prime $v\in\vlnon$ with $v|p$. Then the $\mathscr{L}$-invariant of $C$ at $v$ computed with respect to the two arithmeticoids $\arithl_{\by_1}, \arithl_{\by_2}$ can be compared by choosing their respective Teichmuller lifts to the ring $W(\O_{\cpt})\tensor_{\Z_p}\O_{L_v}\subset B_{\cpt,L_v}$ in the sense of \constrtwo{Proposition }{pr:teichmuller-lift-A} and, in general, such lifts need not necessarily agree.
\item In particular,  $\yadl,\yadl^{max}$ and $\yadl^\R$ are  parameter spaces of (non-trivial) arithmetic knots at all primes $v\in\vlnon$.
\eenum
\ethm
\bp 
The first two  assertions are immediate consequences of the theory developed in this paper! The third assertion is proved using the method of Teichmuller liftings of
\constrtwo{Proposition }{pr:teichmuller-lift-A}. Write  $B=B_{\cpt,L_v}$, $\by_{1,v}=y_1\in\sY_{\cpt,L_v}$ and $\by_{2,v}=y_2\in\sY_{\cpt,L_v}$. Let $\fm_j\subset B$ with $(j=1,2)$ be the maximal ideals of $B$ corresponding to $y_j$. Let $\mathscr{L}_j\subset L_v\subset K_{v,j}$ be the $\mathscr{L}$-invariant at $v$ of $C/\arithl_{\by_j}$ for $j=1,2$. Let $[x_j]\in W(\O_{\cpt})\tensor_{\Z_p}\O_{L_v}$ be a Teichmuller lift of $\mathscr{L}_j$  for $j=1,2$. In general, the valuations  $\abs{\mathscr{L}_j}_{K_{v,j}}<1$ (and hence $x_j\in \fm_j$) and $[x_j]$ will have norms equal to $\abs{\mathscr{L}_j}_{K_{v,j}}<1$ for any of the Frechet norms on $W(\O_{\cpt})\tensor_{\Z_p}\O_{L_v}$  (\cite[Chapitre 2]{fargues-fontaine}). But $\abs{\mathscr{L}_j}_{K_{v,j}}$ depend on $\abs{p}_{K_{v,j}}$ which need not always agree for arbitrary $y_{v,j}\in\sY_{\cpt,L_v}$. Hence  $[x_1]\neq[x_2]$ (in general). This proves {\bf(3)}.

In the analogy between knot theory and arithmetic espoused in \cite{morishita2002} (and in \cite{mazur2012}), the absolute Galois groups $G_{L_v;K_{y_{v,1}}}$ (resp. $G_{L_v;K_{y_{v,2}}}$) correspond to the `peripheral groups' of the two knots \cite[1.6]{morishita2002}; and the representation of absolute group $G_{L}$ on the $p_v$-adic Tate-module provide  representations of the arithmetic knot groups and the $\mathscr{L}$-invariant appears as an invariant similar to an invariant attached to a topological knot (see \cite[8.3 and 8.6]{morishita2002}). By {\bf(3)}, these two arithmetic knots have dependency on the respective arithmeticoids which can be measured using the lifts of $\mathscr{L}$-invariants (in the sense of {\bf(3)}) and this means that these two arithmetic knots must be considered as topologically distinct! This completes the proof.
\ep

\brem 
In the analogy between knots and primes proposed by \cite{morishita2002}, \cite{mazur2012}, the fundamental group   $\pi_1^{et}(\Spec(\O_L)-\{v\})$ (computed using the geometric base-point given by $\O_L\into K_v$), is considered to be the analog for the arithmetic prime knot at $v$,  of the fundamental group of the complement of a classical knot (i.e. of the classical knot-group). Thus one may view \Cref{th:distinct-knots} as asserting that an arithmetic knot given by a prime is not distinguished by the fundamental group of the arithmetic knot-complement. This is analogous to the fact that in knot theory the trefoil knot and its mirror image have isomorphic knot groups, but are distinct knots.
\erem

\section{Appendix: The proof of geometric Szpiro inequality revisited}\label{se:appendix}
\subparat{The proofs of geometric Szpiro inequality} In \cite{bogomolov00}, the authors gave a proof of geometric Szpiro inequality for elliptic fibrations over a genus zero curve and in \cite{zhang01} this was generalized to elliptic fibrations over base curves of genus $g\geq 0$. In Spring 2018, I was on a sabbatical to RIMS, Kyoto with Shinichi Mochizuki and he gave me a few private lectures on his approach ( \cite{mochizuki-bogomolov}) to these works. As Mochizuki informs us in \cite{mochizuki-bogomolov}, Ivan Fesenko had pointed out to him the existence of \cite{bogomolov00}, \cite{zhang01}.

Note that \cite{bogomolov00}, \cite{zhang01} and \cite{mochizuki-bogomolov} deal with  the same theorem (namely the geometric Szpiro inequality) from three different perspectives. 

Here I will record a fourth distinct approach to understanding these three approaches to geometric Szpiro inequality which I discovered in the course of my work on arithmetic Teichmuller spaces. Notably the perspective presented below exhibits all the key features of \iut, and also the perspective of this paper and a reading of this is essential for anyone who wishes to understand \iut, my work and the relationship between these two theories. In particular,  I provide a formulation (and a proof) of \cite[Corollary 3.12]{mochizuki-iut3} in the context of the geometric Szpiro inequality (this is not formulated in \cite{mochizuki-bogomolov}).

\textcolor{red}{This section should be read in tandem with \cite{zhang01}, \cite{bogomolov00}.} (\cite{mochizuki-bogomolov} largely follows these references with added commentary).

\subparat{The setup} Let $C$ be a connected, smooth, projective curve over $\C$ of genus $g\geq 0$. Let $f:X\to C$ 
be a proper, generically smooth morphism with geometric generic fiber of genus one.  This fibration will be fixed for the remainder of this paper. Additionally I will assume that this is not an isotrivial fibration.

\subparat{The universal cover of $\slr$ and its properties} Let $\slr$ be the topological group of $2\times2$ matrices with real entries and of determinant one. Let $\tslr$ be the universal cover of $\slr$. Then $\tslr$ is described as a central extension (\cite[Lemma 3.5]{zhang01})
\be\label{eq:slr-cover}
1\to \langle z^2\rangle \to\tslr \to \slr \to 1
\ee
where $z\in\tslr$ is a certain element, which generates the center of $\tslr$ and which is described explicitly in \cite{zhang01} and so I will not recall that description here.

Let $$\tslz\subset \tslr$$ be the inverse image of $\slz\subset\slr$. Let $$\tslq\subset \tslr$$ be the inverse image of $\slq\subset \slr$.

\brem 
Note that by classical results or the explicit description in \cite{bogomolov00} or \cite{zhang01}, $\tslz$ is isomorphic to the braid group on three letters (see \cite[Theorem 16]{rawnsley12}).
\erem

\subparat{Classical Teichmuller space in genus one} Let $\fH$ be the complex upper half plane. Then one has a natural mapping 
$$\slr\to \fH$$
given by $$\begin{bmatrix} a & b \\ c & d \end{bmatrix}\mapsto \frac{a\cdot i+b}{c\cdot i +d}.
$$

If $\tau\in \fH$, I will write $E_\tau$ for the elliptic curve given by the lattice $[1,\tau]\subset \C$.

Recall from \cite{imayoshi-book} that $\fH$ is the Teichmuller space of genus one Riemann surfaces and the usual modular group $\slz\subset \slr$ can be identified with the Teichmuller modular group (in genus one).

\subparat{Global (geometric) Fargues-Fontaine curves $\syi$, $\sxi$, and the global Frobenius morphism $\vphi_\infty$}\label{app:global-frobenius} Let me now illustrate how Fargues-Fontaine curve type objects with a similar geometry appears in the context of \cite{zhang01}, \cite{bogomolov00}. This is an important point in understanding my approach and Mochizuki's approach to global Diophantine problems via our respective theories.
\bdefn\label{def:global-frob-classical-case} 
Let $\syi=\tslr$ and $\sxi=\slr$. I will refer to $\syi$ as the \textit{incomplete global Fargues-Fontaine curve} and $\sxi=\slr$ as the \textit{complete global Fargues-Fontaine curve.} I will write  
\be\label{eq:global-frob-classical-case}\vphi_\infty=z^2\in\tslr\ee for the generator of the central extension \eqref{eq:slr-cover} and refer to $\vphi_\infty$ as the \textit{global Frobenius morphism of $\syi$}. In particular,  with these definitions, and because of \eqref{eq:slr-cover}, one has
$$\syi/\vphi_\infty^\Z\mapright{\isom} \sxi$$
is a quotient by the global Frobenius morphism.
\edefn

\brem\
\benumlab
\item  Let me spell out the analogy  to the situation of usual Fargues-Fontaine curves where one has a quotient of adic spaces
$$\syfe/\vphi^\Z\mapright{\isom}\sxfe$$
with $\vphi$ the Frobenius morphism of $\syfe$. 
\item Note however that the usual Fargues-Fontaine curve $\syfe,\sxfe$ are local objects from the the point of view of theory of the present series of papers. 
\item On the other hand,  $\syi,\sxi$ are global objects in the context of \cite{zhang01}, \cite{bogomolov00}.
\item Similarly  $\yadl$ and the quotient (topological stack) $[\yadl/L^*]$ are global analogs of $\syi,\sxi$ in the context of the present series of papers (and hence in the context of \iut).
\item Notably $\yadl$ equipped with the action of $L^*$  or its global Frobenius  morphism $\boldsymbol{\varphi}$ given by \Cref{th:galois-action-on-adelic-ff} and \Cref{cor:global-frobenius} are the global arithmetic analogs of $\syi$ and its Frobenius morphism $\varphi_\infty$.
\eenum
\erem

\subparat{Existence of $\flog$-links} In \iut\ Mochizuki works with $\flog$-links. Let me discuss how this appears in the context of \cite{zhang01}, \cite{bogomolov00}. Let $g_\tau\in\slr$ be a lift of $\tau\in \fH$ to $\slr$. Let $\tilde{g}_\tau\in\tslr$ be the lift of $g_\tau$ to $\tslr$.
\bdefn\label{de:geom-log-links}
Let $\tilde{g}_{\tau,0}\in\tslr$ be a lift of $\tau\in\fH$ to $\tslr$. If $\tilde{g}'_{\tau,0}$ is another lift of $\tau$ then I will say that $\tilde{g}_{\tau,0}$ and $\tilde{g}_{\tau,0}'$ are $\flog$-linked if $$\tilde{g}'_{\tau,0}=\tilde{g}_{\tau,0}\cdot \vphi_\infty.$$
I will call $$\left\{\tilde g_{\tau,n}=\tilde{g}_{\tau,0}\cdot \vphi_\infty^n: n\in \Z \right\}$$
a chain of $\flog$-links.
\edefn

\brem
Note that a chain of $\flog$-links is the fiber over $\syi\to\sxi$ over the common image of $\tilde{g}_{\tau,n}$ in $\sxi=\slr$.
\erem

\subparat{The Schottky parameter}\label{ss:schottky-parameterization} Let me begin with the following 
\blem\label{le:link-lemma}
Let $\tau\in\fH$. For any  $\alpha\in\Q^{>0}$,  the elliptic curves $E_1=\C/\Lambda_1$ and $E_{\tau,\alpha}=\C/\Lambda_{\tau,\alpha}$ given by the lattices $\Lambda_{\tau,1}=[1,\tau]\subset \C$ and $\Lambda_{\tau,\alpha}=[1,\alpha\cdot\tau]$ are $\C$-isogenous.
\elem
\bp 
It is sufficient to note that 
$\begin{pmatrix}
	1 & 0\\
	0 & \alpha
\end{pmatrix}\in\glqp$ and
$$\begin{pmatrix}
1 & 0\\
0 & \alpha
\end{pmatrix} \cdot \begin{pmatrix}
1 \\ \tau
\end{pmatrix} =\begin{pmatrix}
1 \\ \alpha\cdot\tau
\end{pmatrix}
$$ takes $\Lambda_1$ to $\Lambda_\alpha$. So the assertion is clear.
\ep

\bdefn 
For an elliptic curve $E$ given by a lattice $[1,\tau]\subset\C$ let $q_E=e^{2\pi \cdot \sqrt{-1}\tau}$ be its Schottky parameter. Then one has the Schottky uniformization
$$E\isom \C^*/q_E^\Z.$$
\edefn
\bp 
This is standard.
\ep 

\blem\label{le:link-schottky-scaling} 
For the elliptic curves $E_1$ and $E_{\tau,\alpha}$ one has the following relationship between their Schottky parameters:
$$q_{E_{\tau,\alpha}}=q_{E_1}^\alpha.$$
\elem
\bp 
This is an elementary calculation:
$$q_{E_{\tau,\alpha}}=e^{2\pi \cdot \sqrt{-1}(\alpha\cdot\tau)}=e^{2\pi \cdot \sqrt{-1}\cdot \alpha\cdot\tau}=q_{E_1}^\alpha.$$
\ep

\subparat{Existence of $\Theta_{gau}$-links}\label{pa:schottky-param}
 Now I want to discuss how Mochizuki's $\Theta_{gau}$-link appears in the geometric case of \cite{zhang01}, \cite{bogomolov00}. Let $\ell\geq 5$ be a prime number, let $\ells=\frac{\ell-1}{2}$. Let $\tau\in \fH$ and  write $q=e^{2\pi\sqrt{-1}\tau}$. 
From now on, for $j=1,2,\ldots,\ells$, for simplicity of notation, I will write $$E_{\tau,j}=E_{\frac{\tau}{2\ell},j^2}.$$

Note that the lattice $[1,\frac{\tau\cdot j^2}{2\ell}]$ for $E_{\tau,j}$ can be obtained from the lattice $[1,\frac{\tau}{2\ell}]$ of $E_{\tau,1}$ by using the matrix
$$\begin{pmatrix}
j & 0\\
0 & \frac{1}{j}
\end{pmatrix}\in \slq$$
which takes operates on $\fH$ by $$\tau\mapsto \begin{pmatrix}
j & 0\\
0 & \frac{1}{j}
\end{pmatrix}\cdot \tau=\frac{j\cdot \tau}{j^{-1}}={j^2\cdot \tau},$$
and hence $\frac{\tau}{2\ell}\mapsto \frac{j^2\cdot \tau}{2\ell}$ under its action.

\brem 
The reason I want to consider $\slq$ in place of $\glqp$ is that $\slq\subset \slr$. Working with $\glqp$ would require replacing $\slr$ by ${\rm GL}_2^+(\R)$ in \cite{zhang01}, \cite{bogomolov00}--this can be done, but presently I have not written down all the changes this necessitates in \cite{zhang01}, \cite{bogomolov00}.
\erem

\blem 
Let $$\alpha_j=\begin{pmatrix}
j & 0\\
0 & \frac{1}{j}
\end{pmatrix}\text{ for }j=1,2,\ldots,\ells.$$
Then the $\ells$-tuple $(\alpha_1,\ldots,\alpha_\ells)\in\slq^\ells$ operates on $\fH^\ells$, by the usual Mobius action on each of the factors, and under this action the diagonal tuple $$(\tau,\tau,\ldots,\tau)\in\fH^\ells$$ is mapped to 
$$(1^2\cdot\tau,2^2\cdot\tau,\ldots,\ells^2\cdot\tau).$$
\elem
\bp 
This is completely clear.
\ep

\bdefn 
Let $\ell\geq5$ be a prime number.
A $\Theta_{gau}$-link is an $\ells$-tuple $$(\tilde{g}_1,\ldots,\tilde{g}_{\ells})\in\tslql\subset\syi^\ells$$ lying over a tuple of elliptic curves $(E_{\tau,1},E_{\tau,2},\ldots,E_{\tau,\ells})$ for some $\tau\in\fH$ (with $E_{\tau,j}$ given by the above convention  and $j=1,2,\ldots,\ells$). The ``link'' here is the rule
$$E_\tau\mapsto (E_{\tau,1},E_{\tau,2},\ldots,E_{\tau,\ells}).$$
\edefn
\brem 
One may think of $E_\tau\mapsto (E_{\tau,1},E_{\tau,2},\ldots,E_{\tau,\ells})$ as a correspondence i.e. a one-to-many function on the moduli of elliptic curves arising from the correspondence on $\fH$  given by the correspondence:
$$\tau\mapsto (1^2\cdot \tau,2^2\cdot \tau, \ldots,\ells^2\cdot \tau)\in\fH^\ells.$$
\erem

\bthm\label{th:theta-tuples}
Let $(\tilde{g}_1,\ldots,\tilde{g}_\ells)\in\syi^\ells$ be a $\Theta_{gau}$-link. Let $E_{\tau,1}, E_{\tau,2}, \ldots,E_{\tau,\ells}$ be the corresponding elliptic curves. Then 
\benumlab
\item The elliptic curves $E_{\tau,1}, E_{\tau,2} ,\ldots , E_{\tau,\ells}$ correspond to distinct points of the genus Teichmuller space $T_1=\fH$.
\item The elliptic curves $$E_{\tau,1}, E_{\tau,2} ,\cdots , E_{\tau,\ells}$$ are all isogenous to $E_\tau$.
\item One has the  tuple of the corresponding Schottky parameters
$$(q_{\tau,1}, q_{\tau,2},\ldots,q_{\tau,\ells}) \in(\C^*)^\ells$$ corresponding to $E_{\tau,1}, E_{\tau,2} ,\ldots , E_{\tau,\ells}$.
\item Let   $P_1$ correspond to the $\ell$-torsion point on the double cover $E_{\tau/2}$ of $E_\tau$ given by $\tau'=\frac{\tau}{\ell}$ and choose a $\Theta$-function as in \cite{mochizuki-theta}. Then
the above tuple of Schottky parameters
$$(q_{\tau}^{1^2/2\ell},q_\tau^{2^2/2\ell},\ldots,q_\tau^{\ells^2/2\ell})\in (\C^{*})^\ells$$
corresponds to the set of values of the chosen $\Theta$-function at the $\ell$-torsion points $$\left\{P_j=j\cdot P_1\in E_{\tau/2}[\ell]: j=1,2,\ldots,\ells \right\}.$$
\eenum
\ethm
\bp 
The proof is clear from the previous lemmas.
\ep

\brem 
Since $\begin{pmatrix}
j & 0\\
0 & \frac{1}{j}
\end{pmatrix}\in\slq\subset \glqp$, this theorem makes it clear that one can think of a $\Theta_{gau}$-Link as arising  from natural correspondences on the upper half plane in a manner similar to Hecke correspondences on classical modular curves. By \Cref{le:link-schottky-scaling}, these correspondences provide a (non-trivial) scaling $$q_\tau\mapsto q_\tau^{j^2}$$ of the Schottky parameters. In the $p$-adic case, \cite{joshi-teich-estimates} takes a similar approach via correspondences on suitable Fargues-Fontaine curves for the construction of $\Theta_{gau}$-links in the context of \iut.
\erem
\newcommand{\bTheta}{\boldsymbol{\Theta}}

\subparat{The $\flog$-$\Theta$-lattice}
Let me note that $\slq,\tslq$ are both countable sets and hence so are $\slq^\ells$ and $\tslql$. Choose an enumerating function to enumerate the elements of $$\slq^\ells=\{\theta_{n}:n\in\Z\}.$$
Then one may  enumerate the set $\tslql$ as 
$$\tslql=\{\theta_{n,m}:(n,m)\in\Z\times \Z\},$$
with the canonical surjection $$\tslql\to \slq^\ells$$ being given by 
$$\theta_{n,m}\mapsto \theta_n.$$
Following \cite{mochizuki-iut1,mochizuki-iut3}, one could represent the  data $\tslql=\{\theta_{n,m}\}$   as a ``lattice'' in the plane with the vertical direction corresponding to the fiber $\{\theta_{n,m}:m\in\Z \}$ of $\tslql\to\slq^\ells$ over $\theta_n\in\slq^\ells$ i.e. over a $\Theta_{gau}$-link. 
In this picture the fiber over $\theta_n$ consists of $\flog$-links i.e. 
while the horizontal direction being represented as $\Theta_{gau}$-links. In \cite{mochizuki-iut3} this sort of presentation of $\Theta_{gau}$-links and $\flog$-links is called the \textit{$\log$-$\Theta$-lattice}:
$$
\xymatrix{
	&  &  & \\ 	
	\cdots \ar@{-}[r] &\theta_{n-1,m+1}\ar@{-}[r]\ar@{-}[d]_{\flog=\vphi_\infty}\ar@{-}[u]^{\flog=\vphi_\infty} & \theta_{n,m+1}\ar@{-}[d] \ar@{-}[u]_{\flog=\vphi_\infty}\ar@{-}[r] & \cdots\\ 
	\cdots \ar@{-}[r] & \theta_{n-1,m} \ar@{-}[r]\ar@{-}[d]_{\flog=\vphi_\infty} &  \theta_{n,m}\ar@{-}[u]_{\flog=\vphi_\infty} \ar@{-}[r]\ar@{-}[d]^{\flog=\vphi_\infty} & \cdots,\\
	&  &  & 
}$$ 
\brem
This should be compared with Mochizuki's discussion of the $\log$-$\Theta$-lattice in \cite[\ssep I, Page 405]{mochizuki-iut3} and its construction  in \cite[\ssep 1, Page 427]{mochizuki-iut3}.
\erem
\subparat{Hodge structures from Schottky parametrization}\label{ss:hodge-str-schottky}
\textcolor{red}{This is an elaboration of \cite[\ssep 22]{joshi-anabelomorphy}.}
Let $E_\tau$ be an elliptic curve with period lattice $[1,\tau]\in\fH$. Then the Schottky parameterization \ssep\ref{pa:schottky-param} provides a mixed Tate Hodge structure 
$$H_{Schottky}(E_\tau)\in \Ext^1_{\Z-MHS}(\Z(0),\Z(1))\isom \C^*$$ given via \cite{deligne-local} by $$H_{Schottky}(E_\tau)=q_\tau=e^{2\pi \sqrt{-1}\tau}\in\C^*.$$  This will be useful in understanding how one may carry out constructions in cohomology theory. In \iut\ cohomology theory means Galois cohomology theory $H^1(G_{L_v},\Q_v(1))$ and one has $$H^1(G_{L_v},\Q_v(1))=\Ext^1_{G_{L_v}}(\Q_v(0),\Q_v(1))$$ for each $v\in\V_L$, so the analogy I make here should not be lost on the readers.

\brem Note that one way to understand the Schottky parametrization  is to recognize that it asserts that the elliptic curve $E_\tau$ is  ordinary at an archimedean prime (and this ordinarity is witnessed by the Schottky Hodge structure $H_{Schottky}(E_\tau)$). \textit{However note that the usual Hodge structure provided by $H^1(E_\tau,\Z)$ is pure of weight one.} Tate's parameterization in the $p$-adic context is the $p$-adic version of Schottky's paramterization.
\erem

Thus one has 
\blem 
Any tuple  $(q_{\tau,1}, q_{\tau,2},\ldots,q_{\tau,\ells}) \in(\C^*)^\ells$ in \Cref{th:theta-tuples} provides a tuple of $\Z$-mixed Hodge structures  $$\left(H_{Schottky}(E_{\tau,1}), H_{Schottky}(E_{\tau,2}) ,\ldots , H_{Schottky}(E_{\tau,\ells})\right) \in \Ext^1_{\Z-MHS}(\Z(0),\Z(1))^\ells.$$
\elem

\subparat{A prelude to the definition of the $\Theta$-values locus} 
Let $$\bTheta_{classical,\tau}=\left\{ (\tilde{g}_1,\ldots,\tilde{g}_\ells)\in\tslql\subset  \syi^\ells: (\tilde{g}_1,\ldots,\tilde{g}_\ells) \text{ provides } E_{\tau} \right\}.$$

Let $S\subset C(\C)$ be the finite set of points over which the fibers of $X\to C$ are singular. 

Let $s\in S$ be a singular point of $f:X\to C$. One would like to understand $$\bTheta_{classical,\tau_s}$$ for $\tau_s\in \fH$ in a sufficiently close i.e. in a small neighborhood of the singular point $s\in S$ of $f:X\to C$ so that the local description of \cite{zhang01} remains valid (for each singular point $s\in S$). 

\subparat{Mochizuki style Theta-values locus for \cite{zhang01},
	\cite{bogomolov00}} So to construct a $\Theta$-values locus in the style of \cite{mochizuki-iut3} one does the following. Let $\syis$ be a copy of $\syi$ indexed by an $s\in S$.

Let \be\label{eq:theta-value-locus}\bTheta_{classical}=\left\{ (\tilde{g}_{s,1},\ldots,\tilde{g}_{s,\ells})_{s\in S}\in   \prod_{s\in S} {\syils} \text{ and } (\tilde{g}_{s,1},\ldots,\tilde{g}_{s,\ells})_{s} \in \bTheta_{classical,\tau_s} \text{ for all } s\in S  \right\}\ee
Each $(\tilde{g}_{s,1},\ldots,\tilde{g}_{s,\ells})_{s\in S}\in\bTheta_{classical}$ is 
the analog of the  \textit{$\Theta$-pilot object}, in the sense of \cite{mochizuki-iut3}, in the context of \cite{zhang01,bogomolov00}. 

One may also define the theta-values set $$\bTheta_{classical}^{val}\subset \Ext^1_{\Z-MHS}(\Z(0),\Z(1))^\ells,$$ in the style of  \iut, for the context of \cite{zhang01, bogomolov00} by taking the tuples of Schottky parameters provided by each tuple $(\tilde{g}_{s,1},\ldots,\tilde{g}_{s,\ells})_{s\in S}\in\bTheta_{classical}$.

\subparat{Stablized height functions $h$ and $H$}\label{ss:height-stabilzation}\
Let $$h:\tslr \to \R$$
be the function defined in \cite[Proof of Theorem 3.3]{zhang01} (also defined in \cite{bogomolov00}) and denoted in \cite{zhang01} by $$\ell:\tslr \to \R.$$ For $g\in\tslr$, write
$$H=e^{h}:\tslr\to \R.$$
One should think of $h$ as the \textit{logarithmic Frobenius-stabilized height function} and $H$ as the \textit{Frobenius-stabilized height function} or more simply as the \textit{stablized logarithmic height function} and \textit{stablized height function} respectively (also see \ssep\ref{ss:height-stabilzation}).

\subparat{Properties of the height function} The central point in the proofs of \cite{bogomolov00,zhang01} is the following property of the height function introduced in \ssep\ref{ss:height-stabilzation}. By \cite{zhang01} one has  for any $\tilde{g}_{1},\tilde{g}_{2}\in\tslr$:
\be\label{eq:bogomolov-zhang} h(\tilde{g}_{1}\cdot\tilde{g}_2)\leq h(\tilde{g}_{1})+h(\tilde{g}_{2}).\ee 
Given a $\Theta$-link $(\tilde{g}_{s,1},\ldots,\tilde{g}_{s,\ells})_{s\in S}\in \bTheta_{classical}$, the quantity  of interest to us is:
$$\sum_{s\in S}h(\tilde{g}_{s,1}\cdot\tilde{g}_{s,2} \cdots \tilde{g}_{s,\ells} ).$$

Hence by \eqref{eq:bogomolov-zhang} one has $$\sum_{s\in S} h(\tilde{g}_{s,1}\cdot\tilde{g}_{s,2} \cdots \tilde{g}_{s,\ells} )\leq \sum_{s\in S}\sum_{j=1}^\ells h(\tilde{g}_j).$$

In particular one can also consider the average
$$\frac{1}{\ells}\sum_{s\in S} h(\tilde{g}_{s,1}\cdot\tilde{g}_{s,2} \cdots \tilde{g}_{s,\ells} )\leq \frac{1}{\ells}\sum_{s\in S}\sum_{j=1}^\ells h(\tilde{g}_j).$$

\brem 
The proof of \cite{zhang01}, \cite{bogomolov00} (and \cite{mochizuki-bogomolov}) can be understood as a sort of averaging over $\tslz$. In what I do here,  this averaging is being replaced here by averaging over $\tslql$. The innovation presented here is the idea of averaging over $\tslql$ which is not considered in \cite{mochizuki-bogomolov}--but corresponds to working simultaneously with Mochizuki's $\Theta_{gau}$-Links,  and $\flog$-links (\cite{mochizuki-iut3}).
\erem

\subparat{The height function $h$ and the global Frobenius morphism (aka $\flog$-Link)}\label{ss:geom-case-height-frobenius}
Let me explain the reason why $h$ is referred to as a stabilized height function in \cref{ss:height-stabilzation}. Let $\tilde{g}'_{\tau,0},\tilde{g}_{\tau,0}\in\tslr$ and assume that $$\tilde{g}'_{\tau,0}=\tilde{g}_{\tau,0}\cdot \vphi_\infty^m$$
for $m\in\Z$. By \Cref{de:geom-log-links}, this equation means that $\tilde{g}'_{\tau,0},\tilde{g}_{\tau,0}$ differ by an iterate of the $\flog$-Link or equivalently differ by the  $m^{th}$ iterate of the global Frobenius $\vphi_\infty$ for some $m\in\Z$. 
Then one has the following key relationship established in \cite[Proof of the Geometric Szpiro Inequality]{zhang01} between the heights  $h(\tilde{g}'_{\tau,0})$ and $h(\tilde{g}_{\tau,0})$ i.e. between the height of $\tilde{g}_{\tau,0}$ and the height of a  log-linked element $\tilde{g}'_{\tau,0}$:
\be h(\tilde{g}'_{\tau,0}) = h(\tilde{g}_{\tau,0}\cdot\vphi_\infty^m) \leq h(\tilde{g}_{\tau,0}) + \pi\cdot m\in\R.\ee
(where $\pi\in\R$ is the familiar real number). 

This property plays a central role in the proof of \cite[Proof of the Geometric Szpiro Inequality]{zhang01} and Mochizuki's discussion of this point  and its relationship to his proof can be found in \cite{mochizuki-bogomolov}.

\subparat{The definition of $\abs{\bTheta_{classical}(D)}$} More generally, keeping in mind \cite{mochizuki-iut4}, let $S\subset D\subset C(\C)$ be a subset containing $S$ such that the local theory of \cite{zhang01} can be applied (for example let $D$ be the disjoint union of small disk $D_s\subset C(\C)$ around each $s\in S$). Then it is possible to define 
$\bTheta_{classical}(D)\supset \bTheta_{classical}$ as follows:
\bdefn\label{def:theta-abs-val-def} 
For a compact domain $S\subset D\subset C(\C)$ chosen as above, let
$$\abs{\bTheta_{classical}(D)}=\sup\left\{\sum_{s\in S} \sum_{j=1}^\ells h(\tilde{g}_{s,j}): (\tilde{g}_{s,1},\tilde{g}_2 \cdots, \tilde{g}_{s,\ells})_{s\in S}\in \bTheta_{classical}(D) \right\}.$$
\edefn

The following lemma will be useful:
\blem\label{le:sum-product-h-prop} 
$$\abs{\bTheta_{classical}(D)}\geq \sup\left\{\sum_{s\in S}h(\tilde{g}_{s,1}\cdot\tilde{g}_2 \cdots \tilde{g}_{s,\ells}): (\tilde{g}_{s,1},\tilde{g}_2 \cdots, \tilde{g}_{s,\ells})_{s\in S}\in \bTheta_{classical}(D) \right\}$$
\elem
\bp 
This is immediate from the above discussed property \eqref{eq:bogomolov-zhang} of $h$ and the definition of $\abs{\bTheta_{classical}(D)}$.
\ep
\brem 
Working with a compact domain $D$ is analogous to working with compactly bounded domains (in the sense of \cite{mochizuki-general-pos}) in the proof of the main theorem of \cite{mochizuki-iut4}.
\erem
\subparat{The monodromy representation and its reduction modulo $\ell$} Let $\pi_1^{top}(C(\C)-S,*)$ be the topological fundamental group of $C(\C)-S$ for some choice of basepoint $*\in C(\C)-S$. This is a group with $2g+\abs{S}$ generators $$a_1,\ldots,b_1,\ldots,a_g,b_g,\{\gamma_s\}_{s\in S}$$ and one relation
$$\prod_{j=1}^g[a_j,b_j]\prod_{s\in S}\gamma_s=1$$
 between these generators.

Let $t\in C(\C)-S$ be a general point and $X_t$ be the fiber of $X\to C(\C)$ over this point and consider the monodromy representation (on the singular homology)
$$\rho:\pi_1^{top}(C(\C)-S,*)\to {Aut}(H_1(X_t,\Z)),$$
and also its ``reduction modulo $\ell$'':
$$\rho_\ell:\pi_1^{top}(C(\C)-S,*)\to {Aut}(H_1(X_t,\Z/\ell)).$$

One can assume without any loss of generality that $\rho_\ell$ is irreducible otherwise by a classical argument due to Faltings, Szpiro's inequality is trivial.

At any point $s\in S$ one has the monodromy around $s$ which is the image of $\rho(\gamma_s)$ and also $\rho_\ell(\gamma_s)$. Let $\tilde{\gamma}^{(j)}_s\in\tslz$ be an arbitrary lift of $\rho_\ell(\gamma_s^{j^2})$ for $j=1,\ldots,\ells$.

\subparat{Mochizuki's Corollary 3.12 in the geometric case} \cite[Corollary 3.12]{mochizuki-iut3} now has a simple formulation with a tautological proof:
\bthm\label{th:geometric-case-of-moccor} 
For a compact domain $D$ as above, and any choice of $(\tilde{g}_1,\ldots,\tilde{g}_\ells)\in\bTheta_{classical}(D)$ , one has
$$\abs{\bTheta_{classical}(D)}\geq \sum_{s\in S} \sum_{j=1}^\ells h(\tilde{\gamma}_{s}^{(j)}) \geq h\left(\prod_{s\in S}\tilde{\gamma}^{(1)}_s\cdot\tilde{\gamma}^{(2)}_s\cdots\tilde{\gamma}^{(\ells)}_s\right).$$
\ethm
\bp 
 The lower bound is immediate from \Cref{le:sum-product-h-prop} and the fact that, by \Cref{def:theta-abs-val-def}, $\abs{\bTheta_{classical}(D)}$ is the supremum of all such sums.
\ep
\brem 
One can think of the product $\prod_{s\in S}\left(\tilde{\gamma}^{(1)}_s\cdot\tilde{\gamma}^{(2)}_s\cdots\tilde{\gamma}^{(\ells)}_s\right)$ and more generally the product $\prod_{s\in S}\left(\tilde{\gamma}^{(1)}_s\cdot\tilde{\gamma}^{(2)}_s\cdots\tilde{\gamma}^{(\ell)}_s\right)$  as a ``simulation'' of a ``multiplicative'' one dimensional subspace of $\rho_\ell$. Since $\rho_\ell$ is irreducible there does not exist any natural one dimensional subspace  (globally) in $\rho_\ell$. This multiplicative subspace point of view is due to \cite{mochizuki-iut3}.
\erem 
\subparat{Additional Remarks} At this point I want to make the following important remarks. 
\benumlab
\item Working with $\slq\subset \glqp$ is essentially equivalent to  Mochizuki's idea of working with $\Theta_{gau}$-links and working with $\vphi_\infty$ is equivalent  $\flog$-links (in \cite{joshi-teich-estimates} I prove both these assertions in the $p$-adic context). 
\item Especially working with $\tslq$ allows us to work with $\Theta_{gau}$-links and $\flog$-links simultaneously. In \iutthr, Mochizuki works with $\Theta_{gau}$-Links and $\flog$-Links separately. My construction of Mochizuki's Ansatz in \cite{joshi-teich-estimates} allows me to work with both these devices simultaneously.
\item The global Frobenius morphism \eqref{def:global-frob-classical-case} plays a central role in \cite{zhang01}, \cite{bogomolov00}.
\item The global relation between generators of the topological fundamental group $$\pi_1^{top}(C(\C)-S,*),$$ namely
$$\prod_{j=1}^g[a_i,b_i]\prod_{s\in S}\gamma_s=1,$$
plays the role of the product formula \eqref{eq:prod-formula1} in the theory of arithmeticoids and their height functions.
\item The function $h$, which is denoted by $\ell$ in \cite{zhang01}, plays  the role of stabilized height function in the  proof of the geometric Szpiro Inequality in \cite{zhang01}.
\eenum

\section{Appendix II: Existence of mutations of varieties over number fields}\label{se:appendix2}
Classical Teichmuller Theory may be understood as exhibiting mutations of the archimedean periods of  Riemann surfaces. In this appendix, I want to provide an independent proof of the existence of mutations of $p$-periods of a variety over a number fields. I will provide a concrete proof, using the case of elliptic curves defined over a number fields, as an example, that mutations of $p$-adic periods also exist. This may be viewed as signaling the existence of an arithmetic Teichmuller Theory. 
\subparat{Transcendence of the Tate parameter} The following consequence of the main result of \cite{barre96} will be used in what follows:
\bthm\label{th:trans-tate} 
Assume $C/\C_p$ is be a Tate elliptic curve over $\C_p$ with its $j$-invariant $$j=j_X\in\bQ$$ (so $C$ is definable over a number field) and with Tate parameter $q$. Then $q$ is transcendental over $\Q$.
\ethm
\bp 
This is a consequence of the main theorem of \cite{barre96} which asserts that whenever $C/\C_p$ is defined by an algebraic $q\in \C_p$ with $0<\abs{q}_{\C_p}<1$ then the $j=j(q)$-invariant of the corresponding  Tate curve is transcendental. So if $q$ is algebraic for $C/\C_p$ then one has arrived at a contradiction to the main theorem of \cite{barre96}.
\ep
\subsection{Existence of arithmetic mutations of an elliptic curve over a number field}
Classical Teichmuller Theory may be understood as exhibiting mutations of the periods of Riemann surfaces. 
The transcendence of the Tate parameters (\Cref{th:trans-tate})  provides a way of uncovering the existence of  the mutations of periods in the  $p$-adic case: The Tate parameter (i.e. the $p$-adic periods) of a Tate elliptic curve can be  subjected to non-trivial arithmetic changes or mutations while the isomorphism class of the elliptic curve stays fixed.  One may call this an \textit{arithmetic mutation of an elliptic curve over a number field}.  The term will be defined by the following theorem:
\bthmdef\label{thdef:mutation}
Let $C/L$ be an elliptic curve over a number field $L$. Let $p$ be a rational prime and suppose that $\wp_1,\wp_2,\ldots,\wp_n$ are primes $L$ lying over $p$ and that $C$ has split multiplicative reduction at each of the $\wp_1,\wp_2,\ldots,\wp_n$. Let $q_1,\ldots,q_n\in\C_p$ be the Tate parameters at the primes $\wp_1,\wp_2,\ldots,\wp_n$. Then there exists a $\Q$-automorphism $\sigma:\C_{p}\to \C_p$ which maps $C/L$ to an isomorphic pair $C'/L'$ but  $\sigma$ does not preserve $q_1,\ldots,q_n$. Let $q_1',\ldots,q_n'\in\C_p$ be the Tate parameters at the primes $\sigma(\wp_1),\sigma(\wp_2),\ldots,\sigma(\wp_n)$ of $L'=\sigma(L)$. Then 
$$(C/L,q_1,\ldots,q_n)\mapsto (C'/L',q_1',\ldots,q_n')$$ will be called the \textit{arithmetic mutation of $C/L$ given by $\sigma$}.
\ethmdef
\bp
This proof draws on the techniques of the proof of the existence of algebraically closed, perfectoid fields given in \cite[Theorem 2.16.1]{joshi-teich}. Since $\C_p$ is algebraically closed, $\C_p$ contains an algebraic closure of its prime subfield $\Q$. In particular $\C_p$ contains $L$ and also an algebraic closure of $L$. By \Cref{th:trans-tate}, the Tate parameters $q_1,\ldots,q_n$ at $\wp_1,\wp_2,\ldots,\wp_n$ are all transcendental over $\Q$,  but they need not all be algebraically independent over $\Q$. After relabeling the primes and the Tate parameters, one may assume  that $r$  is largest integer $1\leq r\leq n$ with the property that $q_1,\ldots,q_r$ are algebraically independent over $\Q$.  Then one can  construct a transcendence basis $T$ of $\C_p/\Q$ containing this maximal algebraically independent subset  $q_1,\ldots,q_r$ of $q_1,\ldots,q_n$.  Choose a $\Q$-automorphism of $\C_p$ which takes $\sigma(q_j)=q_j^{-1}$ for $1\leq j\leq r$ and $\sigma$ is arbitrarily chosen at other elements of this basis $T$ (this similar to the proof of \cite[Theorem 2.16.1]{joshi-teich}). 

Let $L'=\sigma(L)$ (so $L\isom L'$ as fields) and let $C'=\sigma(C)$ be the elliptic curve $C$ viewed as an elliptic curve over $L'$ via $L\isom L'$. Of course, one has $C/L\isom C'/L'$.  Let $\abs{-}_v$ be the pull-back of $\abs{-}_{\C_p}$ as in the proof of \cite[Theorem 2.16.1]{joshi-teich}. This means $\abs{x}_v=\abs{\sigma(x)}_K$ for all $x\in K$.

As $\sigma(q_j)=q_j^{-1}$, for $1\leq j\leq r$,
by the construction of $v$, one has $$\abs{q_j}_v=\abs{\sigma(q_j)}_{\C_p}=\abs{q_j^{-1}}_{\C_p}>1.$$ So one sees that no $\O_{\C_p}^*$-multiple of any of the $q_1,\ldots,q_r$ can be the Tate parameter of $C'/L'$ at any of the primes $\sigma(\wp_1), \ldots,\sigma(\wp_n)$ of $L'$. So $C'/L'$ is is equipped with a new set of Tate parameters $q_1',\ldots,q_n'\in \C_p$ with $\abs{q_j'}_v<1$ for all $1\leq j\leq n$. 

Hence no $\O_{\C_p}^*$-multiple of any of the $q_1,\ldots,q_r$  appears in the list, $q_1',\ldots,q_n'\in \C_p$, of Tate parameters of $C'/L'$, as $\abs{q_j}_v>1$ for $1\leq j\leq r$, while $\abs{q_j'}_v<1$ for $1\leq j\leq n$. Thus the tuple of Tate periods is not preserved by $\sigma$. That is the Tate periods have mutated non-trivially under $\sigma$ while the isomorphism class of the curve $C/L$ has remained fixed.
\ep

\brem 
The main theorems of \cite{joshi-teich} assert, in the notation of the above proof, that (1) the analytic spaces $C^{an}/(K,\abs{-}_v)$ and $C^{an}/(K,\abs{-}_K)$ obtained from the two algebraically closed, complete valued fields are not $\Q_p$-isomorphic (both  $(K,\abs{-}_v), (K,\abs{-}_K)$, contain their own private copies of $\Q_p$ i.e. the embeddings of $\Q_p$ in each valued field may be incomparable, but $\Q_p$ is the minimal complete valued field isometrically contained in both over which one can compare the two analytic spaces) and (2) both these analytic spaces have topologically isomorphic tempered fundamental groups. Hence the mutation described in \Cref{thdef:mutation} can be viewed as arising from distinct analytic spaces while preserving fundamental groups--a fact reminiscent of classical Teichmuller Theory.
\erem

\bibliography{../../master/masterofallbibs.bib}
\end{document}